%% file: main.tex
\documentclass[11pt]{article}

\usepackage[margin=1in]{geometry}
\usepackage{amsmath,amssymb,amsthm,amsxtra}
\usepackage{thmtools,thm-restate}
\usepackage{mathrsfs, mathtools}
\usepackage{bm}
\usepackage[shortlabels]{enumitem}
\usepackage{microtype}
\microtypesetup{expansion=false}
\usepackage[T1]{fontenc}
\usepackage{xcolor}

\usepackage{float}
\usepackage[numbers]{natbib}
\usepackage{hyperref}
\usepackage{cleveref}
\usepackage{graphicx}
\usepackage{tikz}
\usepackage{algorithm}
\usepackage{todonotes}
\usepackage{dsfont}
\usepackage[noend]{algpseudocode}
\makeatletter
\renewcommand{\theHALG@line}{\thealgorithm.\arabic{ALG@line}}
\makeatother

\hypersetup{hidelinks,citebordercolor={0 1 0}}
\makeatletter
\AddToHook{cmd/hyper@natlinkstart/before}{%
    \begingroup
    \ifnum\NAT@ctype=\z@
        \hypersetup{pdfborder={0 0 1}}%
    \fi
}
\AddToHook{cmd/hyper@natlinkend/after}{\endgroup}
\makeatother
\numberwithin{equation}{section}

\declaretheorem[name=Theorem,numberwithin=section,
    refname={theorem,theorems},Refname={Theorem,Theorems}]{theorem}

\declaretheorem[name=Definition,sibling=theorem,
    refname={definition,definitions},Refname={Definition,Definitions}]{definition}
\declaretheorem[name=Lemma,sibling=theorem,
    refname={lemma,lemmas},Refname={Lemma,Lemmas}]{lemma}
\declaretheorem[name=Corollary,sibling=theorem,
    refname={corollary,corollaries},Refname={Corollary,Corollaries}]{corollary}

\declaretheorem[name=Observation,sibling=theorem,
    refname={observation,observations},Refname={Observation,Observations}]{observation}

\newcommand{\EE}{\mathbb{E}}

\newcommand{\eps}{\epsilon}

\newcommand{\1}{\mathbf{1}}

\makeatletter
\DeclareRobustCommand\Equiv{\mathrel{%
 \mathchoice
 {\Equiv@\textfont\displaystyle{.45}}
 {\Equiv@\textfont\textstyle{.45}}
 {\Equiv@\scriptfont\scriptstyle{.5}}
 {\Equiv@\scriptscriptfont\scriptscriptstyle{.55}}
}}
\newcommand{\Equiv@}[3]{%
 \rlap{\raisebox{#3\fontdimen5#12}{$\m@th#2 = $}}%
 \raisebox{-#3\fontdimen5#12}{$\m@th#2 = $}%
}
\makeatother

\def\eps{{\varepsilon}}

\newcommand{\dirichlet}{\mathcal{E}}

\newcommand{\Var}{{\operatorname{Var}}}
\newcommand{\Cov}{\operatorname{Cov}}

\newcommand{\Expectedsub}[2]{ \mathbb{E}_{#1}\left[#2\right] }
\newcommand{\Expectedsubnop}[2]{ \mathbb{E}_{#1}[#2] }
\newcommand{\Expected}[1]{ \mathbb{E}\left[#1\right] }
\newcommand{\Expectedbig}[1]{ \mathbb{E}\big[#1\big] }
\newcommand{\ExpectedBig}[1]{ \mathbb{E}\Big[#1\Big] }
\newcommand{\Expectednop}[1]{ \mathbb{E}[#1] }
\newcommand{\sigmabf}{ \boldsymbol{\sigma} }
\newcommand{\gradvec}[1]{{\nabla}^{#1}_{\!\!\sigmabf}f }
\newcommand{\portproj}[1]{\Pi_{#1}}
\newcommand{\portprob}[2]{p_{#1}^{\sigmabf}(#2)}
\newcommand{\portgrad}[1]{\widetilde{\nabla}^{#1}_{\!\!B, \sigmabf}f}
\newcommand{\Atil}{A_{\partial}}
\newcommand{\Gtil}{G_{\partial}}
\newcommand{\gradveccomplete}[1]{\widehat{\nabla}_{\!\!\sigmabf}^{#1} f}
\newcommand{\covop}[1]{\Lambda_{#1, \taubf}}
\newcommand{\covoppure}[1]{\Lambda_{#1}}
\newcommand{\scalarcov}{K_{B, \taubf}}

\renewcommand{\Var}[1]{ \mathbb{V}\mathrm{ar}\left[#1\right] }
\newcommand{\Varbig}[1]{ \mathbb{V}\mathrm{ar}\big[#1\big] }

\newcommand{\Varsub}[2]{ \mathbb{V}\mathrm{ar}_{#1}\left[#2\right] }
\newcommand{\Varnop}[1]{ \mathbb{V}\mathrm{ar}[#1] }

\renewcommand{\Pr}[1]{ \mathbb{P}[#1] }
\newcommand{\Prlr}[1]{ \mathbb{P}\left[#1\right] }
\newcommand{\Prsub}[2]{ \mathbb{P}_{#1}[#2] }

\newcommand{\taubf}{\boldsymbol{\tau}}

\newcommand{\DTV}[2]{\-D_{\mathrm{TV}}\left({#1},{#2}\right)}
\newcommand{\dist}{\mathrm{dist}}

\renewcommand{\epsilon}{\varepsilon}
\renewcommand{\emptyset}{\varnothing}

\renewcommand{\eps}{\varepsilon}

\newcommand{\swap}[1]{ D_{u \rightarrow v} }

\newcommand{\Ltil}{ \widetilde{L} }

\newcommand{\gap}[1]{ \mathrm{gap}({#1}) }

\newcommand{\port}[1]{\rho_{B,s}}

\newcommand{\lambdamax}{\lambda_{\mathrm{max}}}
\newcommand{\lambdamin}{\lambda_{\mathrm{min}}}

\newcommand{\lowdegree}{\overline{\Delta}}
\newcommand{\lowdegreedef}{d + d^{\frac{1}{2} + \eta}}
\newcommand{\numcolors}{d + d^{\frac{1}{2} + \varepsilon}}

\title{Sampling Colorings on $G(n, d/n)$ via Glauber Dynamics}
\author{Leon Schiller\thanks{Email: \texttt{leon.schiller@hpi.de}, 
Hasso Plattner Institute, University of Potsdam, Germany.}} 
\date{}

\begin{document}
\maketitle

\begin{abstract}

We study the mixing time of Glauber dynamics for sampling uniformly from the proper $q$-colorings of sparse Erd\H{o}s--R\'enyi random graphs $G(n,d/n)$. Our main result establishes that with probability $1-o(1)$ over the draw of $G\sim G(n, d/n)$, Glauber dynamics is rapidly mixing for every $q\geq d+d^{1/2+\varepsilon}$, given any fixed $0<\varepsilon\leq1/2$ and any sufficiently large fixed $d$. In particular, the surplus of colors above $d$ is sublinear and close to $\sqrt{d}$. The achieved color range therefore approaches the tree uniqueness threshold and improves over the previously best achievable number of colors $q > 1.763 d$ due to Efthymiou, Hayes, Stefankovic, and Vigoda \cite{efthymiou2018sampling}. Moreover, our mixing-time bound is $O(n^{2+o_d(1)})$, and hence almost quadratic as $d\rightarrow\infty$.

Our proof is based on establishing a spectral gap for block dynamics resulting from a suitable decomposition of $G(n,d/n)$ into blocks and singleton vertices, ensuring that every high-degree vertex is encapsulated within a block. Our main technical contribution is a novel way of establishing rapid mixing of block dynamics via an integrated Bochner-Bakry-\'{E}mery approach, which we expect to have applications beyond sampling colorings. We further devise a new block decomposition of $G(n,d/n)$ that strengthens and simplifies similar constructions used in previous work. Finally, to compare block and Glauber dynamics, we establish a novel recursive separator decomposition of tree-like blocks, together with a coupling argument that avoids the need for nearly uniform color marginals at vertices adjacent to hubs. 



\end{abstract}

\clearpage
\tableofcontents
\clearpage 

\newcommand{\tmix}{\tau_{\mathrm{mix}}}

\section{Introduction}

Sampling from distributions defined by local constraints is a fundamental problem in randomized algorithms, probability, and statistical mechanics. Proper graph colorings provide a canonical example: given a graph $G$ and $q$ colors, the task is to generate a coloring uniformly over all assignments in which adjacent vertices receive different colors.

This problem is widely studied, most commonly on graphs of maximal degree $\Delta$ and with some number of colors $q > \Delta$. For example, we mention the work of Jerrum~\cite{jerrum1995simple} who gave an efficient sampler for $q>2\Delta$; the work of Vigoda~\cite{vigoda2000improved} who improved this range to $q>11\Delta/6$, and the work of Carlson and Vigoda~\cite{carlson2025flip} who further reduced it to $q\geq1.809\Delta$ for $\Delta\geq125$. 
With the additional assumption that the graph has girth at least five, Chen and Liu~\cite{chen2026spectral} very recently established rapid mixing of Glauber dynamics for $q\geq(1+\varepsilon)\Delta$, for every fixed $\varepsilon\in(0,1)$ and sufficiently large $\Delta$ depending on $\varepsilon$.

In contrast to graphs of bounded degree, sparse random graphs pose another natural setting for studying how the geometry of the constraints affects the computational feasibility of sampling proper colorings. This setting is the focus of this work. Specifically, we consider the widely known Erd\H{o}s--R\'enyi graph $G(n,d/n)$, in which each edge is present independently with probability $d/n$ and $d$ is fixed as $n$ grows. Its average degree is close to $d$, but its maximum degree grows with $n$ such that classical mixing results that require a number of colors proportional to the maximum degree do not yield any fixed color palette on this class of graphs. The challenge is therefore to obtain a sampler that works under a number of colors governed by the average degree $d$ instead of the maximum degree $\Delta$.

The arguably most natural and most widely used sampling algorithm is the so-called \emph{Glauber dynamics} from the class of \emph{Markov Chain Monte Carlo} (MCMC) methods dating back to the seminal work of \citeauthor{metropolis1953equation}~\cite{metropolis1953equation}. Here, at each step, we choose a vertex uniformly at random and resample its color among those absent from its neighborhood. The uniform distribution on proper colorings is the stationary measure for this Markov chain. Its mixing time measures how quickly these local updates forget the initial coloring and approach equilibrium. This makes Glauber dynamics interesting both as a sampling algorithm and as a model of relaxation under local interactions. Rapid mixing from every initial coloring also gives information about whether the space of solutions can be explored efficiently through local moves.

Moreover, the statistical physics literature predicts several phase transitions regarding the space of proper colorings associated to $G(n, d/n)$ as the average degree $d$ increases~\cite{zdeborova2007phase,krzakala2007gibbs}. For large $q$, the cavity method predicts a \emph{clustering transition} at $d\sim q\log q$, where typical colorings split into exponentially many well-separated clusters. At a larger average degree, a \emph{condensation transition} concentrates most of the uniform measure on a small number of dominant clusters, before proper colorings eventually disappear near $d\sim2q\log q$. These changes in equilibrium geometry suggest potential barriers to efficient local sampling e.g. via Glauber dynamics. A particularly natural such barrier is \emph{Gibbs uniqueness}, which measures for which number of colors local color marginals become insensitive to distant boundary conditions; on infinite $d$-regular trees this leads to the scale $q\approx d$~\cite{krzakala2007gibbs}. This motivates the question of whether a Glauber dynamics chain on $G(n,d/n)$ remains rapidly mixing whenever $q$ is in the uniqueness regime.

Considerable progress has been made towards proving rapid mixing for Glauber dynamics with a number of colors depending only on $d$. This line of work starts with Mossel and Sly~\cite{mossel2010gibbs} who established polynomial mixing for a sufficiently large number $q=q(d)$ independent of $n$. Efthymiou~\cite{efthymiou2014mcmc} reduced the color requirement to a quantity linear in $d$, using a weighted block decomposition that shields high-degree vertices. The work of Efthymiou, Hayes, \v{S}tefankovi\v{c}, and Vigoda~\cite{efthymiou2018sampling} combined this approach with local uniformity estimates to obtain rapid mixing for $q\geq(\alpha_*+\varepsilon)d$, for every fixed $\varepsilon>0$ and sufficiently large $d$, where $\alpha_*\approx1.7632$ solves $\alpha_*=\exp(1/\alpha_*)$.

Different, non-MCMC-based sampling algorithms that work with fewer colors already exist. In particular, Efthymiou~\cite{efthymiou2016simple} gave a polynomial-time algorithm for $q\geq(1+\varepsilon)d$ whose output is within $n^{-\Omega(1)}$ total variation distance of uniform, with high probability over the graph. However, this still leaves open whether Glauber dynamics mixes rapidly whenever $q$ is above the uniqueness threshold. Moreover, besides being more natural and simpler algorithms, MCMC methods possess some further advantages over the method studied in \cite{efthymiou2016simple}. In particular, the stated $n^{-\Omega(1)}$ total variation error achieved by the algorithm in \cite{efthymiou2016simple} is tied to the instance size, whereas MCMC-based methods permit an arbitrarily prescribed accuracy by running the chain longer.

In this work, we improve upon the color requirement for rapid mixing of Glauber dynamics introduced in \cite{efthymiou2018sampling}.

\subsection{Main result}

Our main result is an almost quadratic mixing time bound for sampling proper colorings on $G(n,d/n)$ with only $q \ge d + d^{\varepsilon}\sqrt{d}$ colors via Glauber dynamics.

\begin{theorem}[Main theorem]
\label{thm:main}
    The following holds for all $0<\varepsilon\leq1/2$, sufficiently large $d$, and all integers $q \ge d + d^{\varepsilon}\sqrt{d}$.
    With probability $1 - o(1)$ over the draw of $G \sim G(n, d/n)$, for every $0<\delta<1$, the mixing time $$
    \tmix(\delta) \coloneqq \max_{\sigmabf} \min \Big\{t \in \mathbb Z_{\geq0} \mid \DTV{P_{\mathrm{GD}}^t(\sigmabf, \cdot)}{ \mu_{q, G} } \le \delta \Big\}
    $$ of Glauber dynamics for proper $q$-colorings on $G$ satisfies
    $$ 
        \tmix(\delta) = O_{d, \varepsilon} \left( n^{2 + \zeta} \bigl(1+\log(1/\delta)\bigr) \right) \qquad \text{with } \zeta \coloneqq \frac{C_\varepsilon\log(q)}{d^{\varepsilon}},
    $$ where $C_\varepsilon > 0$ is a constant that only depends on $\varepsilon$. Here, $P_{\mathrm{GD}}^t(\sigmabf,\cdot)$ denotes the distribution obtained after running Glauber dynamics for $t$ steps from the initial coloring $\sigmabf$.
\end{theorem}

In particular, this improves upon the previous best color range of $q > 1.7632 d$ obtained by \citeauthor{efthymiou2018sampling} \cite{efthymiou2018sampling}. Very recently, during the preparation of this manuscript, the work \cite{li2026rapid} by \citeauthor{li2026rapid} appeared in which rapid mixing on $G(n,d/n)$ for $q \ge (1 + \eps)d$ colors and $d$ sufficiently large was established independently. In contrast, our \Cref{thm:main} uses asymptotically fewer colors ($d + d^{\varepsilon} \sqrt{d}$ vs $(1 + \varepsilon)d$ whenever $\varepsilon < 1/2$) and obtains a better asymptotic mixing time of $n^{2 + o_d(1)}$ (assuming $q = O(d)$). We give more details as to how our result and our techinques relate to \cite{li2026rapid} in \Cref{sec:relatedwork}.

\subsection{Technical overview}

The proof of \Cref{thm:main} establishes a spectral gap for Glauber dynamics by studying suitably defined block dynamics. To this end, we establish a decomposition of the vertices of $G(n, d/n)$ into blocks of one or multiple vertices such that vertices in each block are updated collectively. Specifically, to establish rapid mixing of Glauber dynamics for only $q \ge d + d^{\varepsilon}\sqrt{d}$ colors, our proof combines the following three main ingredients.
\begin{enumerate}
    \setlength{\itemsep}{0cm}
    \item \textbf{Block decomposition of $G(n, d/n)$ that encapsulates hubs.} To define the block dynamics, we use a modified version of the decomposition employed by \citeauthor{efthymiou2018sampling} \cite{efthymiou2018sampling}. We identify a set $S$ of \emph{singleton vertices}, each of degree at most $\lowdegree \coloneqq \lowdegreedef < q$, with $\eta = \varepsilon/2$. After removing these vertices, the graph decomposes into small connected components, each of which forms a block $B$. Moreover, the singleton vertices in $S$ have the \emph{breakpoint property}, which implies that any path emanating from $s \in S$ contains hubs (i.e., vertices of degree $> \lowdegree$) only sparsely and sufficiently far away from $s$. This ensures that all hubs are encapsulated deep within a block $B$ and separated from the vertices in $S$ by a long buffer of low-degree vertices. Furthermore, every block $B$ is tree-like and contains at most one cycle. The resulting block decomposition is similar to that of \cite{efthymiou2018sampling}, but uses different breakpoint weights to obtain a stronger block diameter bound at a smaller cutoff degree $\lowdegree$.
    \item \textbf{Spectral gap for block dynamics via an integrated Bochner-type inequality.} To establish a spectral gap for the resulting block dynamics, we use an approach based on the \emph{Bochner--Bakry--\'{E}mery} framework. The underlying idea is to expand the continuous-time generator $L$ of block dynamics into a sum of singleton and block update generators to establish that $L^2 - \kappa L \succeq 0$ for some $\kappa > 0$. This approach is commonly used to establish spectral gaps for Markov diffusion operators \cite{bakry2013analysis,bakry2006diffusions} and has very recently attracted considerable attention in the study of discrete Markov chains \cite{goebel2026simple,guo2026approximating,wang2026optimal,chen2026spectral,boban2026rank,chen2026spectral2,miao2026sampling,chen2026high,haqi2026sampling}.

    Our approach is explicitly inspired by the work of \citeauthor{haqi2026sampling} \cite{haqi2026sampling}, who establishes a spectral gap for Glauber dynamics for proper colorings on line graphs. To adapt these ideas to block dynamics, the most important step is to devise a \emph{block elimination lemma} showing that, at the level of the Bochner inequality, a block approximately behaves like a collection of low-degree singleton vertices (the \emph{ports}) without internal edges. To show this, we exploit novel correlation decay estimates that take the form of an eigenvalue bound on a certain block covariance operator. This bound exploits the tree-like structure of blocks together with the breakpoint property of the boundary vertices in $S$.

    After block elimination, the remaining terms in the expansion of $L$ can be aggregated into a quadratic form involving the adjacency matrix of a low-degree subgraph of $G$. Establishing $L^2 - \kappa L \succeq 0$ then reduces to a lower bound on the smallest eigenvalue of this matrix.
    \item \textbf{Comparison of block and Glauber dynamics.} Finally, we transfer rapid mixing from block dynamics to Glauber dynamics. The main estimate needed for this step is a bound on the within-block relaxation time that is uniform over all boundary conditions. To this end, we recursively decompose each block into a tree of separators and study a form of block dynamics that updates the vertices within a separator individually and its child components collectively. Mixing-time estimates at each level of the separator tree can then be combined recursively to obtain a bound on the within-block relaxation time for Glauber dynamics.

    Overall, this approach is similar to the comparison used in \cite{efthymiou2018sampling}. The main difference is how we establish a mixing-time estimate for the block dynamics corresponding to a given level of the separator hierarchy. While the argument in \cite{efthymiou2018sampling} exploits nearly uniform marginals for the children of a hub, we work with significantly fewer colors, so this property does not carry over directly. Instead, we devise a different separator construction together with a carefully constructed coupling argument that avoids the need for color marginals that are close to uniform.
    We construct a common target coloring for the separator and couple an auxiliary Markov chain whose updates follow a carefully chosen order. This allows us to bound the probability of agreement from below. The auxiliary chain can then be compared to the block dynamics at the corresponding level of the separator hierarchy to yield a bound on the desired relaxation time.
\end{enumerate}
Hence the overall idea is that the graph is decomposed into buffered, tree-like blocks whose properties enable an elimination of its internal edges at the level of the Bochner inequality. After elimination, the remaining interactions are controlled via a bound on the smallest eigenvalue of a low-degree subgraph of $G$. Block and Glauber dynamics are then compared via a novel, recursively applied coupling argument.

While steps $1$ and $3$ also require new ideas over the previous work \cite{efthymiou2018sampling}, 
the main technical contribution of this paper is the Bochner-type inequality we devise in step $2$. We further believe the ideas developed in this step to be useful beyond the scope of colorings, i.e., for studying block dynamics associated to other spin systems such as e.g. the hardcore or antiferromagnetic Potts model. For this reason, the remainder of this section is concerned with unpacking how the Bochner inequality is assembled and in particular the block elimination step. 

\subsubsection{The direct Bochner method for graph colorings}

We start by outlining how the Bochner-type method can establish spectral gaps using the example of graph colorings as considered in \cite{haqi2026sampling}, which is the basis of our proof.

\paragraph{Introduction to the Bochner method}
The overall idea underlying the Bochner method is to establish a spectral gap by showing that an average gradient energy (or \emph{Dirichlet form}) associated to our Markov chain decays exponentially in the time $t$ we run it. More precisely, for a reversible chain on $\Omega$ with stationary measure $\mu$, transition rates $r(\sigmabf,\sigmabf')$, and generator $L$, this energy is
\[
    \dirichlet_L(f,f)=\langle f,Lf\rangle_\mu
       =\frac12\Expectedsub{\mu}{
          {\sum}_{\sigmabf'\in\Omega}r(\sigmabf,\sigmabf')
          \bigl(f(\sigmabf')-f(\sigmabf)\bigr)^2}.
\]
Thus it measures squared changes in $f$ across possible transitions, weighted by their rates and averaged over a stationary configuration $\sigmabf\sim\mu$. Writing $f_t=P_tf=e^{-tL}f$, reversibility gives $\frac{\mathrm d}{\mathrm dt}\dirichlet_L(f_t,f_t)=-2\langle Lf_t,Lf_t\rangle_\mu$. To obtain exponential decay of this energy (via Grönwall's lemma), it therefore suffices to show that
$$
    L^2 - \kappa L \succeq 0 \text{ or equivalently } \langle Lf, Lf \rangle \geq\kappa \langle f,Lf \rangle
$$ for some $\kappa > 0$. 
This gives $\dirichlet_L(f_t,f_t)\leq e^{-2\kappa t}\dirichlet_L(f,f)$. If the kernel consists only of constants, it also implies $\gap{L}\geq\kappa$. 

\paragraph{Direct Bochner method for graph colorings}
In discrete settings such as Glauber dynamics for spin systems on graphs, the generator can be written as a sum of local terms $L=\sum_v L_v$ where each $L_v$ corresponds to locally updating the spin at vertex $v$. 
For example, in case of graph colorings, we have $L = \sum_{v \in V} L_v$ with
\begin{align*}
    L_vf(\sigmabf)
        =f(\sigmabf)-\Expected{f(\sigmabf)\mid\sigmabf_{V\setminus\{v\}}}
          =-a_v(\sigmabf)^{-1}\sum_{c\in A_v(\sigmabf)}\gradvec{c}(v),
\end{align*}
where $\gradvec{c}(v) = f(\sigmabf_{v\leftarrow c})-f(\sigmabf)$, $A_v(\sigmabf)=[q]\setminus\{\sigmabf_u:u\in N_G(v)\}$ is the set of available colors, $a_v(\sigmabf)=|A_v(\sigmabf)|$, and $\sigmabf_{v\leftarrow c}$ is the configuration obtained by recoloring $v$ with $c$. Expanding $\langle Lf, Lf\rangle$ and using $L_v^2 = L_v$ then yields \begin{align}\label{eq:bochnerexpand}
    \langle Lf, Lf\rangle = \sum_{v} \langle f, L_vf \rangle + \sum_{u \neq v} \langle L_uf, L_vf \rangle.    
\end{align}
While the first sum is always nonnegative and equal to $\langle f, L f\rangle$, the second sum contains potentially negative terms. 
To establish a spectral gap, the goal is now to compare the sum over cross terms of the form $\langle L_uf, L_vf \rangle$ to the diagonal terms in the first sum which is equal to $\langle f , L f \rangle$: once we show that $\sum_{u \neq v} \langle L_uf, L_vf \rangle \ge -(1-\kappa) \langle f, L f \rangle$, it follows that $L^2 - \kappa L \succeq 0$.

Since the local generator $L_v$ can be expressed in terms of \emph{gradients} $\gradvec{c}$, the cross and diagonal terms above can be re-written and re-organized in terms of such gradients as well. Because it simplifies the argument, this is done for the closely related accelerated generator
\begin{align}\label{eq:acceleratedsingletongen}
    \Ltil=\sum_{v\in V}H_v,\qquad
    H_v=a_vL_v,\qquad
    H_vf(\sigmabf)=-\sum_{c\in[q]}\gradvec{c}(v).
\end{align}
and under setting $\gradvec{c}(v)=0$ for unavailable colors. In this setting, the identity of \citeauthor{haqi2026sampling} \cite[Lemma~4]{haqi2026sampling} yields the following \emph{Bochner-type inequality}:
\begin{align}\label{eq:bochnerexample}
    \langle \Ltil f, \Ltil f \rangle - \kappa \langle f, \Ltil f \rangle \ge \frac{a - \kappa}{2} \sum_{c \in [q]} \Expected{ \| \gradvec{c} \|^2 } + \sum_{c \in [q]} \Expected{ (\gradvec{c})^{\!\top} \! \! A_{G} (\gradvec{c})  }    
\end{align}
where $a$ is a uniform lower bound on the number of available colors at every vertex under any exterior coloring and $A_G$ is the adjacency matrix of the graph $G$. As this is the basis for our block Bochner inequality, we re-prove this inequality in \Cref{sec:singletonstuff} and include all details omitted here.

The use of an inequality like \eqref{eq:bochnerexample} is as follows: to compare the norm of $\gradvec{c}$ to the quadratic form $(\gradvec{c})^{\!\top} \! \! A_{G} (\gradvec{c})$ we can now exploit a lower bound on the smallest eigenvalue $\lambdamin(A_G)$. In particular, if $\lambdamin(A_G) > -(a - \kappa)/2$ then a spectral gap follows. This is useful for colorings on line graphs as studied in \cite{haqi2026sampling} as $\lambdamin(A_G) \ge -2$. It also establishes a spectral gap for colorings on random $d$-regular graphs with $q \ge d + C\sqrt{d}$ colors and $C$ being sufficiently large, since $\lambdamin(A_G) \ge -O(\sqrt{d})$. Similar arguments also apply to the hardcore-model as shown in \cite{goebel2026simple}.

Thus, the Bochner method allows us to take into account the signed structure of pairwise interactions between spins and aggregate them suitably, while taking their collective structure (manifesting in the form of a bound on $\lambdamin(A_G)$) into account. This is the reason why it is particularly effective in average-case settings observed e.g. in \cite{boban2026rank,goebel2026simple,wang2026optimal}.

\subsubsection{ Block dynamics on $G(n, d/n)$}

We proceed by outlining why the above direct analysis does not work for $G(n, d/n)$ and why switching to block dynamics can help overcome the described challenges.

\paragraph{Challenges on $G(n, d/n)$}
While the argument described so far yields a simple proof for rapid mixing of colorings on random $d$-regular graphs with only $d+C\sqrt d$ colors for a sufficiently large constant $C$, the case of $G(n, d/n)$ is much more involved. In particular, we cannot guarantee sufficiently many available colors $a$ uniformly over all vertices and proper colorings without using a number of colors that grows with $n$. After all, $G(n, d/n)$ contains vertices of degree up to $\Theta\big(\frac{\log(n)}{\log\log(n)}\big)$ with high probability, while we wish to work with only $q \sim d$ colors. Consequently, a direct Bochner-type inequality as established by Haqi \cite{haqi2026sampling} and outlined above does not yield a useful spectral-gap estimate in this setting. 

Moreover, even if we had a uniform lower bound on $a$ as for random regular graphs, the smallest eigenvalue of $A_G$ for $G\sim G(n,d/n)$ is not bounded below by $-O(\sqrt d)$ uniformly in $n$. Indeed, for every fixed $d>0$, the extreme-eigenvalue estimates of \citeauthor{benaych2019largest} \cite{benaych2019largest} imply that, with high probability,
\[
    \lambdamin(A_G)=-(1+o(1))\sqrt{\frac{\log n}{\log\log n}}
       \qquad\text{as }n\to\infty.
\]
This behavior is caused by vertices of atypically high degree and their star-like neighborhoods: a star with $k$ leaves has smallest adjacency eigenvalue $-\sqrt{k}$, and the extreme eigenvalues of the sparse random graph are governed by such local structures. Thus the smallest eigenvalue tends to $-\infty$ even though the average degree remains bounded, giving a second obstruction to applying the above spectral criterion directly, whenever we work with a number of colors $q$ that does not grow with $n$.

\paragraph{Moving to block dynamics}
We resolve both problems by placing every hub, together with a buffer of low-degree vertices (vertices of degree at most $\lowdegreedef$), inside a block and resampling all colors in the block jointly. Therefore, a hub is no longer required to have many available colors when its neighbors are held fixed; this is only required for low-degree singletons and vertices within blocks that have neighbors outside of it (these are called \emph{ports}). At those vertices, at least $a = q-\lowdegree$ colors remain available, which remains uniformly bounded from below once we use $q \ge \numcolors$ colors. 

A block-elimination argument then replaces the internal block structure by interactions at the ports, using the diagonal block update terms and decay of correlations through the buffers. The remaining adjacency quadratic form contains only edges joining two singletons or a singleton to a port. Its smallest eigenvalue is still bounded below by $-O(\sqrt d)$ with high probability as we show in \Cref{lem:smallesteigenvalue}, so the hub-induced spectral outliers no longer enter the global estimate.

\paragraph{Block decomposition}
Our decomposition selects a set $S$ of low-degree vertices as singletons and takes the connected components of $G[V\setminus S]$ as blocks. Here low-degree means degree at most $\lowdegree=d+d^{1/2+\eta}$, for a fixed $0<\eta<\varepsilon$. Vertices of larger degree are called hubs. The decomposition has the following properties, formalized in \Cref{sec:blockassumptions}.
\begin{enumerate}
    \setlength{\itemsep}{0cm}
    \item \textbf{Tree-like blocks.} Every block is a tree or contains exactly one cycle.
    \item \textbf{Breakpoint property of singletons.} Every singleton is a \emph{global breakpoint}. To define this, give low-degree vertices weight $W(v)=\alpha^{-1}$ and hubs weight $W(v)=\alpha^A\deg_G(v)/d$, where $\alpha=C\log d>1$, $C>0$ is fixed, and $A$ is large enough. A vertex is a breakpoint if the product of weights along every self-avoiding path starting there is at most $1$. Informally, the large weight of each hub must be compensated by many low-degree vertices along the path.
    \item \textbf{Unique attachment to blocks.} Every singleton has at most one neighbor in any given block. Vertices inside a block that are adjacent to singletons are called \emph{ports}. By the breakpoint property of singleton vertices, all ports are of degree $\le \lowdegree$.
    \item \textbf{Small diameter.} Any two vertices of a block can be joined within that block by a path of at most $D\log n/d^{2\eta}$ edges, where $D$ is a constant independent of $d,n$.
\end{enumerate}
In \cref{sec:blocks}, we show that such a decomposition actually exists on $G(n,d/n)$.

\subsubsection{The Bochner method for block dynamics}
Using the block decomposition described above, we now proceed by explaining how the Bochner method allows us to establish a spectral gap for block dynamics. For convenience, we consider the accelerated generator of block dynamics
\[
    \Ltil=\sum_{s\in S}H_s+q\sum_{B\in\mathcal B}L_B,
    \qquad H_s=a_sL_s,
\]
where $\mathcal B$ is the collection of blocks and $a_s$ is the number of available colors at $s$. Here $L_Bf(\sigmabf)=f(\sigmabf)-\Expected{f\mid\sigmabf_{V\setminus B}}$ resamples all colors in $B$ uniformly conditional on the exterior, while $L_s$ is the single-site update introduced above. Thus singleton updates have rate $a_s$ and block updates have rate $q$. We study this choice of rates because it makes the Bochner expansion more convenient, similarly to the recoloring dynamics used by \citeauthor{haqi2026sampling} \cite{haqi2026sampling}.

Expanding the squared generator analogous to \eqref{eq:bochnerexpand} separates the following diagonal and cross terms:
\begin{align}
    \langle\Ltil f,\Ltil f\rangle
       &\ge \sum_{s\in S}\|H_sf\|^2
           +q^2\sum_{B\in\mathcal B}\|L_Bf\|^2 \notag\\
       &\quad+2\sum_{\{s,t\}\in E(G[S])}\langle H_sf,H_tf\rangle
           +2q\sum_{B\in\mathcal B}\sum_{s\in \partial B}
                    \langle H_sf,L_Bf\rangle,
    \label{eq:overviewblockexpansion}
\end{align}
where we omitted cross terms corresponding to nonadjacent singletons and all pairs of blocks because these can easily be shown to be non-negative.
The singleton interactions $\|H_sf\|^2$ and $\langle H_sf,H_tf\rangle$ can be related to the same gradients $\gradvec{c}$ as appear in \eqref{eq:bochnerexample} (after restricting to vertices in $S$). In fact, we can prove (cf. \Cref{cor:lowdegreebochner} for a closely related, formal statement)
\begin{align}\label{eq:bochnerexample2}
    \sum_{s\in S}\|H_sf\|^2 + 2\sum_{\{s,t\}\in E(G[S])}\langle H_sf,H_tf\rangle \ge \frac{a}{2}\sum_{c \in [q]} \Expected{ \| \gradvec{c} \|_2^2 } + \sum_{c \in [q]} \Expected{ (\gradvec{c})^{\!\top} A_{G[S]} (\gradvec{c}) }    
\end{align}
which is analogous to \eqref{eq:bochnerexample}.

\paragraph{Block-singleton interactions lead to port gradients} To obtain a similar decomposition for the terms in \eqref{eq:overviewblockexpansion} that can be unified with \eqref{eq:bochnerexample2}, we wish to relate terms of the form $\|L_Bf\|^2$ and $\langle H_sf,L_Bf\rangle$ to some kind of gradient representation as well. 

As all interaction between a block $B$ and the rest of the graph happens via a port $\rho$ of $B$, it turns out that this can be done via gradients that only depend on the spin at the ports of $B$. To see how such gradients arise, define the projection $\portproj{\rho}f=\Expected{f\mid\sigmabf_{V\setminus B},\sigmabf_\rho}$ which eliminates all dependence of $f$ on spins within $B$ except for the port $\rho \in B$. 
As all interaction between a block $B$ and a boundary vertex $s \in \partial B$ happens over one unique port $\rho = \rho_{B,s}$ (recall that every $s \in \partial B$ is incident to a unique port in $B$), we get (\Cref{lem:portproj})
$$
    \langle H_sf,L_Bf\rangle \ge \langle H_s\portproj{\rho}f,L_B \portproj{\rho}f \rangle
$$ and we can in fact express 
$$
    L_B\portproj{\rho}f(\sigmabf)
       =-\frac{1}{q}\sum_{c\in[q]}\portgrad{c}(\rho) \text{ with } \portgrad{c}(\rho)
       =q\,\portprob{\rho}{c}
          \left(\portproj{\rho}f( \sigmabf_{\rho\leftarrow c})
                    -\portproj{\rho}f(\sigmabf)\right).
$$ where $\portprob{\rho}{c}=\Pr{\sigmabf_\rho=c\mid\sigmabf_{V\setminus B}}$ is the conditional color probability at the port, and the gradient is set to zero whenever this probability is zero. The vector $\portgrad{c} \in \mathbb{R}^{|P_B|}$ is a \emph{port gradient}, indexed by the set $P_B$ of ports in $B$. $L_B\portproj{\rho}f(\sigmabf)$ can therefore be expressed as a sum over gradients similar to the accelerated singleton generator $H_v$ as stated in \eqref{eq:acceleratedsingletongen}.

Interaction terms of the form $\langle H_sf,L_Bf\rangle$ can therefore be expressed in terms of products of gradients and port gradients in the sense that (see \Cref{lem:completesquares}, \Cref{lem:approxcompletesquares} for formal statements)
$$
    \langle H_sf,L_Bf\rangle \gtrsim \frac{1}{q}\sum_{c \in [q]} \Expected{ \gradvec{c}(s) \portgrad{c}(\rho) } \text{ analogously to } \langle H_sf,H_tf\rangle \ge \sum_{c \in [q]} \Expected{ \gradvec{c}(s) \gradvec{c}(t) }.
$$
Here the first comparison is schematic: it represents a lower bound obtained by discarding nonnegative terms and allowing a controlled negative error that can be charged agains the diagonal terms later on. The formal treatment of this error is presented in \Cref{sec:blockcontribution} and afterwards.

\paragraph{Block Bochner identity}
The fact that a block only interacts via its ports and the role of port gradients in our above expression of $\langle H_sf,L_Bf\rangle$ suggest that all terms in \eqref{eq:overviewblockexpansion} could be expressed in terms of gradients and port gradients. Our goal would then be a Bochner inequality of the form 
\begin{align}\label{eq:blockbochnergoal}
    \langle\Ltil f,\Ltil f\rangle - \kappa \langle f,\Ltil f\rangle \ge  \alpha \sum_{c \in [q]} \ExpectedBig{ \| \gradveccomplete{c} \|_2^2 } + \beta \sum_{c \in [q]} \Expected{ \big(\gradveccomplete{c}\big)^{\!\top} \!\! \Atil \big(\gradveccomplete{c}\big) },
\end{align} for coefficients $\alpha>0$ and $\beta\geq0$, the combined gradient vector
$$
    \gradveccomplete{c}(v) = \begin{cases}
        \gradvec{c}(v) & \text{if } v \in S\\
        \portgrad{c}(v) & \text{if } v \text{ is a port in some block } B\\
        0 & \text{otherwise}
    \end{cases}
$$ and with $A_\partial$ being the adjacency matrix of the graph consisting of all edges that have at least one endpoint in $S$. This would be the analog of the Bochner inequality \eqref{eq:bochnerexample} proved by Haqi \cite{haqi2026sampling}. A lower bound on the smallest eigenvalue of $A_\partial$ that ensures $\alpha+\beta\lambdamin(A_\partial)\geq0$ makes the right-hand side nonnegative and thus implies a spectral gap for $\kappa>0$. We establish this inequality and the required eigenvalue bound formally in \Cref{sec:assembling}.

\paragraph{Establishing \eqref{eq:blockbochnergoal} via block elimination}

An intuitive way of understanding \eqref{eq:blockbochnergoal} is that every block $B$ acts approximately as if it had no internal edges. This is depicted in \Cref{fig:block-elimination}. 
One way of formalizing this statement that (together with what we have already established) leads to an inequality of the form \eqref{eq:blockbochnergoal} would be to show that the block diagonal term $\|L_B f\|^2 = \langle f, L_B f \rangle$ satisfies
$$
    \langle f, L_B f \rangle \ge \gamma \sum_{c \in [q]} \Expected{\| \portgrad{c} \|^2}
$$ for some $\gamma > 0$. It is not hard to see that the right term is closely related to the expected variance of $\portproj{\rho} f$ in the sense that
$
    \sum_{c \in [q]} \Expectednop{ \portgrad{c}(\rho)^2} \lesssim \Expected{ \Var{ \portproj{\rho} f \mid \sigmabf_{V\setminus B} } },
$ as shown in \Cref{lem:portgradientstovariance} (here, the hidden constant is not absolute and may depend on $q$ and $a$). Therefore, it essentially suffices to show 
\begin{align}\label{eq:blockelim}
    \sum_{\rho \in P_B} \Expected{ \Var{ \portproj{\rho} f \mid \sigmabf_{V\setminus B} } } \le (1 + \varepsilon) \langle f, L_B f \rangle = (1 + \varepsilon) \Expected{\Var{ f \mid \sigmabf_{V\setminus B} }}.
\end{align}
Intuitively, this amounts to showing that the variance of $f$ incurred by resampling the entire block $B$ is not much smaller than the sum of the individual variances explained by fluctuations at every port. It follows by showing that the spins at each port are sufficiently un-correlated and thus formalizes the idea that a port acts as if it had no internal edges. This is what we mean by \emph{block elimination}. It poses the main work necessary for establishing \eqref{eq:blockbochnergoal}. 

\input{figures/block-elimination.tex}

\paragraph{Block elimination by bounding the covariance operator}
To see how \eqref{eq:blockelim} connects to de-correlation of ports, suppose that no internal edges exist in $B$. Then, conditional on any coloring outside of $B$, ports are independent and \eqref{eq:blockelim} holds with $\eps = 0$. To see this, consider a centered $f$ and note that   
$$
    0 \le \Expected{ \Big(f - \sum_{\rho} \portproj{\rho}f \Big)^2 } = \Expected{f^2} + \sum_{\rho} \Expected{(\portproj{\rho}f)^2} - 2 \sum_{\rho} \Expected{f \ \portproj{\rho}f} = \Var{f} - \sum_{\rho} \Var{\portproj{\rho}f}
$$ where $\Expected{\cdot}$ is over the coloring in $B$ conditional on any fixed exterior coloring, and we used $\Expected{f \ \portproj{\rho}f} = \Expected{(\portproj{\rho}f)^2}$ as well as independence of the ports $\rho$.

More generally, the left and right hand side of \eqref{eq:blockelim} can be related via the largest eigenvalue of a \emph{port covariance operator} $\covop{B}$. This operator measures correlations between any pair of ports in $B$ under the boundary condition imposed by $\taubf$.
It can be represented as an operator-valued matrix of size $|P_B| \times |P_B|$ such that $(\covop{B})_{i, j}$ is an operator that measures the covariance between any pair of centered functions that only depend on the spin at $i$ and $j$, respectively.
Concretely, for $f_i$ depending only on the spin at $i$ and $f_j$ depending only on the spin at $j$, we define 
$$
    (\covop{B})_{i, j} f_j (c) = \Expected{ f_j(\sigmabf_j) \mid \sigmabf_i = c}  \text{ such that } \langle f_i, (\covop{B})_{i, j} f_j \rangle = \Expected{f_i(\sigmabf_i) f_j(\sigmabf_j)}
$$ where expectations are taken conditional on $\taubf$. The formal definition is found in \Cref{sec:blockelim}.
As we show in \Cref{lem:variancetoeigenvalue}, we have 
$$
    \sum_{\rho \in P_B} \Expected{ \Var{ \portproj{\rho} f \mid \sigmabf_{V\setminus B} } } \le \max_{\taubf}\lambdamax(\covop{B}) \cdot \langle f, L_B f \rangle
$$ and thus establish our block elimination statement \eqref{eq:blockelim} once we bound $\lambdamax(\covop{B}) \le 1 + \varepsilon$.

\paragraph{Covariance norm analysis via correlation decay on trees and breakpoint estimates}
We bound the norm of $\covop{B}$ in \Cref{sec:covarianceanalysis} by exploiting the fact that on a tree, correlations between ports $i, j$ factor into products of edge operators along the unique path from $i$ to $j$. An edge joining two low-degree vertices contracts correlations by a factor at most $1/a \le 1/\sqrt{d}$. For a tree block $B$, the norm of $\covop{B}$ can therefore be bounded by the norm of the scalar matrix 
$$
        (\scalarcov)_{i, j} = \begin{cases}
            \prod_{e \in P(i,j)}\beta_e & \text{if } i \neq j \\
            1 & \text{otherwise,}
        \end{cases}
$$ where $\beta_e$ measures how much correlations contract along edge $e$ and $P(i,j)$ is the unique path from $i$ to $j$. The spectral norm of this matrix can be bounded by rooting our tree block at an arbitrary vertex and controlling an associated sum over root-to-port paths. The cost of this path sum can then be bounded by taking advantage of the breakpoint properties of $\partial B$ as guaranteed by our block decomposition. The argument is carried out in \Cref{lem:treeeigenvalue}.
We subsequently extend it to unicyclic blocks in \Cref{lem:unicycliceigenvaluelowcycle} and \Cref{lem:unicycliceigenvaluehubcycle}.


\subsection{Comparison to concurrent work}\label{sec:relatedwork}

\input{related_work.tex}

\subsection{AI disclosure}

The main ideas of the proofs were developed using long interactions with GPT-5.6 Sol and GPT-6.1 Sol after instructing it to devise a Bochner-type inequality for block dynamics for proper colorings. The same models were also used for assisting with the writing process of this document. All mathematical claims and calculations were independently verified and digested by the author who takes full responsibility for the final content.

\section{Preliminaries}\label{sec:preliminaries}

\input{prelims.tex}

\section{Spectral gap for block dynamics}\label{sec:blockbochner}

\input{blockdynamics.tex}

\section{Blocks and breakpoints in $G(n,d/n)$}\label{sec:blocks}

\input{breakpoints.tex}

\section{Comparing block dynamics to Glauber dynamics}

\input{comparison.tex}

\bibliographystyle{alphanat}
\begingroup
\footnotesize
\raggedright
\bibliography{references}
\endgroup


\end{document}

%% file: figures/block-elimination.tex
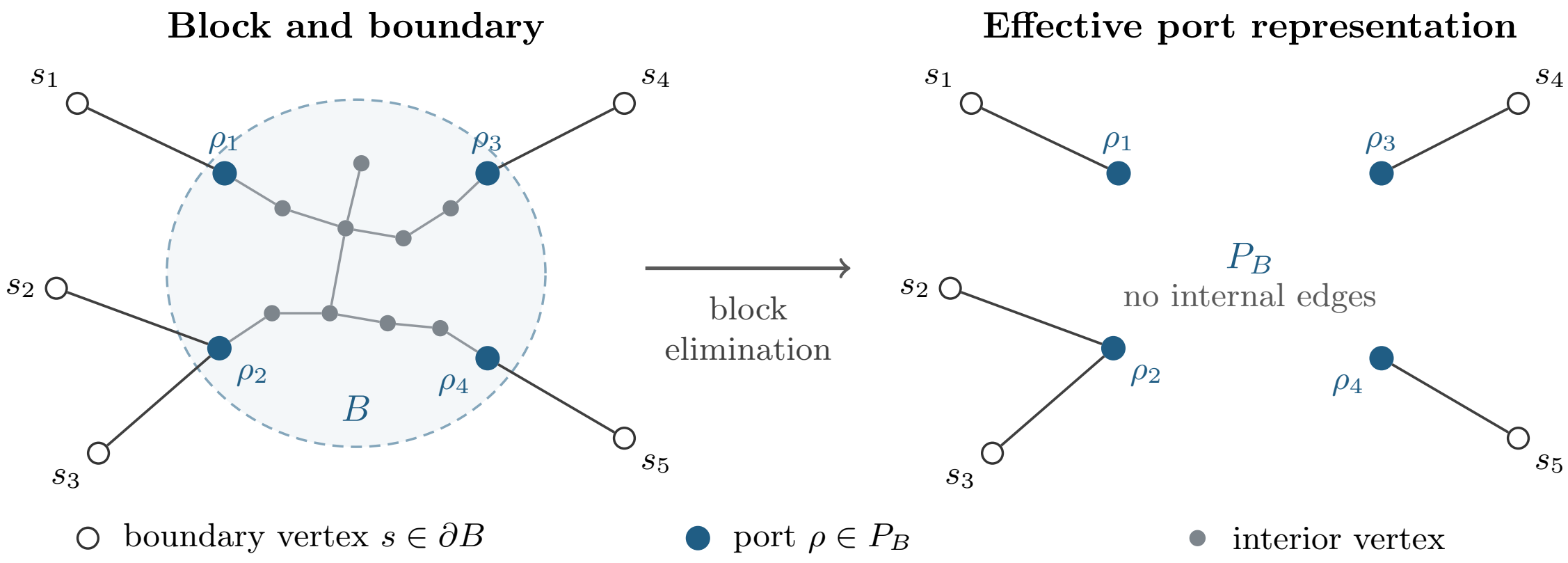
\begin{figure}[htbp]
    \centering
    \input{figures/block-elimination.tikz}
    \caption{Schematic illustration of our block elimination. The shaded block $B$ contains
    the ports $P_B$ and interior vertices. Its boundary $\partial B$ consists
    of singleton vertices. At the level of the Bochner inequality, the
    block approximately behaves as if it had no internal edges. Formally, this is achieved by establishing a statement like \eqref{eq:blockelim} via a spectral norm bound on the port covariance operator. The norm bound relies on the breakpoint properties of the block's boundary as guaranteed by our block decomposition.}
    \label{fig:block-elimination}
\end{figure}

%% file: figures/block-elimination.tikz
\begingroup
\definecolor{beBlue}{RGB}{32,93,132}
\definecolor{beGray}{RGB}{125,133,140}
\begin{tikzpicture}[
    x=1cm,y=.95cm,
    font=\small,
    be/boundary/.style={circle,draw=black!80,fill=white,
        line width=.65pt,minimum size=5.6pt,inner sep=0pt},
    be/port/.style={circle,draw=beBlue,fill=beBlue,
        minimum size=6.2pt,inner sep=0pt},
    be/interior/.style={circle,draw=beGray,fill=beGray,
        minimum size=4pt,inner sep=0pt},
    be/exterior edge/.style={draw=black!75,line width=.7pt},
    be/interior edge/.style={draw=beGray!85,line width=.65pt},
    be/label/.style={font=\footnotesize,inner sep=2pt}
]
    \filldraw[fill=beBlue!5,draw=beBlue!55,line width=.65pt,
        dash pattern=on 3pt off 2pt]
        (0,0) ellipse [x radius=1.8cm,y radius=1.65cm];
    \node[text=beBlue] at (0,-1.35) {$B$};

    \foreach \shift/\tag in {0/L,8.5/R} {
        \begin{scope}[xshift=\shift cm]
            \node[be/boundary,label={[be/label]above left:$s_1$}]
                (\tag-s1) at (-2.65,1.7) {};
            \node[be/boundary,label={[be/label]left:$s_2$}]
                (\tag-s2) at (-2.85,-.15) {};
            \node[be/boundary,label={[be/label]below left:$s_3$}]
                (\tag-s3) at (-2.45,-1.8) {};
            \node[be/boundary,label={[be/label]above right:$s_4$}]
                (\tag-s4) at (2.55,1.7) {};
            \node[be/boundary,label={[be/label]below right:$s_5$}]
                (\tag-s5) at (2.55,-1.65) {};

            \node[be/port,label={[be/label,text=beBlue]above:$\rho_1$}]
                (\tag-r1) at (-1.25,1) {};
            \node[be/port,label={[be/label,text=beBlue]below right:$\rho_2$}]
                (\tag-r2) at (-1.3,-.75) {};
            \node[be/port,label={[be/label,text=beBlue]above:$\rho_3$}]
                (\tag-r3) at (1.25,1) {};
            \node[be/port,label={[be/label,text=beBlue]below left:$\rho_4$}]
                (\tag-r4) at (1.25,-.85) {};

            \foreach \s/\r in {1/1,2/2,3/2,4/3,5/4}
                \draw[be/exterior edge] (\tag-s\s)--(\tag-r\r);
        \end{scope}
    }

    \foreach \name/\xx/\yy in {
        a/-.7/.65,b/-.1/.45,c/.45/.35,d/.9/.65,
        e/-.8/-.4,f/-.25/-.4,g/.8/-.55,h/.3/-.5,i/.05/1.1}
        \node[be/interior] (L-\name) at (\xx,\yy) {};
    \draw[be/interior edge] (L-r1)--(L-a)--(L-b)--(L-c)--(L-d)--(L-r3);
    \draw[be/interior edge] (L-r2)--(L-e)--(L-f)--(L-b);
    \draw[be/interior edge] (L-r4)--(L-g)--(L-h)--(L-f);
    \draw[be/interior edge] (L-b)--(L-i);

    \draw[->,draw=black!65,line width=1pt] (2.75,.05)--(4.7,.05);
    \node[align=center,font=\footnotesize,text=black!75] at (3.73,-.55)
        {block\\elimination};

    \node[text=beBlue] at (8.5,.15) {$P_B$};
    \node[font=\footnotesize,text=black!65] at (8.5,-.25)
        {no internal edges};

    \node[font=\small\bfseries] at (0,2.45)
        {Block and boundary};
    \node[font=\small\bfseries] at (8.5,2.45)
        {Effective port representation};

    \node[be/boundary] at (-2.55,-2.65) {};
    \node[anchor=west,font=\footnotesize] at (-2.35,-2.65)
        {boundary vertex $s\in\partial B$};
    \node[be/port] at (3.25,-2.65) {};
    \node[anchor=west,font=\footnotesize] at (3.45,-2.65)
        {port $\rho\in P_B$};
    \node[be/interior] at (8,-2.65) {};
    \node[anchor=west,font=\footnotesize] at (8.2,-2.65)
        {interior vertex};
\end{tikzpicture}
\endgroup

%% file: related_work.tex
Independently and concurrently, Li and Yang~\cite{li2026rapid} established rapid mixing of single-site Glauber dynamics near the tree-uniqueness threshold. More precisely, their Theorem~1.1 states that, for every fixed $\varepsilon\in(0,1)$ and sufficiently large $d=d(\varepsilon)$, all $q\geq(1+\varepsilon)d$ satisfy, with high probability,
\[
    \tmix(\delta)
       \leq n^{C(\varepsilon,d,q)}\bigl(1+\log(1/\delta)\bigr),
       \qquad 0<\delta<1.
\]
Both their result and ours therefore approach a color-to-degree ratio of one as $d$ grows.

Our \Cref{thm:main} achieves a smaller asymptotic number of colors for every fixed $0<\varepsilon<1/2$ as $d \rightarrow \infty$. Moreover, it yields a stronger explicit mixing-time guarantee than the unspecified polynomial exponent in \cite[Theorem~1.1]{li2026rapid}.
Our mixing time further improves the excess above the quadratic exponent in~\cite{efthymiou2018sampling}, since $\log d/d^\varepsilon=o(1/\log d)$ for every fixed $\varepsilon>0$ and $q = O(d)$.

\paragraph{Comparison of proof methods and their capabilities}

Remarkably, the techniques used in \cite{li2026rapid} are quite different from ours. The only common part is that both proofs build on block dynamics and a comparison of block to Glauber dynamics. However, the block decomposition in \cite{li2026rapid} is very different from ours as it uses overlapping blocks while our blocks are all disjoint. Moreover, \citeauthor{li2026rapid} establish rapid mixing of block dynamics via spectral- and coupling independence by extending the tree correlation-decay techniques from~\cite{chen2023strong}, while we use a Bochner-type inequality. Even the comparison between block and Glauber dynamics is substantially different: Li–Yang bound within-block relaxation times using connectivity and the size of the configuration space, while we exploit recursive separators and a sweep coupling to obtain sharper estimates.

\paragraph{Extensibility to other spin systems.}

We view it as a further strength of our method that the block decomposition we devise and the Bochner-type approach for establishing a spectral gap of block dynamics are not inherently limited to proper colorings, but form a concrete route towards the study of block and Glauber dynamics for different spin systems on random graphs.

The hard-core model on $G(n,d/n)$ is a natural candidate, and it is reasonable to expect that our proof leads to rapid mixing for activities far beyond the uniqueness threshold of the corresponding $d$-regular tree. This is the case for random regular graphs as shown in \cite{chen2025rapid} via localization schemes and in \cite{goebel2026simple} via a direct Bochner proof that exploits the smallest eigenvalue of $G$. Our block elimination argument enables the use of a similar eigenvalue estimate on a low-degree subgraph of $G(n,d/n)$ that has smallest eigenvalue $\ge - O(\sqrt{d})$ like a random $d$-regular graph. This supports the claim that it should be able to use spectral information to extend beyond uniqueness. Interactions with GPT 6.1 Sol further confirm that a simplified form of our argument plausibly yields rapid mixing for Glauber dynamics on the hard-core model up to activities of $\lambda \le O(d^{-1/2 - \varepsilon})$, which nearly reaches the $O(1/\sqrt{d})$ scale achieved for random regular graphs. Similar bounds seem to be attainable for the inverse temperature in the antiferromagnetic Potts model with a constant number of colors.

We therefore view it as a particular strength of our approach in comparison to \cite{li2026rapid} that spectral information remains attainable. This does not appear to be the case for the proof by \citeauthor{li2026rapid} \cite{li2026rapid}, which seems incapable of extending beyond uniqueness for the hard-core model.

%% file: prelims.tex
We introduce the notation used throughout the paper.

\subsection{Graphs}

For an integer $m\geq1$, write $[m]=\{1,\ldots,m\}$.
For a set $A$ and an element $v$, we sometimes write $A+v$ for $A\cup\{v\}$
and $A-v$ for $A\setminus\{v\}$.
All graphs are
finite, simple, and undirected. For $G=(V,E)$, let $n=|V|$,
$N_G(v)=\{u\in V:\{u,v\}\in E\}$, and $\deg_G(v)=|N_G(v)|$.
For $U\subseteq V$, write $G[U]$ for the induced subgraph and
\[
    \partial U=\{v\in V\setminus U:N_G(v)\cap U\neq\emptyset\}
\]
for its exterior vertex boundary. In particular, $\partial B$ consists
of singleton vertices when $B$ is a block in our partition.
Unless another graph is specified, degrees are always measured in $G$.

A self-avoiding path $P=(v_0,\ldots,v_\ell)$ has distinct vertices,
length $\ell$, and $\ell+1$ vertices; the single-vertex path is allowed.
Write $\dist_G(u,v)$ for graph distance and
$B_G(v,r)=\{u\in V:\dist_G(u,v)\leq r\}$ for the radius-$r$ ball.
The diameter of a connected subgraph is measured using paths within
that subgraph. A connected graph is \emph{unicyclic} if it has exactly
one cycle. In $G(n,d/n)$, each pair of vertices in $[n]$ forms an edge
independently with probability $d/n$. An event holds \emph{with high
probability} if its probability is $1-o(1)$ as $n\to\infty$, with $d$
and the other fixed parameters held constant.

\subsection{Colorings and conditional measures}

Let $\Omega_{q,G}$ be the set of proper $q$-colorings
$\sigmabf:V\to[q]$, and let $\mu_{q,G}$, abbreviated $\mu$, be the
uniform measure on this set whenever it is nonempty. We write
$\sigmabf_U$ for the restriction to $U\subseteq V$ and
$\sigmabf_{v\leftarrow c}$ for the configuration obtained by assigning
color $c$ to $v$. The available colors at $v$ are
\[
    A_v(\sigmabf)=[q]\setminus
                   \{\sigmabf_u:u\in N_G(v)\}.
\]
A boundary condition $\taubf$ on $V\setminus U$ is \emph{feasible}
if it extends to a proper coloring of $G$. For such $\taubf$, denote
by $\nu_{U,\taubf}$ the conditional law of $\sigmabf_U$ given
$\sigmabf_{V\setminus U}=\taubf$; we also write $\nu_{\taubf}$ when
the region is clear.

For a probability measure $\nu$ on a finite set, $L^2(\nu)$ denotes
the real functions on its support, with
\[
    \langle f,g\rangle_\nu=\Expectedsub{\nu}{fg},\qquad
    \|f\|_{2,\nu}^2=\langle f,f\rangle_\nu,\qquad
    \Varsub{\nu}{f}=\Expectedsub{\nu}{f^2}
                         -\Expectedsub{\nu}{f}^{\!2}.
\]
We omit the measure when it is clear, and write
$\Cov_\nu(f,g)=\Expectedsub{\nu}{fg}
-\Expectedsub{\nu}{f}\Expectedsub{\nu}{g}$.
For $U\subseteq V$, the notation
$\Varsub{U}{f}=\Varsub{\mu}{f\mid\sigmabf_{V\setminus U}}$
denotes a conditional variance, which is a function of the exterior
colors; for $U=\{v\}$ we write $\Varsub{v}{f}$.
For self-adjoint operators $A,B$ on the same space, $A\succeq B$
means $\langle f,(A-B)f\rangle\geq0$ for every $f$.

\subsection{Markov chains, spectral gaps, and mixing times}

A transition matrix $P$ on a finite state space $\Omega$ acts on
functions by $Pf(\sigmabf)=\sum_{\sigmabf'\in\Omega}
P(\sigmabf,\sigmabf')f(\sigmabf')$. A probability measure $\mu$ is
stationary if $\mu P=\mu$, and $P$ is \emph{reversible} with respect
to $\mu$ if
$\mu(\sigmabf)P(\sigmabf,\sigmabf')
=\mu(\sigmabf')P(\sigmabf',\sigmabf)$ for all states.
Reversibility makes $P$ self-adjoint on $L^2(\mu)$.
A chain is \emph{irreducible} if every state can be reached from every
other state. When discussing convergence, we assume irreducibility
and, in discrete time, aperiodicity. Heat-bath chains have a positive
holding probability and hence are aperiodic.

We use the \emph{positive generator} convention: a continuous-time
chain has semigroup $e^{-tL}$, so its usual infinitesimal generator
is $-L$. For a discrete-time chain, the corresponding positive
operator is $I-P$. The Dirichlet form and spectral gap are
\[
    \dirichlet_L(f,f)=\langle f,Lf\rangle_\mu,
    \qquad
    \gap{L}=\inf_{\Varsub{\mu}{f}>0}
                 \frac{\dirichlet_L(f,f)}{\Varsub{\mu}{f}}.
\]
Thus $\gap{L}\geq\gamma$ is the Poincar\'e inequality
$\Varsub{\mu}{f}\leq\gamma^{-1}\dirichlet_L(f,f)$ for every $f$.
For an irreducible reversible chain with at least two states, the gap
is the smallest positive eigenvalue of $L$, and its reciprocal is the
relaxation time. Our within-block relaxation time uses the harmless
truncation at $1$ specified in \Cref{def:withinblockrelaxation}.

For colorings, let $P_Uf=\Expectedsub{\mu}{f\mid
\sigmabf_{V\setminus U}}$ be the heat-bath update of $U$.
This is an orthogonal projection, and
$\langle f,(I-P_U)f\rangle_\mu
=\Expectedsub{\mu}{\Varsub{U}{f}}$.
In particular, discrete-time Glauber dynamics and its continuous-time
version with rate $1$ at each vertex satisfy
\[
    P_{\mathrm{GD}}=\frac1n\sum_{v\in V}P_{\{v\}},\qquad
    L_{\mathrm{Gl}}=\sum_{v\in V}(I-P_{\{v\}})
                  =n(I-P_{\mathrm{GD}}).
\]

For probability measures $\nu,\mu$ on $\Omega$, their total variation
distance is $\DTV{\nu}{\mu}=\frac12\sum_{\sigmabf\in\Omega}
|\nu(\sigmabf)-\mu(\sigmabf)|$. For $0<\delta<1$, define
\[
    \tmix(P,\delta)=\min\left\{t\in\mathbb Z_{\geq0}:
       \max_{\sigmabf\in\Omega}
       \DTV{P^t(\sigmabf,\cdot)}{\mu}\leq\delta\right\}.
\]
We abbreviate $\tmix(P)=\tmix(P,1/4)$ and write
$\tmix(\delta)=\tmix(P_{\mathrm{GD}},\delta)$ in
\Cref{thm:main}. In continuous time, the
same definition uses $t\geq0$ and $e^{-tL}$ in place of $P^t$.

%% file: blockdynamics.tex
In this section, we show that after choosing a suitable partition of $G$ into singleton vertices and blocks, the resulting block dynamics for colorings has a spectral gap. We start by introducing the assumptions on our block partition and afterwards study the resulting block dynamics via a Bochner-type inequality. Showing that a suitable block partition actually exists in $G(n, d/n)$ is the subject of \Cref{sec:blocks}.

\subsection{Block partition}

Our block partition is obtained by considering the set $\widehat{S}$ of global \emph{breakpoint vertices} as defined below. After identifying a set $S \subseteq \widehat{S}$ as the set of singleton vertices (i.e. vertices that are not part of any block and are updated individually), the set $\mathcal{B}$ of blocks is exactly the set of connected components left in $G$ after removing all vertices in $S$.
Global breakpoints in $\widehat S\setminus S$ are block vertices and are updated with their blocks.

To introduce this more formally, we first define the vertex weights and the breakpoint condition. To this end, we fix constants
\(0<\eta<1/2\), \(A>0\), and \(C>0\), and define
\[
    \lowdegree \coloneqq d+d^{1/2+\eta} , \text{ and }
    \alpha \coloneqq C\log d>1.
\]
With this, we define associate a vertex weight $W(v)$ with every vertex $v$ given by
\begin{align}\label{eq:weights}
    W(v) \coloneqq
    \begin{cases}
        \alpha^{-1}, & \deg_G(v)\leq \lowdegree,\\[1mm]
        \alpha^A\deg_G(v)/d, & \deg_G(v)>\lowdegree.
    \end{cases}
\end{align}
Based on these vertex weights, we can decide when a vertex $v$ is a breakpoint as captured by the following definition.
\begin{definition}[Breakpoints]\label{def:r-breakpoint}
    Let \(r\geq0\) be an integer. A vertex \(v\in V\) is an
    \emph{\(r\)-breakpoint} if, for every self-avoiding path
    \(P=(v_0,\ldots,v_\ell)\) with \(v_0=v\) and \(0\leq\ell\leq r\),
    \[
        W(P)\coloneqq\prod_{i=0}^{\ell}W(v_i)\leq1.
    \]
    The single-vertex path \(P=(v)\) is explicitly included. We say that $v$ is a \emph{global breakpoint} if it is an $r$-breakpoint for all $r$. 
\end{definition}

With this, we consider the following block partition.

\begin{definition}[Valid block partition]\label{def:validblockpartion}
    Given a graph $G$ and parameters $d, \eta, A, C$, fix also a diameter constant $D=D(\eta,A,C)>0$, independent of $d$ and $n$. Let $\widehat{S}$ be the set of global breakpoints in $G$. Choose a subset $S \subseteq \widehat{S}$ and let $\mathcal{B}$ the set of connected components of $G$ after removing all vertices in $S$. We refer to every component in $\mathcal{B}$ as a \emph{block} and to the vertices in $S$ as \emph{singletons}. This block partition resulting from $S$ is called \emph{valid} if the following conditions are met.
    \begin{enumerate}
        \setlength{\itemsep}{0cm}
        \item Every block $B \in \mathcal{B}$ induces a tree with at most one extra edge. If $B$ is a tree, we call it a \emph{tree block}, otherwise we call it \emph{unicyclic}.
        \item For every $v \in S$ and block $B \in \mathcal{B}$, there is at most one edge between $v$ and $B$. In other words, every vertex in the boundary $\partial B = \{v \in S \mid \{u, v\} \in E \text{ for some } u \in B\}$ has exactly one neighbor in $B$. We call vertices in $B$ that are adjacent to (one or multiple) vertices in $\partial B$ \emph{ports}.
        \item The diameter of every block $B$ is at most $D
            \log(n)/d^{2\eta}.
        $
    \end{enumerate}
\end{definition}

We show in \Cref{sec:blocks} that $G(n, d/n)$ admits such a valid partition, with a suitable diameter constant $D$, for every fixed $A\geq0$, $C>0$, and $0<\eta<1/2$, provided that $d$ is sufficiently large.

\subsection{Preliminaries on block dynamics and its generator}

Assume that $G$ admits a valid block partition as given by \Cref{def:validblockpartion}. We consider block dynamics with generator
$$
    L = \sum_{s \in S} L_s + \sum_{B \in \mathcal{B}} L_B
$$
where $$
    L_sf (\sigmabf) = f( \sigmabf ) - \Expectednop{ f \mid \sigmabf_{ V - s} } \text{ and }  L_Bf (\sigmabf) = f(\sigmabf) - \Expectednop{f \mid \sigmabf_{V \setminus B}}
$$ are the generators corresponding to re-sampling the color at vertex $s$, and the entire block $B$, respectively. The goal of this section is to establish the following theorem.

\begin{restatable}[Spectral gap for coloring block dynamics]{theorem}{spectralgapblockcoloring}
\label{thm:spectralgapblockcoloring}
    Fix any $0 \le \varepsilon \le 1/2$ and $0<\eta<\varepsilon$. Then for sufficiently large $A, C$, $d$ and with high probability as
    $n\to\infty$ over $G\sim G(n,d/n)$, every partition satisfying
    the block assumptions from \Cref{def:validblockpartion} results in
    \[
        \gap{L}=\Omega_{\varepsilon,\eta}
            \big( d^{\varepsilon-1/2} \big)
    \]
    whenever the number of colors $q$ satisfies $q\geq d+d^{1/2+\varepsilon}$.
\end{restatable}

To show the above theorem, our goal is to show that $L^2 - \kappa L \succeq 0$ for some $\kappa > 0$, implying that $\gap{L} \ge \kappa$. It is actually convenient to work with the so called \emph{accelerated generator}
$$
    \Ltil = \sum_{s \in S} H_s + q\sum_{B \in \mathcal{B}} L_B, \text{ where } H_s = a_s L_s
$$ and $a_s(\sigmabf)$ is the number of available colors at $s$ given that the colors outside of $s$ are given by $\sigmabf_{V - s}$.
We will show that $\Ltil^2 - \kappa q \Ltil \succeq 0$, which is sufficient as shown by the following lemma.
\begin{lemma}\label{lem:acceleratedgap}
    If $\Ltil^2 - \kappa q \Ltil \succeq 0$ and the kernel of $L$ consists only of constants, then $\gap{L} \ge \kappa$.
\end{lemma}
\begin{proof}
    Since $1\leq a_s\leq q$, the Dirichlet forms satisfy
    $L\preceq\Ltil\preceq qL$, so the two operators have the same
    kernel $\mathcal K$. The spectral theorem and the assumed
    inequality imply that every nonzero eigenvalue of $\Ltil$ is
    at least $\kappa q$. Hence, for $f\perp\mathcal K$,
    \[
        \langle f,Lf\rangle
            \geq q^{-1}\langle f,\Ltil f\rangle
            \geq\kappa\|f\|_2^2,
    \]
    which proves the claim by the variational characterization of the gap.
\end{proof}

Expanding the Bochner operator $B_\kappa = \Ltil^2 - \kappa q \Ltil$, we obtain (using $L_B^2 = L_B$ and $H_s^2 = a_sH_s$),
\begin{align*}
    \langle f, B_\kappa f\rangle &= \langle f, (\Ltil^2 - \kappa q \Ltil) f\rangle = \sum_{s \in S} \langle f,  (a_s - \kappa q)H_sf \rangle + q^2(1-\kappa)\sum_{B \in \mathcal{B}} \langle f, L_B f \rangle \\
    & \hspace{1cm} + 2 \sum_{s < t \in S} \langle H_sf, H_t f \rangle + 2q \sum_{s \in S, B \in \mathcal{B}} \langle H_sf, L_B f\rangle + 2q^2 \sum_{\{B, C\} \subseteq \mathcal{B}} \langle f , L_BL_C f \rangle. 
\end{align*}
For showing $B_\kappa \succeq 0$, most of the above terms are harmless and can be ignored. Concretely, non-adjacent singleton-singleton pairs $s, t$, non-adjacent singleton-block pairs $s, B$, and all block-block pairs $B, C$ can be ignored because the corresponding operators commute due to conditional independence. For block-block pairs this is the case because none of them are adjacent by construction. Hence, the identity simplifies to 
\begin{align}\label{eq:bochner}
    \langle f, B_\kappa f \rangle &\ge  \sum_{s \in S} \langle f,  (a_s - \kappa q)H_sf \rangle + q^2(1-\kappa)\sum_{B \in \mathcal{B}} \langle f, L_B f \rangle \nonumber \\
    & \hspace{2cm} + 2 \sum_{\{s, t\} \in E(G[S])} \langle H_sf, H_t f\rangle + 2q \sum_{B \in \mathcal{B}, s \in \partial B} \langle H_s f, L_B f \rangle,
\end{align} where $\partial B$ denotes the set of singleton vertices bordering $B$. 

\subsection{Handling the contribution of singletons}
\label{sec:singletonstuff}

\newcommand{\grad}[2]{\nabla_{\!\!#1 \leftarrow #2}f}
\newcommand{\gradvecff}[1]{{\nabla}^{#1}_{\!\!\sigmabf} }
\newcommand{\recolor}[2]{\sigmabf_{#1 \leftarrow #2}}
\newcommand{\availablecolors}[2]{A_{#1}(#2)}

To show that $B_\kappa$ is PSD, our first step consists in controlling only the contribution of the vertices in $S$. To this end, we closely follow the ideas presented in \cite{haqi2026sampling} and express the generator $H_s$ conveniently in terms of \emph{recoloring gradients}. To this end, denote by $\availablecolors{s}{\sigmabf}$ the set of available colors at vertex $s$ assuming the colors in its neighborhood are given by the configuration $\sigmabf$, and define the recoloring gradient $\gradvec{c} \in \mathbb{R}^{|S|}$ as 
$$
    \gradvec{c}(s) = \begin{cases}
        f(\recolor{s}{c}) - f(\sigmabf) &\text{ if } c \in \availablecolors{s}{\sigmabf}\\
        0 & \text{otherwise.}
    \end{cases}
$$
With this, we can conveniently re-write 
\begin{align}\label{eq:Hgrads}
    H_sf(\sigmabf) = a_s(\sigmabf)L_sf(\sigmabf) = -\sum_{c \in [q]} \gradvec{c}(s),
\end{align}
and setting $a \coloneqq q - \lowdegree$, we immediately get 
$$ 
    \sum_{s \in S} \langle f,  (a_s - \kappa q)H_sf \rangle \ge \frac{a- \kappa q}{2}\sum_{c \in [q]} \Expected{ \| \gradvec{c} \|_2^2 }
$$
and 
\begin{equation}\label{eq:singletoncrosstermsgradient}
    \langle H_sf, H_t f\rangle = \sum_{(c_1, c_2) \in [q]^2} \Expectednop{ \gradvec{c_1}(s) \gradvec{c_2}(t) }.
\end{equation}
The remarkable fact that allows us to compare these two terms is that the sum over all $c_1 \neq c_2$ is non-negative.
\begin{lemma}[Eliminating complete squares]\label{lem:completesquares}
    For any $\{s, t\} \in E(G[S])$, we have
    $$
       \langle H_sf, H_t f\rangle \ge \sum_{c \in [q]} \Expectednop{ \gradvec{c}(s) \gradvec{c}(t) }.
    $$
\end{lemma}
\begin{proof}
    Fix a feasible coloring $\taubf$ outside of $s, t$. For colors $c_1, c_2, c_1', c_2'$, let the colorings $\sigmabf_{00}, \sigmabf_{10}, \sigmabf_{01}, \sigmabf_{11}$ be defined such that the color of $s$ is $c_1'$ whenever the first index is zero and $c_1$ otherwise while the color of $t$ is $c_2'$ whenever the second index is zero and $c_2$ otherwise. All other colors in $\sigmabf_{00}, \sigmabf_{10}, \sigmabf_{01}, \sigmabf_{11}$ are the same as in $\taubf$. Now, define $f_{ij} = f(\sigmabf_{ij})$ and notice that for given $c_1, c_2$, 
    \begin{align*}
        &\sum_{ \substack{(c_1, c_2) \in [q]^2 \\ c_1 \neq c_2}} \Expectednop{ \gradvec{c_1}(s) \gradvec{c_2}(t) } = \sum_{\taubf} \sum_{ (\{c_1, c_1'\}, \{c_2, c_2'\}) } \varphi(\tau, c_1, c_2, c_1', c_2')
    \end{align*}
    where the sum over $(\{c_1, c_1'\}, \{c_2, c_2'\})$ is over all ordered pairs of unordered pairs of colors (where $c_1 = c_1'$ and $c_2 = c_2'$ are explicitly allowed) and \begin{align*}
        \varphi(\tau, c_1, c_2, c_1', c_2') &\coloneqq \mu(\sigmabf_{00}) (f_{10} - f_{00})(f_{01} - f_{00}) + \mu(\sigmabf_{10}) (f_{00} - f_{10})(f_{11} - f_{10}) \\ 
        &\hspace{1cm} + \mu(\sigmabf_{01}) (f_{11} - f_{01})(f_{00} - f_{01}) + \mu(\sigmabf_{11}) (f_{01} - f_{11})(f_{10} - f_{11}).
    \end{align*} Since $\mu(\sigmabf_{00}) = \mu(\sigmabf_{10}) = \mu(\sigmabf_{01}) = \mu(\sigmabf_{11})$, we get 
    $$
        \varphi(\tau, c_1, c_2, c_1', c_2') = \mu(\sigmabf_{00}) (f_{00} - f_{10} - f_{01} + f_{11})^2 \geq 0
    $$ and the lemma follows.
    All sums in this proof are restricted to feasible exterior colorings and color pairs with four proper corners $\sigmabf_{00}, \sigmabf_{01}, \sigmabf_{10}, \sigmabf_{11}$, which excludes only terms with zero weight or zero recoloring gradient.
\end{proof}
Aggregating all terms, we now get the following useful corollary, which is essentially a re-statement of the Bochner-type identity in \cite[Lemma 4]{haqi2026sampling}.
\begin{corollary}\label{cor:lowdegreebochner}
    With $a \coloneqq q - \lowdegree$, we get 
    \begin{align*}
        &\sum_{s \in S} \langle f,  (a_s - \kappa q)H_sf \rangle + 2\sum_{\{s, t\} \in E(G[S])}\langle H_sf, H_t f\rangle \\
        &\hspace{2cm}\ge \frac{a- \kappa q}{2}\sum_{c \in [q]} \Expected{ \| \gradvec{c} \|_2^2 } + \sum_{c \in [q]} \Expected{ (\gradvec{c})^{\!\top} A_{G[S]} (\gradvec{c}) }.
    \end{align*}
\end{corollary}

\subsection{Handling the contribution of blocks}\label{sec:blockcontribution}

The more challenging part is to handle the block generators $L_B$ and their interaction with singleton vertices in the global Bochner-identity. As before, we would like to compare the contribution of the cross-terms $\langle H_sf, L_B f \rangle$ against the (non-negative) diagonal terms $\langle f, H_s f\rangle$ and $\langle f, L_B f\rangle$. 

Fortunately, it turns out that these diagonal and cross terms can be expressed in terms of similar gradients as is the case for singleton-singleton terms (cf. \Cref{cor:lowdegreebochner}). In the following, we formally carry out this work.

\subsubsection{Projection to port gradients}\label{sec:blockassumptions}

Since blocks only interact via ports, it is convenient to define the following projection of $f$ to the port $\rho$, i.e., to average out all its dependence on the colors in $B - \rho$. 

\newcommand{\portprobs}[3]{p_{#1}^{#3}(#2)}

\begin{definition}[Port projection operator]
    Given a block $B$ and a port $\rho$ in $B$, define the operator $\portproj{\rho}$ as 
    $$
        \portproj{\rho} f(\sigmabf) = \Expectednop{ f \mid \sigmabf_{(V \setminus B) +\rho} }.
    $$
\end{definition}

Then, all interaction terms of a block only depend on gradients on $\portproj{\rho}f$ for some port $\varrho$. More concretely, our goal is to express the remaining terms in \eqref{eq:bochner} in terms of the following gradients.
\begin{definition}[Port gradients]
    Given a block $B$, a port $\rho$ in $B$, and a coloring $\sigmabf$, let 
    $$
        \portprob{\rho}{c} = \Pr{\sigmabf_\rho = c \mid \sigmabf_{V \setminus B} } 
    $$
    be the cavity probability of assigning color $c$ to $\rho$ assuming a block boundary given by $\sigmabf_{V \setminus B}$. Define the \emph{port gradient} $\portgrad{c} \in \mathbb{R}^{|P_B|}$ as 
    $$   
        \portgrad{c}(\rho) = q\portprob{\rho}{c} \left( \portproj{\rho}f(\recolor{\rho}{c}) - \portproj{\rho}f(\sigmabf)\right).
    $$
    Note that, even if $\portprob{\rho}{c}>0$, the literal recoloring $\recolor{\rho}{c}$ may be improper. However, the projected value $\portproj{\rho}f(\recolor{\rho}{c})$ still averages $f$ over feasible colorings with exterior $\sigmabf_{V\setminus B}$ and color $c$ at $\rho$.
\end{definition}

\subsubsection{Relating cross terms to port gradients}

Our first simple lemma that relates the cross terms from \eqref{eq:bochner} to the port projection $\portproj{\rho}f$ is the following.
\begin{lemma}\label{lem:portproj}
    Given any block $B$, any $s \in \partial B$, and abbreviating $\rho = \port{s}$, we have 
    $$
        \langle H_s f, L_B f \rangle \ge \langle H_s \portproj{\rho}f, L_B \portproj{\rho}f \rangle .
    $$
\end{lemma}
\begin{proof}
    Write $P_Bf=\Expectednop{f\mid\sigmabf_{V\setminus B}}$, so
    $L_B=I-P_B$.
    The projection $\portproj{\rho}$ commutes with $H_s$ because it
    leaves the color at $\rho$ unchanged and $s$ is conditionally
    independent of $B\setminus\{\rho\}$ given the colors on
    $(V\setminus(B+s))+\rho$. We further get $P_B\portproj{\rho}=\portproj{\rho}P_B=P_B$ because $P_B$ projects to a subset of the set $\portproj{\rho}$ projects to. Hence,
    $\portproj{\rho}$ commutes with $L_B=I-P_B$ as well. Moreover, since $P_B\portproj{\rho}f=P_Bf$, we have
    \begin{align*}
        L_B(I - \portproj{\rho})f &= (I - \portproj{\rho})f - P_B(I - \portproj{\rho})f \\ 
        &= (I - \portproj{\rho})f - P_Bf + P_Bf = (I - \portproj{\rho})f.    
    \end{align*}
    After setting $f = \portproj{\rho}f  + (I - \portproj{\rho})f$ and expanding $\langle H_s f, L_B f \rangle$, we therefore get 
    \begin{align*}
        \langle H_s f,L_Bf\rangle
          &=\langle H_s\portproj{\rho}f,L_B\portproj{\rho}f\rangle
             +\langle H_s(I-\portproj{\rho})f,
                       L_B(I-\portproj{\rho})f\rangle\\
          &\hspace{1cm}+\langle H_s\portproj{\rho}f,
                                 L_B(I-\portproj{\rho})f\rangle
             +\langle H_s(I-\portproj{\rho})f,
                                 L_B\portproj{\rho}f\rangle.
    \end{align*}
    Since $\portproj{\rho}$ is an orthogonal projection, we get $\langle \portproj{\rho} f, (I - \portproj{\rho})g \rangle = 0$ for all $f,g$ and by our commutation properties together with $L_B(I - \portproj{\rho})f = (I - \portproj{\rho})f$ the cross terms in the second line above vanish.

    $$
        \langle H_s f, L_B f \rangle =  \langle H_s \portproj{\rho}f, L_B \portproj{\rho}f \rangle + \langle H_s (I- \portproj{\rho})f, L_B(I - \portproj{\rho})f \rangle.
    $$ 
    Since $H_s$ is PSD and $L_B(I - \portproj{\rho})f = (I - \portproj{\rho})f$, we get $\langle H_s (I- \portproj{\rho})f, L_B(I - \portproj{\rho})f \rangle = \langle H_s (I- \portproj{\rho})f, (I - \portproj{\rho})f \rangle \ge 0$ and the lemma follows.
\end{proof}
With this, the cross terms $\langle H_s \portproj{\rho}f, L_B \portproj{\rho}f \rangle$ now behave similar to the singleton cross pairs and can be related to port gradients analogously to \Cref{eq:singletoncrosstermsgradient}. 
\begin{lemma}[Expressing cross terms via port gradients, analog to \eqref{eq:singletoncrosstermsgradient}]\label{lem:blockcrosstermsgradient}
    Given any block $B$, any $s \in \partial B$, and abbreviating $\rho = \port{s}$, we have 
    $$
        q\langle H_s \portproj{\rho}f, L_B \portproj{\rho}f \rangle  = \sum_{c_1, c_2 \in [q]} \Expected{ \gradvec{c_1}(s) \portgrad{c_2}(\rho) }.
    $$
\end{lemma}
\begin{proof}
     We have 
    \begin{align*}
        L_B\portproj{\rho}f(\sigmabf) &= \portproj{\rho}f(\sigmabf) - \Expectednop{ \portproj{\rho}f \mid \sigmabf_{V \setminus B}} \\
        &= \portproj{\rho}f(\sigmabf) - \sum_{c \in [q]} \portprob{\rho}{c} \portproj{\rho}f(\recolor{\rho}{c}) = -\frac{1}{q} \sum_{c \in [q]} \portgrad{c}(\rho).
    \end{align*}
    Moreover, using \eqref{eq:Hgrads}, the definition of $\portproj{\rho}$, and conditional independence of $s$ and $B - \rho$, we get
    $$
        H_s \portproj{\rho} f (\sigmabf) = \sum_{c \in \availablecolors{s}{\sigmabf}} (\portproj{\rho}f(\sigmabf) - \portproj{\rho}f(\recolor{s}{c})) = -\sum_{c \in [q]} \Expected{ \gradvecff{c} f(s) \mid \sigmabf_{(V \setminus B) + \rho}}.
    $$
    Taking the expectation then yields the lemma.
\end{proof}
As before, we would like to get rid of terms corresponding to $c_1 \neq c_2$ like in \Cref{lem:completesquares}. While a complete elimination is not possible anymore, an approximate version of this statement still holds, with an error that can be expressed in terms of the conditional variance of $\portproj{\rho}$. The reason for this error is that once we try to eliminate complete squares as in \Cref{lem:completesquares}, the weight of the four configurations $\sigmabf_{00}, \sigmabf_{01}, \sigmabf_{10}, \sigmabf_{11}$ is not identical anymore as the distribution of the color at a port $\rho$ depends on the color at the boundary vertex $s$. However, this fluctuation can be estimated in terms of the expected variance of $\portproj{\rho}f$, which later can be compared to $\langle f, L_B  f \rangle$. 
\begin{lemma}[Approximately eliminating complete squares, analog to \Cref{lem:completesquares}]\label{lem:approxcompletesquares}
    Given any block $B$, any $s \in \partial B$, and abbreviating $\rho = \port{s}$, we have
    $$
        q\langle H_s f, L_B f \rangle \ge \sum_{c \in [q]} \Expected{ \gradvec{c}(s) \portgrad{c}(\rho) } - \frac{4q^2}{a^2} \Expected{ \Varnop{ \portproj{\rho}f \mid \sigmabf_{V \setminus B} } }.
    $$
\end{lemma}
\begin{proof}

    \newcommand{\nutau}{\nu_{\taubf}}
    \newcommand{\pitau}{\pi_{\taubf}}
    \newcommand{\pitaucond}[1]{\pi^{#1}_{\taubf}}

    For any feasible coloring $\taubf$ on $V \setminus (B + s)$, we define $
        \pitau$ as the joint distribution of the colors at $s$ and $\rho$ conditioned on the colors given by $\taubf$. The value $\pitau(a, b)$ is the probability that $\sigmabf_s = a$ and $\sigmabf_\rho = b$ under this measure. We further define 
    $
        \pitaucond{\rho} (b \mid a) 
    $ as the probability that $\sigmabf_\rho = b$ conditioned on $\taubf$ and $\sigmabf_s = a$.
    
    For colors $c_1,c_1',c_2,c_2'$, denote by $\sigmabf_{ij}$ with $i, j \in \{0,1\}$ the coloring on $(V\setminus B) + \rho$ obtained by setting the color at $s$ to $c_1$ if $i = 0$ and $c_1'$ otherwise, while setting   the color at $\rho$ to $c_2$ if $j = 0$ and to $c_2'$ otherwise.
    Now, just like in \Cref{lem:completesquares}, we re-write the sum using unordered pairs of colors:
    \begin{align*}
        &\frac{1}{q}\sum_{ \substack{(c_1, c_2) \in [q]^2 \\ c_1 \neq c_2}} \Expectednop{ \gradvec{c_1}(s) \portgrad{c_2}(\rho) } = \sum_{\taubf} \Pr{\sigmabf_{V \setminus (B + s)} = \taubf} \sum_{ (\{c_1, c_1'\}, \{c_2, c_2'\}) } \varphi(\taubf, c_1, c_2, c_1', c_2')
    \end{align*}
    with
    \begin{align*}
        \varphi(\taubf, c_1, c_2, c_1', c_2') &= \pi_{00} (f_{10} - f_{00})(f_{01} - f_{00}) + \pi_{10} (f_{00} - f_{10})(f_{11} - f_{10}) \\
        &\hspace{.7cm} + \pi_{01} (f_{11} - f_{01})(f_{00} - f_{01}) + \pi_{11} (f_{01} - f_{11})(f_{10} - f_{11})
    \end{align*} where $f_{ij} = \portproj{\rho}f(\sigmabf_{ij})$ and 
    \begin{align*}
        \pi_{00} &= \pitaucond{\rho}(c_2' \mid c_1) \pitau(c_1,c_2), \hspace{.2cm} \pi_{10} = \pitaucond{\rho}(c_2' \mid c_1')\pitau(c_1',c_2), \hspace{.2cm}\\
        \pi_{01} &= \pitaucond{\rho}(c_2 \mid c_1)\pitau(c_1,c_2'), \hspace{.2cm} \pi_{11} = \pitaucond{\rho}(c_2 \mid c_1')\pitau(c_1',c_2').
    \end{align*}
    We remark that every sum in this proof over exterior colorings is over feasible $\taubf$. Moreover, every sum over $(\{c_1,c_1'\},\{c_2,c_2'\})$ is restricted to color pairs for which all four partial colorings $\sigmabf_{ij}$ extend to proper colorings of $G$ like in the proof of \Cref{lem:completesquares}. These restrictions exclude only zero terms.

    Now, it is easy to see that $\pi_{00} = \pi_{01}$ and $\pi_{10} = \pi_{11}$. We therefore write $\pi_0 \coloneqq \pi_{00} = \pi_{01}$ and $\pi_1 \coloneqq \pi_{10} = \pi_{11}$.
    Then, setting $\alpha = f_{00} - f_{01} - f_{10} + f_{11}$ and $\beta = f_{00} - f_{01}$, we have
    \begin{align*}
        \varphi \coloneqq \varphi(\taubf, c_1, c_2, c_1', c_2') &= \pi_1\alpha^2 + (\pi_0 - \pi_1) \alpha \beta \\
        &=
\pi_1\left(
\alpha+\frac{\pi_0-\pi_1}{2\pi_1}\beta
\right)^2
-
\frac{(\pi_0-\pi_1)^2}{4\pi_1}\beta^2.
    \end{align*}

    It remains to bound the negative term. The breakpoint
    condition (with $A\geq1$) makes $\rho$ low-degree. Conditioning first
    on all other vertex colors leaves at least $a$ available colors at
    $\rho$, so $\pitaucond{\rho}(c\mid c_1)\leq1/a$ for every feasible
    $c_1$ and every $c$. We may assume $a\geq2$ for sufficiently large $d$.

    Whenever all four $\sigmabf_{ij}$ are feasible, recoloring only $s$
    between $c_1$ and $c_1'$ preserves properness and the weight of each
    coloring if $\rho$ avoids both colors. Thus
    $\pitau(c_1,c_2)=\pitau(c_1',c_2)$, and the two conditional port laws
    agree after conditioning on this avoidance. Consequently,
    \[
        \frac{\pi_0}{\pi_1}
          =\frac{\pitaucond{\rho}(c_2'\mid c_1)}
                 {\pitaucond{\rho}(c_2'\mid c_1')}
          =\frac{1-\pitaucond{\rho}(c_1'\mid c_1)}
                 {1-\pitaucond{\rho}(c_1\mid c_1')}.
    \]
    The marginal bound therefore gives
    \[
        \left|\frac{\pi_0}{\pi_1}-1\right|
            \leq\frac1{a-1}\leq\frac2a,
        \qquad \frac{\pi_0}{\pi_1}\geq1-\frac1a\geq\frac12.
    \]
    Discarding the first, nonnegative term in the displayed identity
    for $\varphi$, we conclude that
    \[
        \varphi\geq-\frac{(\pi_0-\pi_1)^2}{4\pi_1}\beta^2
                   \geq-\frac{2\pi_0}{a^2}\beta^2.
    \]

    To sum this error, write
    $\pitaucond{s}(c_1)=\sum_{c\in[q]}\pitau(c_1,c)$ and let
    $\mathcal A=\{c_1\in[q]:\pitaucond{s}(c_1)>0\}$.
    By definition,
    \[
        \pi_0=\pitaucond{s}(c_1)
                 \pitaucond{\rho}(c_2\mid c_1)
                 \pitaucond{\rho}(c_2'\mid c_1).
    \]
    For fixed $c_1\in\mathcal A$, summing over ordered port colors gives
    \[
        \sum_{\substack{c_2,c_2'\in[q]\\
                        \pitau(c_1,c_2)\pitau(c_1,c_2')>0}}\pi_0\beta^2
          =2\pitaucond{s}(c_1)
              \Varsub{\pitaucond{\rho}(\cdot\mid c_1)}{\portproj{\rho}f},
    \]
    since $\beta=f_{00}-f_{01}$ is the difference of the values at two
    independent port colors with $s$ fixed to $c_1$. It does not depend
    on $c_1'$, for which there are at most $q$ choices. Enlarging the
    negative error sum from the unordered color pairs above to all
    these ordered choices therefore gives
    \begin{align*}
        &\sum_{ (\{c_1, c_1'\}, \{c_2, c_2'\}) }\varphi(\taubf, c_1, c_2, c_1', c_2') \ge
          \\
          &\hspace{2cm}-\frac{4q}{a^2}\sum_{c_1\in\mathcal A}\pitaucond{s}(c_1)
               \Varsub{\pitaucond{\rho}(\cdot\mid c_1)}{\portproj{\rho}f} =-\frac{4q}{a^2}\Expected{
              \Varnop{\portproj{\rho}f\mid\sigmabf_{V\setminus B}}
              \mid\sigmabf_{V\setminus(B+s)}=\taubf},
    \end{align*}
    using conditional expectation over the singleton color.
    By the unique-port property, projecting a singleton gradient onto
    $\portproj{\rho}$ gives the corresponding gradient of
    $\portproj{\rho}f$; since the port gradient is already measurable
    with respect to $\sigmabf_{(V\setminus B)+\rho}$, this projection
    preserves the expected products in the initial decomposition.
    Averaging over $\taubf$ and multiplying by $q$ bounds its
    off-diagonal contribution by the claimed error. Together with
    \Cref{lem:portproj,lem:blockcrosstermsgradient}, this proves the lemma.
\end{proof}

\subsubsection{Block elimination}\label{sec:blockelim}

Having reduced the singleton-block interaction terms in \eqref{eq:bochner} to port gradients, it remains to also express the diagonal terms in terms of port gradients which will counterbalance the potentially negative terms remaining after \Cref{lem:blockcrosstermsgradient}. To this end, we first relate squared port gradients to the expected conditional variance of $\portproj{\rho}$.

\begin{lemma}\label{lem:portgradientstovariance}
    For any block $B$ with port $\rho \in P_B$,
    \begin{align*}
        \sum_{c \in [q]}\Expected{ \portgrad{c}(\rho)^2 } \le \frac{2q^2}{a} \Expected{ \Varnop{ \portproj{\rho}f \mid \sigmabf_{V \setminus B} } }.
    \end{align*}
\end{lemma} 
\begin{proof}
    Given a coloring $\tau$ on $V \setminus B$, define $f_\rho^{\taubf} (c) = \Expectednop{ f(\sigmabf) \mid \sigmabf_{V\setminus B } = \taubf, \sigmabf_\rho = c}$.
    Then, 
    \begin{align*}
        \sum_{c \in [q]}\Expected{ \portgrad{c}(\rho)^2 \mid \sigmabf_{V \setminus B} = \taubf } &= q^2 \sum_{(c_1, c_2) \in [q]^2} \portprob{\rho}{c_2}\portprob{\rho}{c_1}^2 \left(  f_\rho^{\taubf} (c_1) - f_\rho^{\taubf} (c_2) \right)^2 
    \end{align*}
    Because $\portprob{\rho}{c} \le 1/a$ for all ports and colors, we get
    \begin{align*}
           \sum_{c \in [q]}\Expected{ \portgrad{c}(\rho)^2 \mid \sigmabf_{V \setminus B} = \taubf } &\le \frac{q^2}{a} \sum_{(c_1, c_2) \in [q]^2} \portprob{\rho}{c_2}\portprob{\rho}{c_1} \ ( f_\rho^{\taubf} (c_1) - f_\rho^{\taubf} (c_2) )^2\\
        &= \frac{2q^2}{a} \Varnop{ \portproj{\rho}f \mid \sigmabf_{V \setminus B } = \taubf }.
    \end{align*} Taking expectation over $\taubf$ yields the lemma.
\end{proof}

Then, we relate the RHS above summed over all ports of $B$ to $\langle f , L_B f\rangle$. This is done by establishing approximate conditional independence of the port spins and manifests in a bound on the largest eigenvalue of a covariance operator $\covop{B}$ defined as follows. 

\newcommand{\centeredspace}[1]{\mathcal{S}_{#1, \taubf}}
\newcommand{\productspace}[1]{\mathcal{H}_{#1, \taubf}}
\newcommand{\nutau}{\nu_{\taubf}}

\begin{definition}[Block covariance operator $\covoppure{B}$]\label{def:covop}
    Given a  block $B$ and a coloring $\taubf$ on $V \setminus B$, let $\nutau$ be the measure of $\sigmabf_B$ conditional on $\taubf$.
    For any port $\rho \in P_B$, let $\mathcal{S}_\rho$ be the space of functions $f: [q] \rightarrow \mathbb{R}$ such that $\Expectednop{ f(\sigmabf_\rho) \mid \sigmabf_{V \setminus B} = \taubf} = 0$.
    We restrict $\mathcal S_\rho$ to functions on
    $\{c\in[q]:\nutau(\sigmabf_\rho=c)>0\}$, so its
    $L^2(\nutau)$ inner product is nondegenerate.
    Consider the Hilbert space
    $$
        \productspace{B} = \bigoplus_{\rho \in P_B} \centeredspace{\rho},
    $$ where $\oplus$ denotes the direct product. Define the operator $\covop{B}: \productspace{B} \rightarrow \productspace{B}$ as   
    $$
      (\covop{B} g)_{\rho} (c) \coloneqq \Expectedsub{\nutau}{ {\sum}_{r \in P_B} \ g_{r} (\sigmabf_{r}) \mid \sigmabf_\rho = c }.
    $$
    Equivalently, we may view $\covop{B}$ as an operator-valued $|P_B| \times |P_B|$ matrix where $(\covop{B})_{\rho_1, \rho_2}$ maps from $\centeredspace{\rho_2}$ to $\centeredspace{\rho_1}$ and 
    $$
        (\covop{B})_{\rho_1, \rho_2}h(c) = \Expected{h(\sigmabf_{\rho_2}) \mid \sigmabf_{\rho_1} = c }.
    $$
\end{definition}
The operator $\covop{B}$ essentially measures the covariance between different ports in the sense that 
$$
    \langle f, \covop{B}f \rangle_{\productspace{B}} = {\sum}_{(\rho_1, \rho_2) \in P_B^2} \ \Expectedsub{\nutau}{ f_{\rho_1}(\sigmabf_{\rho_1}) f_{\rho_2}(\sigmabf_{\rho_2}) }
$$
where $f$ is a vector of $|P_B|$ centered functions that only depend on the spin at a single port each. Because of that, the largest eigenvalue of $\covop{B}$ in $\productspace{B}$ allows us to compare the right hand side of \Cref{lem:portgradientstovariance} to $\langle f, L_B f\rangle$. 

\begin{lemma}\label{lem:variancetoeigenvalue}
    For any block $B$,
    \begin{align*}
        \sum_{\rho \in P_B}\Expected{ \Varnop{ \portproj{\rho}f \mid \sigmabf_{V \setminus B} } } \le \sup_{\taubf}\lambdamax ( \covop{B}) \cdot \langle f, L_B f \rangle.
    \end{align*}
\end{lemma}
\begin{proof}
    \newcommand{\utau}[1]{u_{\taubf, #1}}
    \newcommand{\utauvec}{u_{\taubf}}
    Given a coloring $\taubf$ on $V \setminus B$, we denote by $\nutau$ the measure of $\sigmabf_B$ conditional on $\taubf$ and for every $\rho \in P_B$, we define the centered version of $\portproj{\rho}f$ as 
    $$
        \utau{\rho}(c) \coloneqq \portproj{\rho}f( \sigmabf_{\tau, c} ) - \Expectedsub{\nutau}{ \portproj{\rho}f} = \Expectedsub{\nutau}{L_B f \mid \sigmabf_\rho= c}
    $$ where $\sigmabf_{\tau, c}$ is the coloring that agrees with $\tau$ on $V \setminus B$ and sets all colors in $B$ to $c$\footnote{Since $\portproj{\rho}f$ only depends on the spin at $(V\setminus B) + \rho$, we could equivalently choose any coloring that agrees with $\tau$ on $V \setminus B$ and sets the color of $\rho$ to $c$.}. We further define 
    $$
        A_{\taubf} \coloneqq {\sum}_{\rho \in P_B} \Expectedsubnop{\nutau}{ \utau{\rho}(\sigmabf_\rho)^2 } =  {\sum}_{\rho \in P_B}\Varnop{ \portproj{\rho}f \mid \sigmabf_{V \setminus B} = \taubf } 
    $$ and the function
    $$
        U_{\taubf}(\sigmabf) \coloneqq \sum_{\rho \in P_B} \utau{\rho}(\sigmabf_\rho). 
    $$ To simplify notation, we omit the dependence on $\taubf$ in the following and simply write $A, U, u,$ and $\covoppure{B}$ instead of $A_{\taubf}, U_{\taubf}, \utauvec, \covop{B}$. Moreover, for the rest of this proof, inner products and norms refer to the Hilbert space $L^2(\nutau)$ unless stated otherwise. Now, $A$ and $U$ are related through $\covoppure{B}$ as 
    $$
         \Expected{ U(\sigmabf)^2 } = \| U \|^2 = \langle u, \covoppure{B} u \rangle_{\productspace{B}} \le \lambdamax( \covoppure{B} ) \| u \|_{\productspace{B}}^2 = \lambdamax( \covoppure{B} ) A.
    $$
    Moreover, we have
    \begin{align*}
        \langle L_Bf, U \rangle &= {\sum}_{\rho \in P_B} \Expectedsub{\nutau}{ u_\rho(\sigmabf_\rho)L_Bf(\sigmabf) }\\
        &= {\sum}_{\rho \in P_B}  \Expectedsub{\nutau}{ u_\rho(\sigmabf_\rho)\Expectednop{L_Bf \mid \sigmabf_{\rho}} } = {\sum}_{\rho \in P_B}  \Expectedsubnop{\nutau}{ u_\rho(\sigmabf_\rho)^2 } = A.
    \end{align*} Hence, by Cauchy-Schwarz, 
    \begin{align*}
        A = \langle L_Bf, U \rangle \le \| L_Bf \| \| U \| \le \sqrt{ \langle L_B f, L_B f \rangle \lambdamax( \covoppure{B} ) A },
    \end{align*} and since $\langle L_B f, L_B f \rangle = \langle f, L_B^2 f \rangle = \langle f, L_B f \rangle$, the lemma follows.
\end{proof}

\subsubsection{Assembling diagonal and cross terms}

Combining all the lemmas, we get the following.
\begin{corollary}[Analog of \Cref{cor:lowdegreebochner}]\label{cor:hubbochner}
    Let $\lambda\geq1$ be an upper bound for
    $\sup_{B,\taubf}\lambdamax(\covop{B})$. Assume $A\geq1$, $a\geq2$,
    and $1-\kappa\geq8\lambda\lowdegree/a^2$. Then,
    \begin{align*}
        &q^2(1-\kappa)\sum_{B \in \mathcal{B}} \langle f, L_B f \rangle + 2q \sum_{B \in \mathcal{B}, s \in \partial B} \langle H_s f, L_B f \rangle\\
        &\hspace{1cm} \ge \left(\frac{a(1-\kappa)}{2\lambda}-\frac{4\lowdegree}{a}\right) \sum_{B \in \mathcal{B}} \sum_{c \in [q]} \Expected{ \| \portgrad{c} \|_2^2 } + 2 \sum_{B \in \mathcal{B}} \sum_{s \in \partial B } \sum_{c \in [q]} \Expected{ \gradvec{c}(s) \portgrad{c}(\rho_{B,s}) }
    \end{align*}
\end{corollary}
\begin{proof}
    Fix a block $B$. For each port $\rho\in P_B$, write
    \begin{align*}
        V_\rho\coloneqq
            \Expected{\Var{\portproj{\rho}f\mid\sigmabf_{V\setminus B}}}, \qquad
        m_B(\rho)=|\{s\in\partial B:\rho_{B,s}=\rho\}|.
    \end{align*}
    The breakpoint condition makes every port low-degree, as in the
    proof of \Cref{lem:approxcompletesquares}. Hence
    $m_B(\rho)\leq\deg_G(\rho)\leq\lowdegree$, and
    \Cref{lem:variancetoeigenvalue} gives
    \begin{align*}
        \sum_{s\in\partial B}V_{\rho_{B,s}}
          &=\sum_{\rho\in P_B}m_B(\rho)V_\rho \leq\lowdegree\sum_{\rho\in P_B}V_\rho
            \leq\lowdegree\lambda\langle f,L_Bf\rangle.
    \end{align*}
    Applying \Cref{lem:approxcompletesquares}
    with its error $4q^2/a^2$, we thus obtain
    \begin{align*}
        &q^2(1-\kappa) \langle f, L_B f \rangle + 2q \sum_{s \in \partial B} \langle H_s f, L_B f \rangle \\
        &\hspace{2cm}\geq q^2\left(1-\kappa-
                      \frac{8\lambda\lowdegree}{a^2}\right)
                       \langle f,L_Bf\rangle
            +2\sum_{s\in\partial B}\sum_{c\in[q]}
                 \Expected{\gradvec{c}(s)\portgrad{c}(\rho_{B,s})}.
    \end{align*}
    On the other hand, \Cref{lem:portgradientstovariance,lem:variancetoeigenvalue}
    imply
    \begin{align*}
        \sum_{c\in[q]}\Expected{\|\portgrad{c}\|_2^2}
           &\leq\frac{2q^2}{a}\sum_{\rho\in P_B}V_\rho
             \leq\frac{2q^2\lambda}{a}\langle f,L_Bf\rangle.
    \end{align*}
    By assumption, $1-\kappa-8\lambda\lowdegree/a^2\geq0$, so
    \begin{align*}
        q^2\left(1-\kappa-\frac{8\lambda\lowdegree}{a^2}\right)
                  \langle f,L_Bf\rangle
          &\geq\left(\frac{a(1-\kappa)}{2\lambda}
                         -\frac{4\lowdegree}{a}\right)
                    \sum_{c\in[q]}\Expected{\|\portgrad{c}\|_2^2}.
    \end{align*}
    Combining the last two bounds with the preceding estimate and
    summing over all blocks $B\in\mathcal B$ proves the statement.
\end{proof}

\subsection{Covariance analysis}\label{sec:covarianceanalysis}

We show that $\lambdamax(\covoppure{B}) \le 1 + o_d(1)$. This is based on exploiting decay of correlations within blocks by exploiting their tree-like nature. We start by carrying this out for tree blocks and afterwards extend the argument to unicyclic blocks. 

Throughout this section, we assume that $A$ is sufficiently large depending on the fixed parameters $\varepsilon,\eta$, and then take $d$ sufficiently large. This enables sufficiently strong breakpoint properties.

\subsubsection{Covariance factors along tree paths}

Our first observation is that covariance factors along tree paths inside a block. Specifically, recall from \Cref{def:covop} that given two ports $\rho_1, \rho_2 \in P_B$, the operator $(\covop{B})_{\rho_1, \rho_2}$ measures the covariance between two ports given the coloring $\taubf$ on $V \setminus B$. This operator can be seen as a chain of individual covariance operators along the unique path from $\rho_1$ to $\rho_2$ in $B$. 

\newcommand{\Kuv}[2]{K_{#1 \leftarrow #2}}

\begin{observation}[Covariance operators factor along tree paths]\label{obs:treecov}
    Given a tree block $B$, a coloring $\taubf$ on $V \setminus B$ and two $\rho_1, \rho_2 \in P_B$, let $P = (\rho_1 = v_0, v_1, \ldots, v_{m-1}, v_m = \rho_2)$ be the unique path from $\rho_1$ to $\rho_2$ in $B$. For any $u, v \in B$ denote by $\Kuv{u}{v}$ the operator 
    $$
        \Kuv{u}{v}: \centeredspace{v} \rightarrow \centeredspace{u}, \Kuv{u}{v} f(c) = \Expectedsub{\nutau}{ f(\sigmabf_{v}) \mid \sigmabf_{u} = c }
    $$ and note that $(\covop{B})_{\rho_1, \rho_2} = K_{\rho_1 \leftarrow \rho_2}$ Then, 
    $$
        (\covop{B})_{\rho_1, \rho_2} =  \Kuv{v_0}{v_1}\Kuv{v_1}{v_2} \cdots \Kuv{v_{m-1}}{v_m}.
    $$
\end{observation}
Moreover, we can bound the norm of each $\Kuv{v_{i}}{v_{i+1}}$ as follows.
\begin{lemma}[Covariance along a single edge]\label{lem:edgecovariance}
    Let $u,v$ be adjacent vertices in a tree block $B$. Then,
    $$
        \| \Kuv{u}{v} \| \le \beta_{\{u,v\}} \coloneqq \begin{cases}
            a^{-1} & \text{if } \deg(u), \deg(v) \le \lowdegree \\
            a^{-1/2} & \text{if exactly one of } \deg(u), \deg(v) > \lowdegree \\
            1 & \text{otherwise.}
        \end{cases} 
    $$ with $a = q - \lowdegree$.
\end{lemma}
\begin{proof}
    If either $\centeredspace{u}$ or $\centeredspace{v}$ is
    trivial, the operator norm is zero and the bound is immediate.
    Otherwise, expectations and inner products below are taken
    under $\nutau$, unless a different measure is specified.

    Remove the edge $\{u,v\}$ while keeping the exterior
    coloring fixed, and let $p_i,r_i$ be the resulting marginal
    probabilities of color $i$ at $u,v$, respectively.
    Since $B$ is a tree, removing this edge separates the
    components containing $u$ and $v$, so their colors are
    independent under the edge-deleted measure. Restoring
    the edge amounts to conditioning on distinct colors.
    Consequently, their joint law under $\nutau$ is
    \[
        \pi(i,j)=\frac{p_i r_j\mathds{1}_{i\ne j}}{1-\theta},
        \qquad \theta=\sum_{i\in[q]}p_i r_i.
    \]

    If $u$ is low-degree, then after removing $\{u,v\}$ and conditioning
    on all other neighbors of $u$, at least
    $q-\deg_G(u)+1\geq a+1$ colors are available at $u$.
    Its conditional color is uniform on these colors, and averaging
    over its neighbors therefore gives
    \[
        p_i\leq\frac1{a+1}\qquad\text{for every }i.
    \]
    Similarly, $r_i\leq1/(a+1)$ if $v$ is low-degree.

    To pass from these marginal bounds to an operator norm bound,
    regard $P=\Kuv{u}{v}$ and $Q=\Kuv{v}{u}$ also as conditional
    expectation kernels on all functions, not just centered ones.
    By the expression for $\pi$, their entries are
    \[
        P(i,j)=\frac{r_j\mathds{1}_{j\ne i}}{1-r_i},
        \qquad
        Q(j,i)=\frac{p_i\mathds{1}_{i\ne j}}{1-p_j}.
    \]
    Here and below we restrict the kernels to colors of positive
    actual marginal probability, so all denominators are positive.
    For any Markov kernel $T$, put
    \[
        \delta(T)=\max_{x,y}\frac12\sum_z|T(x,z)-T(y,z)|,
        \qquad
        \operatorname{osc}(h)=\max_x h(x)-\min_x h(x).
    \]
    Subtracting $(\max h+\min h)/2$ from $h$ in the difference of
    two rows shows that
    \[
        \operatorname{osc}(Th)\leq\delta(T)\operatorname{osc}(h).
    \]

    Now, for two distinct supported colors $j,k$ at $v$, a direct
    calculation gives
    \[
        \frac12\sum_i|Q(j,i)-Q(k,i)|
          =\max\left\{\frac{p_j}{1-p_k},
                       \frac{p_k}{1-p_j}\right\}.
    \]
    Indeed, if $p_j\geq p_k$, the row difference is nonnegative
    at every color except $j$, and its negative part has mass
    $p_j/(1-p_k)$; the other case is symmetric.
    Consequently, if $u$ is low-degree, then
    \[
        \delta(Q)\leq
          \frac{1/(a+1)}{1-1/(a+1)}=\frac1a.
    \]
    The analogous bound $\delta(P)\leq1/a$ holds if $v$ is
    low-degree. Without a degree assumption we always have
    $\delta(P),\delta(Q)\leq1$. Thus, writing
    \[
        b_x=\begin{cases}
             a^{-1},&\deg_G(x)\leq\lowdegree,\\
             1,&\deg_G(x)>\lowdegree,
             \end{cases}
    \]
    we get $\delta(Q)\leq b_u$ and $\delta(P)\leq b_v$.

    Finally, $P$ and $Q$ are adjoints with respect to the actual
    marginal laws of $u,v$, since both corresponding inner
    products equal $\Expectedsub{\nutau}{f(\sigmabf_u)g(\sigmabf_v)}$.
    Hence $PQ$ is self-adjoint and positive semidefinite on
    $\centeredspace{u}$, with largest eigenvalue
    $\lambda=\|\Kuv{u}{v}\|^2$.
    Take a nonzero centered eigenfunction $h$ for
    $\lambda$. It is nonconstant on the support, and therefore
    \[
        \lambda\operatorname{osc}(h)
          =\operatorname{osc}(PQh)
          \leq\delta(P)\delta(Q)\operatorname{osc}(h)
          \leq b_u b_v\operatorname{osc}(h).
    \]
    Dividing by $\operatorname{osc}(h)>0$ yields
    \[
        \|\Kuv{u}{v}\|\leq\sqrt{b_u b_v}
          =\beta_{\{u,v\}},
    \] as desired.
\end{proof}

\subsubsection{Tree block eigenvalue bound}

In this section, we prove the following bound on $\lambdamax(\covop{B})$ for tree blocks.
\begin{lemma}[Tree block covariance operator eigenvalue bound]\label{lem:treeeigenvalue}
    Consider any tree block $B$. Then, for every coloring $\taubf$ on $V \setminus B$ and every $\delta > 0$, we have
    $$
        \lambdamax(\covop{B}) \le (1 + \delta).
    $$ provided that $d$ and $A$ are large enough.
\end{lemma}
As our first step, we show that it suffices to study the largest eigenvalue of a simpler scalar matrix instead of the operator valued matrix $\covop{B}$.
\newcommand{\parent}[1]{\mathrm{par}(#1)}
\newcommand{\anc}[1]{\mathrm{anc}(#1)}

\begin{lemma}[Reduction to a scalar covariance matrix]\label{lem:scalarop}
    Let $B$ be a tree block and $\taubf$ be a coloring on $V \setminus B$.
    For two ports $i, j \in P_B$ denote by $P(i, j)$ the unique path from $i$ to $j$ in $B$.
    Define $\scalarcov \in \mathbb{R}^{|P_B| \times |P_B|}$ as 
    $$
        (\scalarcov)_{i, j} = \begin{cases}
            \prod_{e \in P(i,j)}\beta_e & \text{if } i \neq j \\
            1 & \text{otherwise,}
        \end{cases}
    $$ where $\beta$ is as defined in \Cref{lem:edgecovariance}. Then,
    \begin{align*}
        \lambdamax(\covop{B}) \le \lambdamax(\scalarcov). 
    \end{align*}
\end{lemma}
\begin{proof}
    By self-adjointness of $\covop{B}$ on $\productspace{B}$, we get that  
    $$
        \lambdamax(\covop{B}) = \sup_{f, \|f\|_{\productspace{B}} \neq 0} \frac{ \langle f, \covop{B} f \rangle_{\productspace{B}} }{ \|f\|_{\productspace{B}}^2 } = \sup_{f, \|f\|_{\productspace{B}} \neq 0} \frac{ \sum_{(i, j) \in P_B^2} \langle f_i, (\covop{B})_{i,j} f_j \rangle }{\sum_{i \in P_B} \|f_i\|^2 }.
    $$
    Now, by Cauchy-Schwarz, \Cref{obs:treecov}, and \Cref{lem:edgecovariance}, we have
    \begin{align*}
        \langle f_i, (\covop{B})_{i,j} f_j \rangle &\le \|(\covop{B})_{i,j}\| \|f_i\| \|f_j\| \\ 
        &\le \|f_i\| \|f_j\| \prod_{\ell=1}^{m-1} \|\Kuv{v_\ell}{v_{\ell+1}}\| \\
        &\le (\scalarcov)_{i,j} \|f_i\| \|f_j\| ,    
    \end{align*}
    where $(v_1=i,\ldots,v_m=j)$ is the unique path from $i$ to $j$ in $B$. Thus, setting $x \in \mathbb{R}^{|P_B|}, x_i = \|f_i\|$, we get 
    \begin{align*}
        \langle f, \covop{B} f \rangle_{\productspace{B}} \le x^\top {\scalarcov} x \le \lambdamax(\scalarcov) \|x\|^2 = \lambdamax(\scalarcov) \|f\|^2_{\productspace{B}}    
    \end{align*}
    and the lemma follows.
\end{proof}

Next, we factor the matrix $\scalarcov$ in the following way. 
\begin{lemma}[Factoring $\scalarcov$]\label{lem:Kfactor}
    Let $B$ be a tree block and $\taubf$ be a coloring on $V \setminus B$. Root the tree at an arbitrary vertex $r$ and
    denote by $\parent{v}$ the parent of a vertex $v \neq r$ in this rooted tree. Moreover, denote by $P(u, v)$ the unique path from a $u$ to $v$ in the tree $B$. 
    Define the diagonal matrix
    $$
        D \in \mathbb{R}^{|B| \times |B|} \text{ with } D_{r,r} = 1, D_{i, i} = 1 - \beta_{i, \parent{i}}^2, \text{ and } D_{i, j} = 0 \text{ whenever } i\neq j
    $$ and 
    $$
        R \in \mathbb{R}^{|B| \times |P_B|} \text{ with } R_{v, \rho} = \begin{cases}
            \prod_{e \in P(\rho, v)} \beta_{e} & \text{if } v \text{ is an ancestor of } \rho \\
            0 & \text{otherwise.}
        \end{cases}
    $$
    Then, $$
        \scalarcov = R^{\top} D R \text{ and } \|\scalarcov\| \le \|R\|^2.
    $$
\end{lemma}
\begin{proof}
    For a vertex $v$ in our rooted tree $B$, denote by $\anc{v}$ the set of ancestors of $v$.
    For two ports $i, j$ with lowest common ancestor $w$, note that $(R^{\top} D R)_{\rho_1,\rho_2}$ is the inner product of the $\rho_1$-th column of $R$ and the $\rho_2$-th column of $R$ in which the $v$-th entry is multiplied by $D_{v,v}$. 
    Hence, $$
        (R^{\top} D R)_{\rho_1,\rho_2} = \prod_{e \in P(\rho_1, \rho_2) } \beta_{e} \sum_{v \in \anc{w} } D_{v,v} \prod_{e \in P(w, v)} \beta_e^2.
    $$ We show inductively that the above sum is equal to one, which finishes the proof of the first statement. To this end, call it $S_w$ and note that $S_r = 1$. For any $w \in B$, we have 
    $$
        S_w = D_{w,w} + \beta_{w, \parent{w}}^2 S_{\parent{w}}, 
    $$ and since $S_{\parent{w}} = 1$ by induction hypothesis, it follows that $S_w = 1$, as desired. To show that $\|\scalarcov\| \le \|R\|^2$, it suffices to note that $0 \preceq D \preceq 1$.
\end{proof}

Hence, it only remains to estimate the norm of $R$.

\newcommand{\treeevbound}{\lambda_{\mathrm{tree}}}
\newcommand{\treeevboundval}{\frac{1}{1 - \zeta} + 
        \frac{\sqrt{d}\zeta^{A-2}}{1 - \omega}}

\begin{lemma}[Reducing the norm of $R$ to path sums]
\label{lem:Rpathsumreduction}
    Let $R$ be as in \Cref{lem:Kfactor}. For $v\in B$ and
    an integer $k\geq1$, define
    \[
        S_k(v)\coloneqq
          \sum_{\substack{\rho\in P_B:\ v\in\anc{\rho}\\
                           \dist_B(v,\rho)=k}}
                    \prod_{e\in P(v,\rho)}\beta_e^2.
    \]
    Then
    \[
        \|R\|\leq1+\sum_{k\geq1}
                         \sqrt{\max_{v\in B}S_k(v)}.
    \]
\end{lemma}
\begin{proof}
    Split $R$ according to the length of the path from a port
    to its ancestor:
    \[
        R=J+\sum_{k\geq1}R_k,\qquad
        J_{v,\rho}=\mathbf1_{\{v=\rho\}},\qquad
        (R_k)_{v,\rho}
            =\mathbf1_{\{\dist_B(v,\rho)=k\}}R_{v,\rho}.
    \]
    The columns of $J$ are distinct coordinate vectors, so
    $\|J\|\leq1$. For a fixed $k$, a port $\rho$ has at most
    one ancestor at distance $k$. Thus every column of $R_k$
    has at most one nonzero entry, and distinct rows have
    disjoint supports. Consequently $R_kR_k^\top$ is
    diagonal, with
    \[
        (R_kR_k^\top)_{v,v}
           =\sum_{\rho\in P_B}(R_k)_{v,\rho}^2=S_k(v).
    \]
    It follows that
    $\|R_k\|^2=\|R_kR_k^\top\|=\max_{v\in B}S_k(v)$.
    The triangle inequality now gives the claimed bound.
    All sums over $k$ are finite, since the corresponding
    matrices vanish beyond the height of the rooted tree.
\end{proof}

The operator-norm estimate has therefore been reduced to a
scalar path sum. We now use the breakpoint weights to estimate
this sum, distinguishing paths that stay low-degree from those
that contain a hub. We state the estimate for general families of
self-avoiding paths, so that it also applies to the cycle-to-port
paths in \Cref{lem:unicycliceigenvaluehubcycle}.

\begin{lemma}[Estimating the path sums]
\label{lem:Rnormbound}
    In the standing coloring regime, take $A$ and then $d$
    sufficiently large, and define
    \[
        \zeta\coloneqq\frac{\sqrt{\lowdegree}}a,\qquad
        \omega\coloneqq(d\zeta^{2A})^{1/(2(A+1))},
        \qquad a=q-\lowdegree.
    \]
    Let $B$ be a tree or unicyclic block, $v\in B$, and
    $k\geq1$ an integer. Let $\mathcal P_k(v)$ be any family
    of distinct self-avoiding paths of length $k$ contained
    in $B$, starting at $v$ and ending at a port in $P_B$.
    Define
    \[
        S_k(v)\coloneqq
           \sum_{P\in\mathcal P_k(v)}\prod_{e\in P}\beta_e^2
           =S_k^{\mathrm{low}}(v)+S_k^{\mathrm{hub}}(v),
    \]
    where the first sum contains paths whose vertices all
    have degree at most $\lowdegree$, and the second contains
    paths with at least one hub. Uniformly over the choice
    of the path family,
    \[
        S_k^{\mathrm{low}}(v)\leq\zeta^{2k},
    \]
    and
    \begin{equation}\label{eq:Rhubpathsum}
        S_k^{\mathrm{hub}}(v)\leq
        \begin{cases}
            0,&k<A-1,\\[1mm]
            \displaystyle
            \left(\frac d{\zeta^2}\right)^{2/(A+1)}
                      \omega^{2k},&k\geq A-1.
        \end{cases}
    \end{equation}
    In particular, these estimates apply to the ancestor-to-port
    paths in \Cref{lem:Rpathsumreduction}. Consequently, for
    a tree block and $R$ as in \Cref{lem:Kfactor},
    \[
        \|R\|\leq\treeevboundval.
    \]
\end{lemma}
\begin{proof}
    \leavevmode
    \paragraph{Low-degree paths}
    There are at most $\lowdegree^k$ paths of length $k$
    starting at $v$ and containing only low-degree vertices.
    Each edge on such a path has $\beta_e\leq1/a$. Hence
    \[
        S_k^{\mathrm{low}}(v)
           \leq\lowdegree^k a^{-2k}=\zeta^{2k}.
    \]

    \paragraph{Paths containing hubs}
    For a path $P=(v=v_0,v_1,\ldots,v_k=\rho)$ in
    $\mathcal P_k(v)$, define
    \[
        Q_v(P)\coloneqq\prod_{j=0}^{k-1}\deg_G(v_j),
        \qquad
        \pi_v(P)\coloneqq\frac1{Q_v(P)}.
    \]
    A simple random walk on $G$, started at $v$, follows
    this particular path during its first $k$ steps with
    probability $\pi_v(P)$. Distinct paths specify disjoint
    events, even if they end at the same port. Thus
    \[
        \sum_{P\in\mathcal P_k(v)}\pi_v(P)\leq1.
    \]
    Since $\pi_v(P)Q_v(P)=1$, it is enough to bound
    $Q_v(P)\prod_{e\in P}\beta_e^2$ for each path containing
    a hub: the path sum is a weighted sum with total weight
    at most one.

    Fix such a path, and let $\ell$ be its number of
    low-degree vertices and $h\geq1$ its number of hubs.
    Thus $\ell+h=k+1$. The endpoint $\rho$ is low-degree,
    by the breakpoint test on an edge from the exterior
    boundary to $\rho$. Set
    \[
        \gamma\coloneqq
            \prod_{\substack{u\in P\\
                             \deg_G(u)>\lowdegree}}\deg_G(u).
    \]
    The degree product $Q_v(P)$ omits $\rho$, so its
    low-degree factors are bounded by $\lowdegree^{\ell-1}$.
    Moreover, $\beta_e^2$ contributes a factor $a^{-1}$
    for each low-degree endpoint of $e$. There are at least
    $2(\ell-1)$ such incidences along the path. Therefore,
    \[
        Q_v(P)\prod_{e\in P}\beta_e^2
           \leq\lowdegree^{\ell-1}a^{-2(\ell-1)}\gamma
           =\zeta^{2(\ell-1)}\gamma.
    \]

    Choose a neighbor $b$ of $\rho$ outside $B$. It is a
    global breakpoint, and prepending $b$ to the reversed path
    gives a self-avoiding path starting at $b$. The breakpoint
    condition and \eqref{eq:weights} imply
    \[
        \alpha^{-(\ell+1)}\alpha^{Ah}\frac{\gamma}{d^h}
                  \leq1,
        \qquad\text{hence}\qquad
        \gamma\leq d^h\alpha^{\ell+1-Ah}.
    \]
    Since every hub has degree greater than $d$ and
    $\alpha>1$, this also gives $Ah\leq\ell+1$. Combining
    it with $\ell+h=k+1$ yields
    \[
        h\leq\frac{k+2}{A+1},\qquad k\geq A-1.
    \]
    In particular, there are no hub paths when $k<A-1$.

    For sufficiently large $d$,
    $d/(\alpha^{A+1}\zeta^2)\geq1$: indeed,
    $d/\lowdegree$ is bounded below, $a^2$ grows polynomially
    in $d$, and $\alpha^{A+1}$ grows only polylogarithmically.
    Using $\ell-1=k-h$, we obtain
    \begin{align*}
        Q_v(P)\prod_{e\in P}\beta_e^2
          &\leq\alpha^2(\alpha\zeta^2)^{k-h}
                              \left(\frac d{\alpha^A}\right)^h\\
          &=\alpha^2(\alpha\zeta^2)^k
                     \left(\frac d{\alpha^{A+1}\zeta^2}\right)^h\\
          &\leq\alpha^2(\alpha\zeta^2)^k
                     \left(\frac d{\alpha^{A+1}\zeta^2}\right)^{(k+2)/(A+1)}\\
          &=\left(\frac d{\zeta^2}\right)^{2/(A+1)}
                          \omega^{2k}.
    \end{align*}
    Averaging this bound with the weights $\pi_v(P)$
    proves \eqref{eq:Rhubpathsum}.

    \paragraph{Recovering the norm bound for tree blocks}
    If $B$ is a tree, apply the path estimates to the family
    of paths from $v$ to its descendant ports at distance $k$.
    By \Cref{lem:Rpathsumreduction} and
    $\sqrt{x+y}\leq\sqrt x+\sqrt y$,
    \[
        \|R\|\leq
           1+\sum_{k\geq1}\zeta^k+
             \left(\frac d{\zeta^2}\right)^{1/(A+1)}
                           \sum_{k\geq A-1}\omega^k.
    \]
    Our choices of $A,d$ ensure $\zeta,\omega<1$. Summing
    the geometric series and using
    \[
        \left(\frac d{\zeta^2}\right)^{1/(A+1)}
                        \omega^{A-1}
             =\sqrt d\,\zeta^{A-2}
    \]
    gives the claimed norm bound.
\end{proof}

\begin{proof}[Proof of \Cref{lem:treeeigenvalue}]
    By \Cref{lem:scalarop}, we have $\lambdamax(\covoppure{B}) \le \lambdamax(\scalarcov)$. By \Cref{lem:Kfactor} and \Cref{lem:Rnormbound}, 
    $$
        \lambdamax(\scalarcov) \le \|R\|^2 \le \left(\treeevboundval\right)^2.
    $$ Because $q \ge \numcolors$ and $\lowdegree = \lowdegreedef$ while $\eta < \varepsilon$, we can choose $d, A$ large enough such that this is at most $1 + \delta$, as desired.
\end{proof}

\subsubsection{Unicyclic block eigenvalue bound}

The essentially same bound holds for unicyclic blocks by re-using the tree proof and conditioning on spins associated to the cycle. 

\begin{lemma}[Unicyclic block covariance operator eigenvalue bound]\label{lem:cycleeigenvalue}
    Consider any unicyclic block $B$. Then, for every coloring $\taubf$ on $V \setminus B$ and every $\delta > 0$, we have
    $$
        \lambdamax(\covop{B}) \le (1 + \delta).
    $$ provided that $d$ and $A$ are large enough.
\end{lemma}
To reuse the tree proof, we consider two cases.

\paragraph{Cycles with a low-degree vertex}
We delete a cycle edge incident to a low-degree vertex and apply the
tree covariance bound from \Cref{lem:treeeigenvalue}, denoted by
$\treeevbound$, to the resulting tree. Restoring the edge amounts to
conditioning its endpoints to have different colors. This event has
high probability even after fixing the color at any one port, which
lets us compare variances under the two measures and transfer the
tree bound back to $B$ with a factor $1+1/a$.

\begin{lemma}\label{lem:unicycliceigenvaluelowcycle}
    Assume $A\geq2$.
    Consider any unicyclic block $B$ in which the unique cycle $C$ contains a vertex of degree $\le \lowdegree$. Then there exists an edge $e$ in $B$ such that 
    $$
        \lambdamax(\covop{B}) \le  \left( 1 + \tfrac{1}{a} \right) \treeevbound
    $$ for every boundary coloring $\taubf$ on $V \setminus B$. Here, $\treeevbound$ is the eigenvalue bound for tree blocks that follows from \Cref{lem:treeeigenvalue}. 
\end{lemma}
\begin{proof}
    Due to self adjointness of $\covop{B}$, we have 
    $$
        \lambdamax(\covop{B}) \le \sup_f \frac{\langle f, \covop{B} f \rangle_{\productspace{B}}}{\|f\|^2_{\productspace{B}}} = \sup_f \frac{\Varsub{\mu}{ F }}{ \sum_{\rho \in P_B} \Varsub{\mu}{f_\rho } }
    $$ where $\mu$ is the distribution of colors in $B$  conditional on $\taubf$, the function $F$ is defined as $F(\sigmabf) = \sum_{\rho} f_\rho(\sigmabf_\rho)$, and where the supremum is over $f \in \productspace{B}$ with $\|f\|_{\productspace{B}} \neq 0$.

    Now, choose a cycle edge $e = \{u, v\}\in E(C)$ incident to a vertex of degree $\le \lowdegree$. Note that removing $e$ from $B$ renders the resulting graph a tree. Denote by $\nu$ the distribution of colorings in $B$ after removing the edge $e$, conditional on $\taubf$. Set $\mathcal E=\{\sigmabf_u\neq\sigmabf_v\}$, so $\mu=\nu(\cdot\mid\mathcal E)$. By the law of total variance, we have
    $$
        \Varsub{\nu}{F} \ge \Prsub{\nu}{\mathcal{E}} \ \Varsub{\mu}{ F }.
    $$ Moreover, 
    $$
        \Prsub{\mu}{\sigmabf_\rho = c} = \Prsub{\nu}{\sigmabf_\rho = c \mid \mathcal{E}} = \Prsub{\nu}{\sigmabf_\rho = c} \frac{\Prsub{\nu}{\mathcal{E} \mid \sigmabf_\rho = c} }{\Prsub{\nu}{\mathcal{E}}}
    $$
    and we claim that
    \begin{align} \label{eq:deletededgeprob}
        \Prsub{\nu}{\mathcal{E} \mid \sigmabf_\rho = c} \ge 1 - \frac{1}{a+1} = \frac{a}{a+1}.    
    \end{align}
    For $A\geq2$, a port $\rho$ cannot have a hub neighbor
    $v\in B$, since for any $s\in\partial B$ adjacent to $\rho$ the
    path $(s,\rho,v)$ would give
    $W(s)W(\rho)W(v)=\alpha^{A-2}\deg_G(v)/d>1$,
    contradicting the breakpoint condition at $s$.
    Indeed, for $e = \{u,v\}$, we can assume that one of the vertices in $e \setminus \{\rho\}$ is of degree $\le \lowdegree$: if $e$ is not incident to the port $\rho$ this is the case by construction, otherwise both endpoints of $e$ are low-degree as otherwise the breakpoint condition of the boundary vertices would be violated. Denote this low-degree endpoint in $e \setminus \{\rho\}$ by $u$. Now, after deleting the edge $e = \{u, v\}$ and revealing the colors of all vertices in $B \setminus \{u\}$, the vertex $u$ has at least $a + 1$ available colors and thus matches the color of $v$ with probability at most $\frac{1}{a+1}$, which yields \eqref{eq:deletededgeprob}. In total, we thus get 
    $$
        \Varsub{\nu}{f_\rho} \le \Expectedsub{\nu}{f_\rho^2} \le \frac{a+1}{a} \Prsub{\nu}{\mathcal{E}} \ \Expectedsub{\mu}{f_\rho^2} = \left( 1 + \tfrac{1}{a}\right) \Prsub{\nu}{\mathcal{E}} \ \Varsub{\mu}{f_\rho}. 
    $$
    Keeping the original degrees, weights and ports, the tree-block proof applies under $\nu$, since deleting $e$ only removes paths from the breakpoint tests.
    Hence, combining everything and applying \Cref{lem:treeeigenvalue} to the tree resulting after removing the edge $e$ from $B$, we get
    \begin{align*}
        \Varsub{\mu}{F} &\le \frac{1}{\Prsub{\nu}{\mathcal{E}}} \Varsub{\nu}{F} \\
        &\le \frac{\treeevbound}{\Prsub{\nu}{\mathcal{E}}} \sum_{\rho \in P_B} \Varsub{\nu}{f_\rho} \le \left( 1 + \tfrac{1}{a}\right)\treeevbound \sum_{\rho \in P_B} \Varsub{\mu}{f_\rho},
    \end{align*} as desired.
\end{proof}

\paragraph{Cycles consisting only of hubs} 

If the cycle in $B$ consists only of high-degree vertices, we condition on the colors in the whole cycle $C$ and use a similar correlation decay implied by the breakpoint property as used for tree-blocks.

\begin{lemma}\label{lem:unicycliceigenvaluehubcycle}
    Consider any unicyclic block $B$ in which the unique cycle $C$ contains only vertices of degree $> \lowdegree$. Then,
    $$
        \lambdamax(\covop{B}) \le \treeevbound + \frac{d\zeta^{2(A-2)}}{1 - \omega^2}
    $$ for all boundary colorings $\taubf$. Here, $\zeta, \omega$ are as defined in \Cref{lem:Rnormbound}, and  $\treeevbound$ is the eigenvalue bound for tree blocks that follows from \Cref{lem:treeeigenvalue}.
\end{lemma}
\begin{proof}
    Given a test function $f \in \productspace{B}$, we define $F(\sigmabf) = \sum_{\rho \in P_B} f_\rho(\sigmabf_\rho)$. By the law of total variance, we get 
    \begin{align}\label{eq:variancedecomp}
        \Var{F} = \Expectedbig{ \Var{F \mid \sigmabf_C} } + \Varbig{ \Expected{F \mid \sigmabf_C} }.
    \end{align}
    After conditioning on $\sigmabf_C$, the remaining block is a
    forest, so the bound from \Cref{lem:treeeigenvalue} remains
    applicable. The tree proof carries over componentwise after conditional centering, retaining the original degrees, weights and ports, since each port keeps its original exterior breakpoint neighbor and the additional pinned colors only modify the lists. Averaging over the cycle coloring gives
    \[
        \Expected{\Var{F\mid\sigmabf_C}}
          \leq\treeevbound\sum_{\rho\in P_B}
                       \Expected{\Var{f_\rho\mid\sigmabf_C}}
          \leq\treeevbound\sum_{\rho\in P_B}\Var{f_\rho},
    \]
    where the last inequality follows from the law of total
    variance. To handle the second term in
    \Cref{eq:variancedecomp}, we write
    $$
        \Expected{F \mid \sigmabf_C} = \sum_{\rho \in P_B} g_\rho(\sigmabf) \text{ with } g_\rho(\sigmabf) = \Expected{f_\rho \mid \sigmabf_C}.
    $$
    Now, 
    $$
        \Varbig{ \Expected{F \mid \sigmabf_C} } = \Big\| {\sum}_{\rho \in P_B} g_\rho \Big\|^2 \le \left( {\sum}_{\rho \in P_B} \|g_\rho\| \right)^2.
    $$ Now, for a port $\rho$, denote by $c(\rho)$ the unique vertex in $C$ at which the tree which $\rho$ is part of is attached. Every edge of the path from $c(\rho)$ to $\rho$ is a bridge, so the edge-deletion proof of \Cref{lem:edgecovariance} and the conditional-independence argument of \Cref{obs:treecov} carry over unchanged. Recall the operator $\Kuv{u}{v}$ from \Cref{obs:treecov} and note that 
    $$
        g_\rho(\sigmabf) = \Kuv{c(\rho)}{\rho} f_\rho (\sigmabf_{c(\rho)}) \text{ and thus } \|g_\rho\| \le \|\Kuv{c(\rho)}{\rho}\|\|f_\rho\| \le \|f_\rho\| \prod_{e \in P(c(\rho), \rho)} \beta_e ,
    $$ where $P(c(\rho), \rho)$ is the unique path from $c(\rho)$ to $\rho$, and where we applied \Cref{lem:edgecovariance} together with \Cref{obs:treecov}. Consequently, by Cauchy-Schwarz, we get 
    \begin{align*}
        \left( {\sum}_{\rho \in P_B} \|g_\rho\| \right)^2 &\le \left( \sum_{\rho \in P_B} \prod_{e \in P(c(\rho), \rho)} \beta_e^2 \right)\left( \sum_{\rho \in P_B} \|f_\rho\|^2 \right) \\
        &= \left( \sum_{v \in C} \sum_{\substack{\rho \in P_B \\ c(\rho) = v} } \prod_{e \in P(c(\rho), \rho)} \beta_e^2 \right)\sum_{\rho \in P_B} \Var{f_\rho}.
    \end{align*} To estimate the factor before the sum of variances, we apply
    \Cref{lem:Rnormbound} directly to paths from the cycle to
    the ports. For $v\in C$ and $k\geq1$, let
    \[
        \mathcal P_k(v)\coloneqq
          \{P(v,\rho):\rho\in P_B,\ c(\rho)=v,\
                                  \dist_B(v,\rho)=k\},
    \]
    where $P(v,\rho)$ is the unique path in the tree attached
    to $C$ at $v$. These are distinct self-avoiding paths in
    $B$ ending at ports, and each contains the hub $v$.
    Therefore the hub-path estimate \eqref{eq:Rhubpathsum}
    gives
    \newcommand{\kcorrelationsum}[2]{\Gamma_{#1}(#2)}
    \[
        \kcorrelationsum{v}{k}
          \coloneqq\sum_{P\in\mathcal P_k(v)}
                              \prod_{e\in P}\beta_e^2
          \leq\left(\frac d{\zeta^2}\right)^{2/(A+1)}
                              \omega^{2k}.
    \] Now, we claim that every port $\rho \in P_B$ has distance at least $A|C| - 1$ from the cycle $C$. To see this, start at an exterior breakpoint adjacent to $\rho$, follow the path through $\rho$ to $c(\rho)$, and then traverse $C$ without taking the closing edge. This path contains at least $|C|$ hubs. By the breakpoint property and because $\deg(v)/d > 1$ for every hub, this implies $\dist_B(\rho, c(\rho)) + 1 \ge A|C|$, as desired. Combining everything, we thus get
    \begin{align*}
        \sum_{v \in C} \sum_{\substack{\rho \in P_B \\ c(\rho) = v} } \prod_{e \in P(c(\rho), \rho)} \beta_e^2 &= \sum_{v \in C}\sum_{k \ge A|C| - 1} \kcorrelationsum{v}{k} \\
        &\le |C| \left( \frac{d}{\zeta^2} \right)^{\frac{2}{A+1}} \frac{\omega^{2(A|C|-1)}}{1 - \omega^2} \le \frac{d \zeta^{2(A-2)}}{1 - \omega^2} |C| \omega^{2A(|C|-1)} \le \frac{d \zeta^{2(A-2)}}{1 - \omega^2}.
    \end{align*} for sufficiently large $d, A$. Combining everything, the lemma follows.
\end{proof}

It only remains to prove \Cref{lem:cycleeigenvalue}.
\begin{proof}[Proof of \Cref{lem:cycleeigenvalue}]
    Since both \Cref{lem:unicycliceigenvaluelowcycle} and \Cref{lem:unicycliceigenvaluehubcycle} yield $\lambdamax(\covop{B}) \le 1 + \delta$ given $d$ and $A$ large enough, the lemma follows.
\end{proof}

\subsection{Assembling the Bochner inequality via the smallest eigenvalue}\label{sec:assembling}

Starting at \eqref{eq:bochner} and using \Cref{cor:lowdegreebochner} and \Cref{cor:hubbochner}, we get the following.
\begin{corollary}\label{cor:assembledbochner}
    Let $\Atil$ be the adjacency matrix of the graph $\Gtil$ on $V$ with
    \[
        E(\Gtil)=\{\{u,v\}\in E(G):u\in S\text{ or }v\in S\}.
    \]
    Recall that $S$ denotes the set of singletons in our block decomposition. For any color $c \in [q]$, let $\gradveccomplete{c} \in \mathbb{R}^{|V|}$ with
    $$
        \gradveccomplete{c}(v) = \begin{cases}
            \gradvec{c}(v) & \text{if } v \in S\\
            \portgrad{c}(v) & \text{if } v \text{ is a port in some block } B\\
            0 & \text{otherwise.}
        \end{cases}
    $$
    Under the assumptions of \Cref{cor:hubbochner}, we have
    \begin{align*}
        &\langle f, B_\kappa f \rangle \ge  \min\left\{ \frac{a- \kappa q}{2},  \frac{a(1-\kappa)}{2\lambda}-\frac{4\lowdegree}{a} \right\} \sum_{c \in [q]} \ExpectedBig{ \| \gradveccomplete{c} \|_2^2 } + \sum_{c \in [q]} \Expected{ \big(\gradveccomplete{c}\big)^{\!\top} \!\! \Atil \big(\gradveccomplete{c}\big) },
    \end{align*} where $\lambda$ is the covariance upper bound from \Cref{cor:hubbochner}.
\end{corollary}

To compare the potentially negative term arising from $(\gradveccomplete{c}\big)^{\!\top} \!\! \Atil \big(\gradveccomplete{c}\big)$ to the strictly positive first term, we use the fact that the smallest eigenvalue of $\Atil$ is sufficiently close to $0$ as it is of order $\sqrt{d}$ in absolute value. This follows because the smallest eigenvalue of $G(n, d/n)$ is $O(\sqrt{d})$ once we remove all high-degree vertices. In particular, this is the case for the graph induced on the singletons in $S$ and all ports. Moreover, once we delete all the edges that are entirely within some block, we delete a union of trees with at most one extra edge. The largest eigenvalue of this deleted graph therefore is also bounded by $O(\sqrt{d})$ such that the smallest eigenvalue of $\Atil$ stays of order $O(\sqrt{d})$ in absolute value. 
\begin{lemma}[Smallest eigenvalue of $\Atil$]\label{lem:smallesteigenvalue}
    Let $G\sim G(n,d/n)$, assume $A\geq1$, and take $d$ sufficiently
    large. There is an absolute constant $C$ such that, with probability
    $1-o(1)$ as $n\to\infty$,
    $$
        \lambdamin(\Atil) \geq -C\sqrt d.
    $$
\end{lemma}
\begin{proof}
    Write $p=d/n$ and let $A_G^{\mathrm{tr}}$ be the adjacency matrix
    obtained by zeroing all rows and columns of $A_G$ corresponding to
    vertices of degree greater than $2d$. The regularization theorem
    of Le, Levina, and Vershynin
    \cite[Theorem~1.1]{le2017concentration} (see also
    \cite{feige2005spectral}) gives, with high probability,
    \[
        \|A_G^{\mathrm{tr}}-p(J-I)\|\leq C_0\sqrt d
    \]
    for an absolute constant $C_0$, where $J$ is the all-ones matrix.
    Let $U=S\cup\bigcup_B P_B$. All vertices in $U$ have degree at most
    $\lowdegree$: this follows from the single-vertex breakpoint test
    for $S$, and, for a port $x$ adjacent to $s\in S$, from the test
    on $(s,x)$, since a high-degree port would give
    $W(s)W(x)=\alpha^{A-1}\deg_G(x)/d>1$.
    For large $d$, $\lowdegree\leq2d$, so $A_{G[U]}$ is a principal
    submatrix of $A_G^{\mathrm{tr}}$. Since $pJ$ is positive
    semidefinite, the preceding estimate implies
    \[
        \lambdamin(A_{G[U]})\geq -C_0\sqrt d-p.
    \]
    This applies even though $U$ depends on $G$.

    Let $F$ consist of the edges of $G[U]$ lying entirely inside a block.
    Its components are trees or unicyclic, and its maximum degree is at
    most $\lowdegree$. Orient every edge of a tree component from a
    child to its parent, after choosing a root. For a unicyclic
    component, orient the unique cycle consistently around the cycle
    and all remaining edges toward the cycle. Each non-root vertex
    in a tree, and every vertex in a unicyclic component, then has
    exactly one outgoing edge; roots have none.

    Define the oriented adjacency matrix $M\in\mathbb R^{U\times U}$ by
    \[
        M_{uv}=
        \begin{cases}
            1,&\text{if the edge $\{u,v\}$ is oriented from $u$ to $v$},\\
            0,&\text{otherwise}.
        \end{cases}
    \]
    Each undirected edge is given exactly one orientation, so
    $(M+M^\top)_{uv}=1$ precisely when $\{u,v\}\in E(F)$, and is
    zero otherwise. Thus $A_F=M+M^\top$.
    The row sum at $u$ is its outdegree, while the column sum at $v$
    is its indegree. Consequently the induced matrix norms satisfy
    \[
        \|M\|_\infty=\max_u\sum_v|M_{uv}|\leq1,
        \qquad
        \|M\|_1=\max_v\sum_u|M_{uv}|\leq\lowdegree.
    \]
    To see why these bounds control the Euclidean operator norm,
    apply Cauchy--Schwarz to each row: for any $z\in\mathbb R^U$,
    \begin{align*}
        \|Mz\|_2^2
            &=\sum_u\left|\sum_v M_{uv}z_v\right|^2\\
            &\leq\sum_u\left(\sum_v M_{uv}\right)
                            \left(\sum_v M_{uv}|z_v|^2\right)\\
            &\leq\|M\|_\infty
                     \sum_v\left(\sum_u M_{uv}\right)|z_v|^2\\
            &\leq\|M\|_\infty\|M\|_1\|z\|_2^2.
    \end{align*}
    Hence $\|M\|\leq\sqrt{\|M\|_1\|M\|_\infty}$. Finally, transposition
    preserves this norm, and the triangle
    inequality gives
    \[
        \|A_F\|=\|M+M^\top\|
            \leq\|M\|+\|M^\top\|
            =2\|M\|
            \leq2\sqrt{\|M\|_1\|M\|_\infty}
            \leq2\sqrt{\lowdegree}.
    \]
    On $U$, $\Atil=A_{G[U]}-A_F$, and outside $U$ its rows and columns
    vanish. Therefore, 
    \[
        \lambdamin(\Atil)
            \geq-C_0\sqrt d-p-2\sqrt{\lowdegree}
            \geq-C\sqrt d,
    \]
    as desired.
\end{proof}

Finally, this yields a spectral gap as it implies that $B_\kappa$ is PSD.
\begin{proof}[Proof of \Cref{thm:spectralgapblockcoloring}]
    Put $q_0=d+d^{1/2+\varepsilon}$ and
    $r_d=q_0-\lowdegree=d^{1/2+\varepsilon}-d^{1/2+\eta}$.
    By \Cref{lem:treeeigenvalue,lem:cycleeigenvalue}, choose $A$ and
    then $d$ large enough that the uniform covariance bound is at
    most $2$. We may therefore use $\lambda=2$ in
    \Cref{cor:assembledbochner}. On the high-probability event from
    \Cref{lem:smallesteigenvalue}, we also have
    $\lambdamin(\Atil)\geq-C\sqrt d$.

    Take $\kappa=\kappa_d=r_d/(4q_0)$. For every $q\geq q_0$,
    $a=q-\lowdegree\geq r_d$, $\kappa\leq1/4$, and
    $\kappa q\leq a/4$, since
    $a/q=1-\lowdegree/q\geq r_d/q_0$.
    Since $r_d/\sqrt d=d^\varepsilon-d^\eta\to\infty$ and
    $\lowdegree/d\to1$, increase $d$ so that
    $r_d^2\geq64\lowdegree$ and $r_d\geq8C\sqrt d$.
    Then, with $\lambda=2$,
    \[
        \frac{8\lambda\lowdegree}{a^2}
           =\frac{16\lowdegree}{a^2}\leq\frac14\leq1-\kappa,
    \]
    so the assumption in \Cref{cor:hubbochner,cor:assembledbochner}
    holds. Moreover,
    \[
        \min\left\{\frac{a-\kappa q}{2},
                     \frac{a(1-\kappa)}{2\lambda}
                          -\frac{4\lowdegree}{a}\right\}
          \geq\min\left\{\frac{3a}{8},
                          \frac{3a}{16}-\frac{4\lowdegree}{a}\right\}
          \geq\frac a8.
    \]
    Thus \Cref{cor:assembledbochner} and the Rayleigh quotient give,
    for every $f$,
    \[
        \langle f,B_\kappa f\rangle
            \geq\left(\frac a8-C\sqrt d\right)
                     \sum_{c\in[q]}
                         \ExpectedBig{\|\gradveccomplete{c}\|_2^2}
            \geq0.
    \]
    Since $q\geq\lfloor\lowdegree\rfloor+2$ for sufficiently large $d$,
    \Cref{lem:comparisonirreducibility} shows that the kernel of $L$
    consists only of constant functions. The form comparison
    $L\preceq\Ltil\preceq qL$ gives the same kernel for $\Ltil$.
    Hence $B_\kappa=\Ltil^2-\kappa q\Ltil\succeq0$, and
    \Cref{lem:acceleratedgap} yields $\gap{L}\geq\kappa_d$.
    Since $\eta<\varepsilon$, we have
    $\kappa_d=\Theta_{\varepsilon,\eta}
        (\min\{1,d^{\varepsilon-1/2}\})$ as $d\to\infty$,
    giving the stated asymptotic bound. \qedhere

\end{proof}

%% file: breakpoints.tex
In this section, we show that $G(n, d/n)$ admits a decomposition into tree/unicyclic blocks with the properties stated in \Cref{def:validblockpartion} and probability $1- o(1)$ once $d$ is large enough. To this end, we first show that the connected components left after removing all global breakpoints $\widehat{S}$ have the desired tree/unicyclic and diameter properties in the following subsection. Afterwards in \Cref{sec:merging}, we show that some of these blocks can be merged such that the resulting partition also satisfy condition 2 from \Cref{def:validblockpartion} ensuring that every singleton has at most one neighbor in every given block. 

To introduce this more formally, recall the definition of breakpoints and weights from \Cref{sec:blockassumptions}: we have some 
\(0<\eta<1/2\), \(A>0\), and \(C>0\), and define
\[
    \lowdegree \coloneqq d+d^{1/2+\eta} , \text{ and }
    \alpha \coloneqq C\log d>1.
\]
Moreover, the vertex weight $W(v)$ associated to every vertex $v$ is given by
\begin{align}\label{eq:modified-breakpoint-weights}
    W(v) \coloneqq
    \begin{cases}
        \alpha^{-1}, & \deg_G(v)\leq \lowdegree,\\[1mm]
        \alpha^A\deg_G(v)/d, & \deg_G(v)>\lowdegree.
    \end{cases}
\end{align}
Throughout this subsection, paths are self-avoiding, their length is their
number of edges, and a single vertex is also a path. For such a path $P$, set
\[
    W(P)=\prod_{v\in V(P)}W(v).
\]
Call $P$ \emph{heavy} if $W(P)>1$. A vertex is a \emph{global
breakpoint} if no heavy path starts there. Write $\widehat S$ for the set
of global breakpoints, and $\widehat S_r$ for the set of vertices from which
no heavy path of at most $r$ edges starts. Thus
$\widehat S\subseteq\widehat S_r$. Let $\widehat{\mathcal B}$ be the
collection of vertex sets of the connected components of
$G[V\setminus\widehat S]$. Diameters of these components are measured
inside the components, rather than in the ambient graph.

The overall goal of this section is to prove the following theorem.
\begin{theorem}[Valid block partition in $G(n, d/n)$]
\label{thm:validblockpartition}
    Given any fixed $A\geq0$, $C>0$, and $0<\eta<1/2$,
    there are constants $K>0$ and $d_0>1$, depending only on
    $\eta,A,C$, with the following property. For all $d\geq d_0$,
    $G\sim G(n,d/n)$, and
    \[
        h=\left\lceil K\frac{\log n}{d^{2\eta}}\right\rceil,
    \]
    there is a set $S \subseteq \widehat{S}$ such that the resulting block partition is valid as specified by \Cref{def:validblockpartion}, while the diameter of every block satisfies
    \[
        \operatorname{diam}(G[B])< 16h
        \leq16(K+1)\frac{\log n}{d^{2\eta}}.
    \] These properties hold simultaneously with probability $1-o(1)$ as $n\to\infty$ over the draw of $G \sim G(n, d/n)$.
\end{theorem}

\subsection{Small components after removing breakpoints}
\label{subsec:modified-breakpoints}

We show the following preliminary block decomposition theorem which is central to our final block decomposition.
\begin{theorem}[Small components after deleting breakpoints]
\label{thm:modified-breakpoint-components}
    Given any fixed $A\geq0$, $C>0$, and $0<\eta<1/2$, there are constants $K>0$ and $d_0>1$, depending only on
    $\eta,A,C$, such that for all $d\geq d_0$,
    $G\sim G(n,d/n)$, and 
    \[
        h=\left\lceil K\frac{\log n}{d^{2\eta}}\right\rceil,
    \]
    the following properties hold
    simultaneously, with probability $1-o(1)$ as $n\to\infty$:
    \begin{enumerate}[(i)]
        \item $G[V\setminus\widehat S]$ contains no self-avoiding path
        on $8h$ vertices. In particular, every
        $B\in\widehat{\mathcal B}$ satisfies
        \[
            \operatorname{diam}(G[B])<8h
            \leq8(K+1)\frac{\log n}{d^{2\eta}}.
        \]
        \item Every $G[B]$, $B\in\widehat{\mathcal B}$, is a tree or
        is unicyclic, meaning that it has exactly one cycle.
        \item $\widehat S_r=\widehat S$ for every integer $r\geq h$.
    \end{enumerate}
    The constants in the diameter bound are independent of $d$ and $n$;
    the limit is taken with $d$ fixed.
\end{theorem}

The proof builds on the following three ingredients.
\begin{enumerate}
    \setlength{\itemsep}{0cm}
    \item First, we show that no heavy path of length $> h$ exists with high probability. This is achieved by standard tail estimates on vertex weights and a repeated application of Cauchy-Schwarz to handle dependencies between degrees.
    \item Only showing that no long heavy paths is insufficient for our desired diameter bounds. After all, a long path $P$ of non-breakpoints must not itself be heavy, it could as well be light and merely contain sufficiently heavy paths branching away from $P$. To show that no long non-breakpoint paths can exist nonetheless, we therefore show that every path $P$ of length $\ge h$ in $G[V \setminus \widehat{S}]$ is \emph{non-elementary} meaning that two of its vertices are connected by some path of length $\le 2h$ outside of $P$. Intuitively this is the case because a light path $P$ of non-breakpoints must contain many heavy witnesses branching away from it. However, using the bounds from step one, we know that these heavy witnesses appear very sparsely in $G$. A light path $P$ can therefore only be breakpoint-free if its vertices \emph{share} some of their heavy witnesses, implying that $P$ is non-elementary.
    \item Finally, to fully rule out paths of length $> 8h$ in $G[V \setminus \widehat{S}]$, we show that every ball of radius $10h$ in $G$ is a tree or unicyclic by a standard first moment calculation. Together with the fact that every long path in $G[V \setminus \widehat{S}]$ is non-elementary, this simple fact rules out any component of diameter $>8h$ in $G[V \setminus \widehat{S}]$.
\end{enumerate}
Since this argument also works whenever we replace the set $\widehat{S}$ of global breakpoints by the set $\widehat{S}_h$ of $h$-breakpoints, we further get that $\widehat{S}_h = \widehat{S}$.

\subsubsection{A product estimate for the vertex weights.}
Let $w(k)$ denote the weight in
\eqref{eq:modified-breakpoint-weights} assigned to a vertex of degree $k$ and define the parameters
$$
    \xi = d^{2\eta},\qquad
    \theta=\frac{1}{128(A+1)},\qquad
    \gamma=\frac{\theta}{2},\qquad
    K=\frac{64}{\gamma},\qquad
    t=\frac{\theta\xi}{\log\alpha}.
$$ In the following, all statements hold for all $\eta < 1/2$, $A, C > 0$, provided that $d$ is large enough as a function of $\eta, A, C$.
We first show the following tail bound on $w(X)$ for $X\sim\operatorname{Bin}(n,d/n)$.
\begin{lemma}[Moment bound for a vertex weight]\label{lem:breakpoint-single-moment}
    For $d$ sufficiently large, all integers $n \geq d$, and $X\sim\operatorname{Bin}(n,d/n)$, we have
    \begin{equation*}
        \EE\bigl[w(3+X)^{2t}\bigr]\leq\mathrm e^{-\theta\xi}.
\end{equation*}
\end{lemma}
\begin{proof}
    Write $D=3+X$ and $\tau = d^{\frac{1}{2} + \eta}$. If $D \le d + \tau$, then
$w(D)^{2t} \le \alpha^{-2t}=\mathrm e^{-2\theta\xi}$. To bound the contribution of the event that $D > d + \tau$, define
\[
    u=\frac{\tau}{2d},\qquad z=u+\frac{2t}{d},
\]
and note that due to $\1_{\{D>d+\tau\}}\leq\mathrm e^{u(D-d-\tau)}, $ $\log x\leq x-1$, and $\EE[\mathrm e^{zX}]\leq\exp(d(\mathrm e^z-1))$, we get 
\begin{align*}
    \EE\left[(D/d)^{2t}\1_{\{D>d+\tau\}}\right]
    &\leq\mathrm e^{-u\tau}\EE\left[\mathrm e^{z(D-d)}\right]\\
    &\leq\exp\left(-\frac{\xi}{2}+3z
                    +d(\mathrm e^z-1-z)\right).
\end{align*}
Since $\eta<1/2$, we have $t=o_d(\tau)$ and $\tau=o_d(d)$, so
\[
    z=\left(\frac12+o_d(1)\right)\frac{\tau}{d},\qquad
    3z=o_d(\xi),\qquad
    d(\mathrm e^z-1-z)=\left(\frac18+o_d(1)\right)\xi.
\]
Consequently the preceding expectation is at most
$\mathrm e^{-\xi/8}$ for sufficiently large $d$. Including the factor
$\alpha^{2At}=\mathrm e^{2A\theta\xi}$ from the high-degree weight gives
\[
    \EE\bigl[w(D)^{2t}\bigr]
    \leq\mathrm e^{-2\theta\xi}
       +\mathrm e^{-(1/8-2A\theta)\xi}
    \leq 2\mathrm e^{-2\theta\xi}
    \leq\mathrm e^{-\theta\xi}.
\]
since $1/8-2A\theta\geq2\theta$ and $\theta\xi\geq\log2$ for sufficiently large $d$.
\end{proof}

As a consequence of the above, we get the following bound on the probability that a fixed path is heavy. 
\begin{lemma}[Tail bound for products of vertex weights]\label{lem:breakpoint-product-tail}
    For $d$ sufficiently large, every fixed graph $F$ with maximum degree at most three and every $T \subseteq V(F)$,
    \begin{align}
        \Prlr{\prod_{v\in T}W(v)>1\,
              \mid F\subseteq G}
    \le \mathrm e^{-\gamma\xi|T|}
    \label{eq:breakpoint-product-tail}
\end{align}
\end{lemma}
\begin{proof}
    Applying Markov's inequality to the $t$-th power gives
    $$
        \Prlr{\prod_{v\in T}W(v)>1\,
              \mid F\subseteq G} \le \Expected{\prod_{v\in T}W(v)^{t}
              \mid F\subseteq G}.
    $$
    The power $t$ is needed to use the $2t$-moment bound from
    \Cref{lem:breakpoint-single-moment}; the unpowered application would
    lose a factor $t$ in the decay exponent. All expectations below
    are taken conditional on $F\subseteq G$.
    Now, not every factor in the above product is independent, but any two factors $W(u), W(v)$ depend on at most one edge that has not been revealed yet. Let $e_1, \ldots, e_m$ be those missing edge indicators for $m = \binom{n}{2} - |E(F)|$. For $v\in T$ and all integers $0 \le k \le m$ define
    \newcommand{\fk}[2]{f_{#1}^{(#2)}}
    $$
        \fk{v}{k} \coloneqq \sqrt{\Expected{W(v)^{2t} \mid e_{k+1}, \ldots, e_m}}.
    $$
    For $v\notin T$, set $\fk{v}{k}=1$; this also handles edges
    with fewer than two endpoints in $T$.
    Then we prove by induction over $k$ that
    $$
        \Expected{\prod_{v\in T}W(v)^{t}
              \mid F\subseteq G} \le \Expected{ \prod_{v\in T} \fk{v}{k}
              }.
    $$ Indeed, for $k=0$ equality holds because
    $\fk{v}{0}=W(v)^t$ for $v\in T$. Moreover, given any $k$, let $u, w$ be the endpoints of the edge indicator $e_{k+1}$. Then, by Cauchy-Schwarz,
    \begin{align*}
        \Expected{ \prod_{v\in T} \fk{v}{k}} &= \Expected{ \Expected{ \prod_{v\in T} \fk{v}{k}
               \mid e_{k+2 }, \ldots, e_m  } } \\
        &= \Expected{ \prod_{v\in T, v \neq u, w} \fk{v}{k} \Expected{ \fk{u}{k}\fk{w}{k}
               \mid e_{k+2 }, \ldots, e_m  } } \\
               &\le \Expected{ \prod_{v\in T, v \neq u, w} \fk{v}{k} \left( \Expected{ (\fk{u}{k})^2 \mid e_{k+2 }, \ldots, e_m } \Expected{(\fk{w}{k})^2
               \mid  e_{k+2 }, \ldots, e_m  } \right)^{1/2} } \\
        &= \Expected{ \fk{u}{k+1}\fk{w}{k+1} \prod_{v\in T, v \neq u, w} \fk{v}{k} } = \Expected{ \prod_{v\in T} \fk{v}{k+1} },
    \end{align*}
    since $\Expected{ (\fk{\ell}{k})^2 \mid e_{k+2 }, \ldots, e_m }^{1/2} = \fk{\ell}{k+1}$ and $\fk{v}{k} = \fk{v}{k+1}$ for $v \neq u, w$.
    Taking $k=m$ in the induction gives
    $$
        \Expected{\prod_{v\in T}W(v)^t
              \mid F\subseteq G} \le \Expected{ \prod_{v\in T} \fk{v}{m}
              } = \prod_{v\in T} \Expected{W(v)^{2t} \mid F\subseteq G }^{1/2}.
    $$
    Conditional on $F\subseteq G$, the degree of $v\in T$ has law
    $\deg_F(v)+\operatorname{Bin}(n-1-\deg_F(v),d/n)$ and is
    stochastically bounded by $3+X$, where
    $X\sim\operatorname{Bin}(n,d/n)$. Since $w$ is nondecreasing,
    \Cref{lem:breakpoint-single-moment} therefore yields
    \[
        \Expected{W(v)^{2t}\mid F\subseteq G}
          \leq\EE\bigl[w(3+X)^{2t}\bigr]
          \leq\mathrm e^{-\theta\xi}.
    \]
    Substituting this in the product above gives
    $\mathrm e^{-\theta\xi|T|/2}=\mathrm e^{-\gamma\xi|T|}$,
    as claimed.
\end{proof}

Now, we can show that every heavy path is small using a simple first moment estimate.
\begin{lemma}[Heavy paths are short]\label{lem:breakpoint-short-witnesses}
    Let $\mathcal E$ be the event that there is no heavy path on $h$ or more vertices. Then for $d$ sufficiently large,
    $$\Pr{\mathcal{E} } \ge 1 - 2n^{-31}.$$
\end{lemma}
\begin{proof}
    A specified path on $m$ vertices forces a graph of
maximum degree two, so \Cref{lem:breakpoint-product-tail} bounds its
conditional probability of being heavy by $\mathrm e^{-\gamma\xi m}$.
There are at most $n^m$ ordered choices of its vertices, and its
$m-1$ edges occur with probability $(d/n)^{m-1}$. Hence
\begin{align*}
    \Pr{\mathcal E^c}
    &\leq\sum_{m=h}^n nd^{m-1}\mathrm e^{-\gamma\xi m}
      \leq\frac{n}{d}
            \frac{(d\mathrm e^{-\gamma\xi})^h}
                 {1-d\mathrm e^{-\gamma\xi}}
      \leq2n\mathrm e^{-\gamma\xi h/2}
      \leq2n^{-31},
\end{align*}
where we used $\xi = d^{2\eta}$, the definition of $h$, and
$\gamma K=64$.
\end{proof}

\subsubsection{Non-breakpoint paths are non-elementary}

Call a path $L=(v_1,\ldots,v_h)$ \emph{elementary} if there is no
other self-avoiding path of at most $2h$ edges joining two distinct
vertices of $L$ with all internal vertices outside $L$. Hence, a non-elementary path $L$  contains an external path of length at most $2h$ between some pair of vertices on $L$. We show that every long, non-breakpoint path is non-elementary with high probability.

\begin{lemma}[Long non-breakpoint paths are non-elementary]\label{lem:breakpoint-elementary-paths}
    \begin{equation}
        \Prlr{\mathcal E\text{ and there is an elementary
            $h$-vertex path in }G[V\setminus\widehat S]}
        \leq n^{-7}.
    \end{equation}    
\end{lemma}
\begin{proof}

Fix an ordered list $L=(v_1,\ldots,v_h)$ of distinct vertices and
condition only on its path edges being present. For each $i$, define
\[
    Q_i \coloneqq \max_{ \substack{v_i, u_1, \ldots, u_\ell} }\prod_{j=1}^{\ell}W(u_j)
\] where the maximum is taken over all paths $(v_i,u_1,\ldots,u_\ell)$ of length at most $h$ and $u_j\notin V(L)$ for all $j$. Include the empty path
with value $1$ among these choices, so that $Q_i\geq1$. Intuitively, $Q$ is the maximum weight of any path of length $\le h$ starting at $v$ that does not intersect $L$.
Set 
$U_i=W(v_i)Q_i$ to include the weight of $v_i$ in this quantity.

Suppose that $\mathcal E$ holds, $L$ is elementary, and every vertex
of $L$ lies outside $\widehat S$. A heavy length $\le h$ path $P_i$ starting at $v_i$ that witnesses the non-breakpoint property of $v_i$ cannot
leave $L$ and return: returning to another vertex would violate
elementarity, and returning to the same vertex would violate
self-avoidance. Its vertices on $L$
therefore form an interval $J_i$ of indices containing $i$, followed
possibly by an external path. If this external path starts at $v_j$,
its weight excluding $v_j$ is at most $Q_j$. Thus,
\[
    1<W(P_i) \le \prod_{k\in J}U_k.
\]
The intervals $J_i$ for all $1 \le i \le |L|$ cover all vertices of $[h]$, but they are not disjoint. Therefore, we now choose a sub-family of the intervals $J_i$ that are pairwise disjoint.   

To this end, we repeatedly choose an interval $I$ of maximum length among the remaining intervals $J_i$
and discard all intervals that meet $I$. The intervals remaining after this procedure are clearly
disjoint. Moreover, their union covers at least $|L|/3$ vertices of $L$: every time we choose an interval of length $\ell$, every subsequently discarded interval has length at most $\ell$ too, so their union has length at most
$3\ell$. Since the original family covers $[h]$, the union
$I$ of the chosen intervals thus covers at least $h/3$ vertices.

Now, if $I\subseteq[h]$ be the union of our chosen intervals, denote for every $i \in I$ by $\widetilde{P}_i$ a path starting at $v_i$ such that $W(\widetilde{P}_i) = U_i$. Elementarity ensures that paths attached at different
vertices of $L$ have disjoint external vertex sets as an intersection
would give an external connection of at most $2h$ edges.

$L$ together with the chosen witnesses $\widetilde{P}_i$ therefore forms a tree of maximum degree $3$. Denoting by $\ell_i$ the number of external vertices used by $\widetilde{P}_i$, applying \Cref{lem:breakpoint-product-tail}, and bounding the number of all possible choices for the external path $\widetilde{P}_i$ times the probability that it exists by $n^{\ell_i} (d/n)^{\ell_i} = d^{\ell_i}$, we therefore get
\begin{align*}
    &\Prlr{L\text{ is elementary and }
        \prod_{i\in I}U_i>1\mid L\subseteq G }\\
    &\hspace{3cm}\leq
        \mathrm e^{-\gamma\xi |I|}
        \prod_{i\in I}\left(\sum_{\ell_i=0}^h
                            (d\mathrm e^{-\gamma\xi})^{\ell_i}\right)
        \leq(2\mathrm e^{-\gamma\xi})^{|I|},
\end{align*}
as the product is smaller than $2^{|I|}$ if $d$ is chosen large enough.
Since there are at most $2^h$ choices for $I$ and $|I| \ge h/3$, we thus get \begin{align*}
    &\Prlr{\mathcal E,\ L\text{ is elementary},
               V(L)\subseteq V\setminus\widehat S
               \mid L\subseteq G } \le 2^h(2\mathrm e^{-\gamma\xi})^{h/3}
       \leq\mathrm e^{-\gamma\xi h/4}.
\end{align*} if $d$ is large enough.
Finally, summing over all possible paths $L$ gives
\[
    \Prlr{\mathcal E\text{ and there is an elementary
            $h$-vertex path in }G[V\setminus\widehat S]} \le nd^{h-1}\mathrm e^{-\gamma\xi h/4}
    \leq n\mathrm e^{-\gamma\xi h/8}
    \leq n^{-7},
\]
as desired.
\end{proof}

\subsubsection{Small balls have at most one cycle}

Moreover, we show the following simple, first moment based lemma. 
\begin{lemma}[Balls have at most one cycle]\label{lem:balls}
    Fix any constant $c>20$ and let $R=\lceil ch\rceil$, so that $R>20h$.
    Denote by $B_G(x,R)=\{y:\operatorname{dist}_G(x,y)\leq R\}$ the ball of radius $R$ in $G$.
    Then, for $d$ sufficiently large depending on $\eta,A,C,c$,
    with probability $1-o(1)$ as $n\to\infty$, every induced ball
    $G[B_G(x,R)]$ is a tree or unicyclic.
\end{lemma}
\begin{proof}
    Recall that $\xi=d^{2\eta}$ and
    $h=\lceil K\log n/\xi\rceil$. Since $c$ is fixed and
    $\xi/\log(\mathrm e d)\to\infty$, we may take $d$ large enough that
    \[
        \frac{4cK\log(\mathrm e d)}{\xi}\leq\frac12.
    \]
    If one of these balls has at least two edges outside a breadth-first
    spanning tree rooted at $x$, choose two such edges. The tree paths
    from their four endpoints to $x$ have at most $R$ edges each.
    Their union is a tree, and adding the two chosen edges therefore
    produces a connected subgraph on $k \le 4R+1$ vertices with exactly
    $k+1$ edges. There are $k^{k-2}$ labeled trees on a specified set
    of $k$ vertices and at most $k^4$ choices for the two additional
    edges. A union bound therefore gives
    \begin{align*}
        &\Prlr{\text{some induced radius-$R$ ball has
                                at least two cycles}}\\
        &\hspace{1cm}\leq
        \sum_{k=4}^{4R+1}\binom nk k^{k-2}k^4(d/n)^{k+1}
        \leq\frac{d}{n}\sum_{k=4}^{4R+1}k^2(\mathrm e d)^k\\
        &\hspace{1cm}\leq
        \frac{d}{n}(4R+1)^3(\mathrm e d)^{4R+1}
        =n^{-1+4cK\log(\mathrm e d)/\xi+o_n(1)}
        \leq n^{-1/2+o_n(1)}
        =o(1).
    \end{align*}
    Here we used $R=cK\log n/\xi+O_c(1)$, and the limit is taken
    with $d$ and $c$ fixed. Thus the argument applies to every fixed
    multiple $c>20$, after choosing $d$ sufficiently large depending on $c$.
\end{proof}

\subsubsection{From paths to components}

Finally, we prove \Cref{thm:modified-breakpoint-components}.
\begin{proof}[Proof of \Cref{thm:modified-breakpoint-components}]
    We work on the simultaneous events that $\mathcal E$ holds, there is no elementary $h$-vertex path in $G[V\setminus\widehat S]$, and every induced radius-$24h$ ball has at most one cycle. The estimates from \Cref{lem:breakpoint-short-witnesses}, \Cref{lem:breakpoint-elementary-paths} and \Cref{lem:balls}, applied with $c=24$ and increasing $d_0$ if necessary, show that their intersection has probability $1-o(1)$.

    Suppose that $G[V\setminus\widehat S]$ contains a path $P$ on $8h$
    vertices, and let $v$ be its first vertex. By assumption, we know that $P$ is not elementary and that the induced graph
    \[
        H=G[B_G(v,10h)]
    \]
    contains $P$ and is a tree or unicyclic. Every external connection
    of at most $2h$ edges between vertices of a subpath of $P$ also
    lies in $H$: each of its vertices has distance at most
    $(8h-1)+2h<10h$ from $v$.

    If $H$ is a tree, then any $h$-vertex subpath of $P$ is elementary,
    because an additional external connection would create a cycle.

    We thus only have to consider the case where $H$ is unicyclic. To this end, let $C$ be its unique cycle. When $|V(C)|>3h$,
    any $h$-vertex subpath is again elementary: an external connection
    of at most $2h$ edges, together with the intervening subpath, would
    create a cycle of length at most $2h+(h-1)<3h$.

    It thus only remains to consider $|V(C)|\leq3h$. A self-avoiding path in a unicyclic graph has at most two portions outside the cycle. Indeed,
it cannot enter a tree attached to the cycle and later return to
the cycle without revisiting the attachment vertex. A path disjoint
from the cycle has just one such portion. At least $8h-3h=5h$
vertices of $P$ lie outside $C$, so one of these portions contains
an $h$-vertex subpath $L$. This $L$ is elementary, since an external
connection would create a cycle containing vertices outside $C$,
contradicting uniqueness of the cycle.

In every case, $P$ contains an elementary $h$-vertex path. This is
a contradiction, proving our diameter bound of $8h$.

Moreover, a shortest path inside
any component has fewer than $8h$ vertices. Thus the entire
component lies in $B_G(u,8h)$ for any one of its vertices $u$.
This ball has at most one cycle, so the connected induced component
also has at most one cycle, which proves that every block is a tree or unicyclic.

Finally, on occurrence of $\mathcal E$, every vertex outside $\widehat S$ has a
heavy witness on fewer than $h$ vertices, and hence on fewer than
$h$ edges. Every test with radius $r\geq h$ detects that witness.
Conversely, a global breakpoint passes every finite-radius test.
Therefore $\widehat S_r=\widehat S$ simultaneously for all integers
$r\geq h$, proving the third property.

For fixed $d$ and sufficiently large $n$, we further have $\log n/ \xi \ge 1$, and
thus
\[
    8h\leq8(K+1)\frac{\log n}{d^{2\eta}},
\]
which verifies the asserted independence of $d$ of the hidden constant in our diameter bound and thus
completes the proof.
\end{proof}

\subsection{Merging components to satisfy unique attachment ports}
\label{sec:merging}

With \Cref{thm:modified-breakpoint-components}, it is not hard to deduce our final desired theorem showing that a decomposition consistent with \Cref{def:validblockpartion} does exist with high probability.

\begin{proof}[Proof of \Cref{thm:validblockpartition}]
    \Cref{thm:modified-breakpoint-components} already establishes that every block resulting after removing the vertices in $\widehat{S}$ is a tree or unicyclic, and that each such block has diameter at most $8h$. It only remains to be shown that we can choose $S$ such that every boundary vertex $v \in S$ is only adjacent to at most one vertex in any block $B \in \mathcal{B}$ resulting after removing all vertices in $S$. 

    To show this, start with $S = \widehat{S}$. Assume that there is a vertex $v \in S$ and a neighboring block $B$ such that $v$ has two neighbors in $B$. We call such vertices \emph{violating}. Let $\mathcal{B}_v$ be the set of blocks adjacent to $v$. Then we merge $v \cup \bigcup_{B \in \mathcal{B}_v}$ into a new block $B'$. This procedure is equivalent to removing $v$ from our set of singletons $S$ and considering the resulting blocks. We repeat this procedure until no singletons with more than one neighbor in the same block are left. 

    Work on the event from \Cref{thm:modified-breakpoint-components}
    and the event that every induced ball of radius $24h$ has at most
    one cycle. By \Cref{lem:balls}, applied with $c=24$, their
    intersection has probability $1-o(1)$ after increasing $d_0$ if
    necessary. Call the components of $G[V\setminus\widehat S]$
    the original blocks. By the first property of
    \Cref{thm:modified-breakpoint-components}, an original block
    contains no self-avoiding path on $8h$ vertices, and hence has
    diameter at most $8h-2$.

    We show by induction that every block created by the merging
    procedure is unicyclic and has diameter less than $16h$.
    In particular, every tree block present at any step is an
    original block. Suppose that the assertion holds for all
    previous merges, and let $v$ be the next violating vertex.
    All current blocks have diameter less than $16h$, so the
    connected induced graph
    \[
        H=G\left[\{v\}\cup\bigcup_{B\in\mathcal B_v}B\right]
    \]
    is contained in $G[B_G(v,16h)]$: every vertex in an adjacent
    block can be reached from $v$ by one edge followed by a path
    of length at most the diameter of that block.
    Thus $H$ has at most one cycle.

    For each $B\in\mathcal B_v$, write
    \[
        r_B=|E(v,B)|,\qquad
        \ell_B=|E(G[B])|-|B|+1\in\{0,1\}.
    \]
    There are no edges between distinct current blocks. Therefore
    the number of edges of $H$ outside a spanning tree is
    \[
        |E(H)|-|V(H)|+1
          =\sum_{B\in\mathcal B_v}(\ell_B+r_B-1)\leq1.
    \]
    All summands are nonnegative, and $v$ is violating, so at least
    one $r_B$ is at least two. It follows that the sum equals one,
    every $\ell_B$ is zero, exactly one block has $r_B=2$, and
    every other adjacent block has $r_B=1$.
    In particular, all blocks absorbed at this step are trees,
    and hence are original blocks. This is why merges cannot
    keep absorbing previously merged blocks.

    The new block $B'$ is consequently unicyclic. Moreover, every
    $x\in B'$ satisfies
    \[
        \operatorname{dist}_{G[B']}(v,x)
          \leq1+(8h-2)=8h-1.
    \]
    Joining two such paths through $v$ gives
    \[
        \operatorname{diam}(G[B'])\leq16h-2<16h,
    \]
    which completes the induction.

    Each merge removes one vertex from $S$, so the procedure
    terminates after at most $|\widehat S|$ steps.
    At termination, every block is a tree or unicyclic, and no
    singleton has more than one neighbor in any block.
    Also, $S\subseteq\widehat S$, so all remaining singletons are
    still global breakpoints; their weights and degrees are
    measured in the original graph and do not change during merging.

    Finally, the choice $K=64/\gamma$ above depends only on $A$.
    Set $D=16(K+1)$, independently of the weight constant $C$.
    For fixed $d$ and sufficiently large $n$, we have
    $\log n/d^{2\eta}\geq1$, and therefore
    \[
        \operatorname{diam}(G[B])<16h
          \leq D\frac{\log n}{d^{2\eta}}.
    \]
    This ensures the precise diameter condition in
    \Cref{def:validblockpartion}, as well as the asserted
    asymptotic bound, for every fixed $C>0$ once $d$ is sufficiently
    large. All conditions for a valid block partition
    now hold, which proves the theorem.
\end{proof}

%% file: comparison.tex
The goal of this section is to prove \Cref{thm:main} by comparing
single-site Glauber dynamics with the block dynamics studied above.
The comparison essentially amounts to estimating the relaxation time inside a block, which we control by
recursively separating a block into pieces whose boundary vertices
are low-degree. This lets us handle hubs without assuming that their
color marginals are close to uniform and constitutes the main difference to the comparison theorems in \cite{efthymiou2018sampling, li2026rapid}.

We first prove the within-block relaxation time estimate, then show how it lets us compare Dirichlet forms, and finally choose the right parameters to yield a proof of
\Cref{thm:main}. 

\subsection{The within-block relaxation time bound}

Recall the parameters $A, C, \eta, d$ governing our block partition from \Cref{sec:blockassumptions} and \Cref{sec:blocks}.
Throughout this section, we work in the regime
\[
    A\geq1,\qquad C>0,\qquad
    0<\eta<\varepsilon\leq\tfrac12,\qquad
    q\geq d+d^{1/2+\varepsilon},
\]
and take $d$ sufficiently large depending on these fixed parameters.
As before, put
\[
    \lowdegree=d+d^{1/2+\eta},\qquad
    \alpha=C\log d,\qquad a=q-\lowdegree.
\]
We take $d$ large enough that $\alpha>1$, $q\geq7$, and $a = q - \lowdegree \geq20$.
Call a vertex a hub if its degree in the original graph $G$ exceeds
$\lowdegree$, and write $\mathcal H$ for the set of hubs. All degrees
and the distinction between hubs and low-degree vertices remain those
of $G$, even after conditioning or passing to a subtree. For now, we work with continuous time Glauber dynamics defined as follows.

\begin{definition}[Continuous-time Glauber dynamics]
\label{def:withinblockglauber}
For a block $B$ and a feasible coloring $\taubf$ of $V\setminus B$,
let $\nu_{B,\taubf}$ be the uniform measure on proper extensions to $B$.
For $v\in B$, let $P_v^{B,\taubf}f (\sigmabf) = \Expected{f(\sigmabf) \mid \sigmabf_{V \setminus B} = \taubf, \sigmabf_{B \setminus \{v\}}}$ be the single-site heat-bath update
and define the positive generator
\[
    L_{\mathrm{Gl}}^{B,\taubf}
       =\sum_{v\in B}(I-P_v^{B,\taubf}).
\]
This is the continuous-time version of Glauber dynamics where every vertex is updated at rate $1$.
\end{definition}
We are concerned with estimating the within-block relaxation time defined below.
\begin{definition}[Within-block relaxation time]
\label{def:withinblockrelaxation}
Write
$\Varsub{v}{f}$ for the conditional variance given all
variable colors except the one at $v$, with the exterior coloring
fixed. We define $\tau_{\mathrm{rel}}(B,\taubf)$ as the smallest constant
$T\geq1$ such that
\[
    \Varsub{\nu_{B,\taubf}}{f}
       \leq T\sum_{v\in B}
                    \Expectedsub{\nu_{B,\taubf}}{\Varsub{v}{f}}
    \qquad\text{for every }f.
\]
The value is infinite if this inequality fails. This convention also
covers a conditional state space consisting of one coloring.
\end{definition}

The main theorem of this subsection is the following.

\begin{theorem}[Uniform within-block relaxation]
\label{thm:withinblockcomparison}
    There are absolute constants $c_0,c_1>0$ and an event of
    probability $1-o(1)$ over $G\sim G(n,d/n)$ on which the following
    holds. For every tree or unicyclic block $B$
    whose exterior boundary consists of global breakpoints, and every
    feasible exterior coloring $\taubf$, writing
    $b_B=\operatorname{diam}(G[B])$, we have
\begin{equation}
        \tau_{\mathrm{rel}}(B,\taubf)
        \leq
        \exp\left\{c_0\left[
          \left(b_B+1+\frac{\log n}{d^{2\eta}}\right)\log q
                    +\frac{\log n}{q}\right]\right\}.
        \label{eq:localrelaxationgeneral}
    \end{equation}
    In particular, for any fixed diameter constant $D$ and all
    sufficiently large $n$, blocks satisfying
    \[
        b_B\leq D\frac{\log n}{d^{2\eta}}
    \]
    satisfy, simultaneously over $B,\taubf$,
\begin{equation}
        T_*:=\max\left\{1,\sup_{B,\taubf}
                         \tau_{\mathrm{rel}}(B,\taubf)\right\}
            \leq n^{\zeta_d},
        \qquad
        \zeta_d=c_1(D+1)\frac{\log q}{d^{2\eta}}.
        \label{eq:withinblockloss}
    \end{equation}
    All these conclusions apply to the partition guaranteed by
    \Cref{thm:validblockpartition}, with $D=16(K+1)$.
    Here $d$ is fixed before taking $n\to\infty$, and the event and
    the lower threshold on $n$ can be chosen independently of $q$.
    In particular, if $q=\Theta(d)$, then
    $\zeta_d=O_D(\log d/d^{2\eta})=o_d(1)$.
\end{theorem}

To prove the above, we first devise a recursive decomposition of a tree/unicyclic block into separators. Then we use this hierarchy of separators to prove a deterministic recursive relaxation time bound. Finally, we bound
the cost of this deterministic estimate using properties that hold simultaneously over all connected
tree/unicyclic vertex sets of $G$. This uniformity is important since
the partition in \Cref{thm:validblockpartition} is chosen after $G$
has been sampled.

\subsubsection{Constructing separators recursively}

We proceed by devising a simple recursive decomposition of a tree/unicyclic block into recursive separators. Each separator is constructed such that after conditioning on its colors, the boundary condition on each remaining component is \emph{admissible}, in the following sense.


For any subset of vertices $T\subseteq V$, a boundary condition
fixes the colors on $V\setminus T$ and is \emph{feasible} if it
extends to a proper coloring of $G$. A feasible boundary condition
is called \emph{admissible}
if every $w\in T$ with a neighbor in $V\setminus T$ satisfies
$\deg_G(w)\leq\lowdegree$. Equivalently, no hub in $T$ has a
fixed exterior neighbor. This is a restriction on which vertices
touch the exterior, not on their particular colors, and all degrees
refer to the original graph $G$.
Admissibility ensures that all vertices touching the exterior are
covered by the marginal upper bound stated below.

\begin{lemma}[Admissible boundaries and low-degree marginals]
\label{lem:comparisonboundary}
    Suppose $A\geq1$. A block whose exterior boundary consists
    of global breakpoints has an admissible boundary condition.
    Moreover, for every subset $T\subseteq V$, every feasible
    coloring $\taubf$ of $V\setminus T$, and every $w\in T$
    with $\deg_G(w)\leq\lowdegree$, the uniform measure on
    proper $q$-colorings satisfies
\begin{equation}
        \Prlr{\sigmabf_w=c\mid\sigmabf_{V\setminus T}=\taubf}
            \leq\frac1a
        \qquad\text{for every }c\in[q],
        \label{eq:comparisonlowmarginal}
\end{equation} where $a = q - \lowdegree$. This marginal bound
    does not require the boundary condition to be admissible.
\end{lemma}
\begin{proof}
    A global breakpoint $s$ is low-degree by the test on the
    single-vertex path. If a neighbor $v$ of $s$ were a hub, then
    \[
        W(s)W(v)=\alpha^{A-1}\frac{\deg_G(v)}d>1,
    \]
    contrary to the breakpoint test on $(s,v)$. Thus every
    vertex in the block having an exterior neighbor is
    low-degree.

    For the marginal bound, first fix
    $\sigmabf_{V\setminus T}=\taubf$ and then reveal all colors
    in $T\setminus\{w\}$. The law at $w$ is uniform on at least
    $q-\deg_G(w)\geq a$ available colors. Averaging over the
    additional conditioning proves
    \eqref{eq:comparisonlowmarginal}.
\end{proof}

To construct separators
starting from a chosen root (or the entire cycle in the unicyclic
case), we repeatedly add every hub adjacent to a vertex already in
the separator, until no such hub remains outside it. Thus we include
whole connected components of hubs attached to the root or cycle,
not just its immediate hub neighbors. We do not continue through a
low-degree vertex outside the initial root or cycle: these vertices
remain in the pieces left behind and become their low-degree roots.

\begin{lemma}[Recursive separator decomposition]
\label{lem:comparisonseparators}
    Let $T\subseteq V$ be such that $G[T]$ is a tree or is
    unicyclic, and fix an admissible boundary condition on $T$.
    There is a connected
    separator $\mathcal C$ with the following properties.
    Enumerate its vertices as $\mathcal C=\{u_1,\ldots,u_k\}$.
    Every component $T_i$ of $T\setminus\mathcal C$ is a tree
    attached to $\mathcal C$ by exactly one edge
    $\{u_{\mathrm{par}_{\mathcal C}(i)},w_i\}$ for some index $\mathrm{par}_{\mathcal C}(i)\in[k]$,
    where $w_i\in T_i$ is low-degree.
    Conditioning on any feasible coloring of $\mathcal C$
    leaves an admissible boundary condition on each $T_i$.

    Recursing on the $T_i$, rooted at $w_i$, partitions $T$
    into separators whose adjacency graph is a rooted tree.
    Every noninitial separator contains exactly one
    low-degree vertex, its root.
\end{lemma}
\begin{proof}
    If $T$ is a tree, choose a root $r$ and let $\mathcal C$
    consist of $r$ and every hub reachable from $r$ along a
    path whose vertices other than $r$ are hubs. If $T$ is
    unicyclic, start with its entire cycle and add every hub
    reachable from it along a path whose vertices outside
    the cycle are hubs.

    The separator is connected and contains the cycle when
    there is one. Each remaining component is therefore a
    tree attached by exactly one edge. Its attachment vertex
    $w_i$ must be low-degree, since otherwise it would have
    been included in $\mathcal C$. The only new fixed
    exterior neighbor of $T_i$ is at $w_i$, so admissibility
    is preserved.

    Repeat the construction with $T_i$ rooted at $w_i$.
    Each later separator consists of its low-degree root
    and all hubs reachable from it along paths containing
    only hubs after the root. Disjointness, termination, and the
    asserted rooted-tree structure follow from the
    one-edge attachment property.
\end{proof}

\subsubsection{A one-level comparison}

The next lemma compares allows us to state the relaxation time in a block as a function of the relaxation times of the children $T_i$ associated to a separator $\mathcal{C}$. 

To introduce it, fix a subset $T\subseteq V$ such that $G[T]$ is a tree or is
unicyclic, a feasible admissible coloring $\taubf$ of $V\setminus T$,
and a separator $\mathcal C$ as supplied by
\Cref{lem:comparisonseparators}. Define
\[
    \nu\coloneqq\nu_{T,\taubf}
       =\mu_{q,G}\bigl(\sigmabf_T\in\cdot\mid
                         \sigmabf_{V\setminus T}=\taubf\bigr)
\]
as the uniform measure over colorings of $T$ conditional on
$\taubf$. Recall that
\[
    \mathcal H=\{v\in V:\deg_G(v)>\lowdegree\}
\]
is the set of hubs, with degrees measured in the original graph.
Enumerate the separator vertices as
$\mathcal C=\{u_1,\ldots,u_k\}$, where $k=|\mathcal C|$;
the order will be chosen in the proof below.
Let $T_1,\ldots,T_m$ be the connected components of
$G[T\setminus\mathcal C]$. For each $i\in[m]$, let $\mathrm{par}_{\mathcal C}(i)\in[k]$
be the index of its attachment vertex in $\mathcal{C}$: the unique edge joining
$T_i$ to the separator is $\{u_{\mathrm{par}_{\mathcal C}(i)},w_i\}$, where
$w_i\in T_i$ is the low-degree root of the child component.
The child index $i\in[m]$ and separator index $\mathrm{par}_{\mathcal C}(i)\in[k]$
are distinct types of indices; different children may have
the same attachment vertex.

For each $v\in\mathcal C$, define
\[
    R_v\coloneqq\bigl|\{i\in[m]:u_{\mathrm{par}_{\mathcal C}(i)}=v\}\bigr|
       =\bigl|N_G(v)\cap(T\setminus\mathcal C)\bigr|
\]
as the number of child components at separator vertex $v$.
Furthermore, define
\[
    \ell=|\mathcal C\setminus\mathcal H|,\qquad
    r_{\mathcal H}=\sum_{v\in\mathcal C\cap\mathcal H}R_v
\]
such that $\ell$ counts the low-degree vertices in our separator and $r_{\mathcal{H}}$ is the number of children attached to hubs in $\mathcal{C}$.
We have $\ell=1$ for every noninitial separator, while the initial
separator satisfies only $\ell\leq|\mathcal C|$, since in the
unicyclic case it may contain many low-degree cycle vertices.

For a function $f$ on the colorings of $T$, write
\begin{align*}
    \mathcal E_0(f)
       &\coloneqq\Expectedsub{\nu}
            {\Varsub{\nu}{f\mid\sigmabf_{\mathcal C}}},\\
    \mathcal E_v(f)
       &\coloneqq\Expectedsub{\nu}
            {\Varsub{\nu}{f\mid\sigmabf_{T\setminus\{v\}}}},
       \qquad v\in T.
\end{align*}

To estimate the relaxation time within a block, we prove the following estimate that lets us aggregate the relaxation time on the children $T_i$ to get a bound on the relaxation time in $T$. Applying it repeatedly will ultimately lead to \Cref{thm:withinblockcomparison}
\begin{lemma}[One-level separator comparison]
\label{lem:separatorcomparison}
    Let $T$, $\mathcal C$, and $T_1,\ldots,T_m$ be as above.
    For every feasible admissible exterior coloring
    $\taubf\in[q]^{V\setminus T}$, set
    $\nu=\nu_{T,\taubf}$. Then, for every
    real-valued function $f$ on the proper extensions of
    $\taubf$ to $T$,
\begin{equation}
        \Varsub{\nu}{f}
          \leq K_{\mathcal C}
              \left(\mathcal E_0(f)
                      +\sum_{v\in\mathcal C}\mathcal E_v(f)\right)
        \qquad \text{with} \qquad
        K_{\mathcal C}
          =3(2k+3)q^k
                      \exp\left(4\ell+\frac{28r_{\mathcal H}}q\right).
        \label{eq:separatorpoincare}
    \end{equation}
    Consequently, suppose that for each $i\in[m]$ there is
    $F_i\geq1$ such that $\tau_{\mathrm{rel}}(T_i,\xi)\leq F_i$
    for every feasible exterior coloring
    $\xi\in[q]^{V\setminus T_i}$ satisfying
    $\xi_{V\setminus T}=\taubf$. Then
\begin{equation}
        \tau_{\mathrm{rel}}(T,\taubf)
             \leq K_{\mathcal C}\max\{1,\max_i F_i\}.
        \label{eq:separatorrecurrence}
    \end{equation}
\end{lemma}
\begin{proof}
    Our goal is to construct a coupling of a Markov chain for colorings of $T$ with the following update operators $P_0, \ldots P_k$ that all keep
    $\sigmabf_{V\setminus T}=\taubf$ fixed.
    Let $P_0$ resample all children together, conditional on
    $\sigmabf_{\mathcal C}$. The children are conditionally
    independent. Order the vertices of $\mathcal C$ as $u_1,\ldots,u_k$ so that
    each vertex $u_i$ has at most two neighbors in $G[\mathcal C]$ that come after $u_i$ in this sequence.
    Such an ordering exists because this graph is a tree or
    unicyclic: we can construct it by repeatedly removing a vertex of
    degree at most two from the remaining induced subgraph and
    listing the vertices in their removal order.
    Let $P_j$ be the heat-bath update at $u_j$.

    Rather than choosing the updates in a random order, we
    study an auxiliary chain that applies them in a
    carefully chosen order such that controlling the success probability of our coupling becomes easier. Our goal is to construct a coupling for this Markov chain which afterwards helps us establish \eqref{eq:separatorpoincare}.

    To this end we introduce a so-called \emph{sweep} by which we mean the fixed sequence of updates
    $P_0,P_1,\ldots,P_k,P_0$: first resample all child
    components, then update each separator vertex once in
    the order $u_1,\ldots,u_k$, and finally resample all
    children again. This is an auxiliary update sequence
    used for comparison, not a single step of Glauber dynamics.

    With this, the proof has three steps:
    \begin{enumerate}
        \setlength{\itemsep}{0cm}
        \item Couple one sweep of updates from any two
            initial configurations, so that the final
            configurations agree with probability at least
            $\epsilon_{\mathcal C}$, uniformly over the
            initial configurations.
        \item Append the updates in reverse order to obtain
            a reversible chain $Q$. The coupling then gives
            $\operatorname{gap}(I-Q)\geq\epsilon_{\mathcal C}$.
        \item Compare the Dirichlet form of $Q$ with those
            of its individual updates. This gives
            \eqref{eq:separatorpoincare}, with
            $K_{\mathcal C}=3(2k+3)/\epsilon_{\mathcal C}$.
    \end{enumerate}
    Finally, we apply the Poincare inequalities in the child
    components to obtain the recursive bound
    \eqref{eq:separatorrecurrence}. The auxiliary chain
    resamples these components exactly; their replacement
    by single-site updates happens in this last step.

    \paragraph{Step 1: coupling a sweep}
    Starting from any two configurations $X,Y$ consistent
    with $\taubf$, we jointly sample the two sweeps so that
    each copy has the prescribed heat-bath transitions.
    Our goal is to make the final configurations agree on
    all of $T$ with a probability bounded below uniformly
    over $X,Y$. We first try to bring both separators to
    the same target coloring. Once they agree, the final
    $P_0$ update can use identical samples in both copies,
    since their conditional child distributions then coincide.
    The resulting agreement probability, quantified below
    by $\epsilon_{\mathcal C}$, will give a spectral-gap
    bound for a reversible version of the sweep.

    For each child $i\in[m]$, let $p_i^X,p_i^Y$ be the
    distributions of the color at $w_i$ after the first
    $P_0$ update in the two copies. Explicitly, for $c\in[q]$,
    \[
        p_i^X(c)\coloneqq
            \nu(\sigmabf_{w_i}=c\mid\sigmabf_{\mathcal C}=X_{\mathcal C}),
    \]
    and analogously with $Y$ in place of $X$. For any palette
    $\mathcal P\subseteq[q]$, write
    \[
        p_i^X(\mathcal P)\coloneqq\sum_{c\in\mathcal P}p_i^X(c)
           =\nu(\sigmabf_{w_i}\in\mathcal P
                          \mid\sigmabf_{\mathcal C}=X_{\mathcal C}),
    \]
    and define $p_i^Y(\mathcal P)$ analogously. Each atom $p_i^X(c),p_i^Y(c)$ is at most $1/a$ by
    \eqref{eq:comparisonlowmarginal}.
    Before sampling the children, choose a seven-element palette
    $\mathcal P\subseteq[q]$ such that
\begin{equation}
        \sum_{\substack{i\in[m]\\u_{\mathrm{par}_{\mathcal C}(i)}\in\mathcal H}}
            \left(p_i^X(\mathcal P)+p_i^Y(\mathcal P)\right)
                      \leq\frac{14r_{\mathcal H}}q.
        \label{eq:paletteload}
    \end{equation}
    To justify this choice, fix $X,Y$ and temporarily choose
    $\mathcal P_{\mathrm{rnd}}$ uniformly among the seven-element
    subsets of $[q]$. Each color $c$ belongs to this random
    palette with probability
    \[
        \Prlr{c\in\mathcal P_{\mathrm{rnd}}}
            =\frac{\binom{q-1}{6}}{\binom q7}=\frac7q.
    \]
    Therefore, since $\sum_{c\in[q]}p_i^X(c)=1$,
    \[
        \Expectedsub{\mathcal P_{\mathrm{rnd}}}
            {p_i^X(\mathcal P_{\mathrm{rnd}})}
          =\sum_{c\in[q]}p_i^X(c)
                        \Prlr{c\in\mathcal P_{\mathrm{rnd}}}
          =\frac7q,
    \]
    and the same holds for $p_i^Y$. There are exactly
    $r_{\mathcal H}$ children attached to hubs, so linearity
    of expectation gives expected total load
    $r_{\mathcal H}(7/q+7/q)=14r_{\mathcal H}/q$.
    At least one palette has load no greater than this average,
    proving \eqref{eq:paletteload}. We fix such a palette
    before sampling the children; it may depend on $X,Y$.

    We now construct $Z$, a common target coloring of the
    separator $\mathcal C$ only. All its colors are chosen
    before sampling the children, using the fixed initial
    configurations $X,Y$, the exterior coloring $\taubf$,
    and the palette $\mathcal P$. We assign the targets in
    reverse order: for $j=k,k-1,\ldots,1$, set $v=u_j$.
    At this point, every $Z_{u_t}$ with $t>j$ has already
    been chosen. Define the forbidden set
    \[
        \mathcal F_v\coloneqq
           \bigcup_{\substack{t>j\\\{v,u_t\}\in E(G)}}
                      \{X_{u_t},Y_{u_t},Z_{u_t}\}.
    \]
    Here ``later'' refers to the forward update order
    $u_1,\ldots,u_k$, even though these vertices are processed
    first when constructing $Z$. There are at most two such
    neighbors, so $|\mathcal F_v|\leq6$.
    Excluding their target colors will make $Z$ proper;
    excluding their initial colors will make the later
    forward recoloring sequence legal in both copies.

    If $v$ is a hub, choose $Z_v$ to be the smallest color
    in $\mathcal P\setminus\mathcal F_v$. This set is nonempty,
    since $|\mathcal P|=7>6$. Admissibility ensures that $v$
    has no neighbor outside $T$, so there is no further
    constraint from the fixed exterior coloring. Constraints
    from its children will be enforced by the clearing event.

    For a low-degree separator vertex $v$, we choose $Z_v$ differently and define
    \[
        \mathcal A_v\coloneqq
            [q]\setminus\{\taubf_x:x\in N_G(v)\setminus T\}
    \]
    as precisely the colors compatible with the fixed
    coloring $\taubf$ outside $T$. In particular, $\mathcal A_v$
    does not yet impose any constraint from neighbors in the
    separator or in its child components. It depends only on $\taubf$, not on the starting
    configurations $X,Y$. The separator constraints are handled
    by excluding the colors described above, and the child
    constraints by the clearing event below.

    The vertex $v$ has $\deg_{G[\mathcal C]}(v)$ neighbors in the
    separator and $R_v$ neighbors in $T\setminus\mathcal C$.
    Every other neighbor lies outside $T$ and forbids at most
    one color. Hence
    \begin{align*}
        |\mathcal A_v|
          &\geq q-|N_G(v)\setminus T|\\
          &=q-\bigl(\deg_G(v)-\deg_{G[\mathcal C]}(v)-R_v\bigr)\\
          &\geq a+\deg_{G[\mathcal C]}(v)+R_v,
    \end{align*}
    where the last step uses $\deg_G(v)\leq\lowdegree$
    and $a=q-\lowdegree$.
    Thus the candidate set $\mathcal A_v\setminus\mathcal F_v$
    has size at least
    $|\mathcal A_v|-6\geq a+R_v-6\geq R_v+14>0$.
    For each color $c$, define its load from the children of $v$ by
    \[
        L_v(c)\coloneqq
          \sum_{\substack{i\in[m]\\
                  u_{\mathrm{par}_{\mathcal C}(i)}=v}}
                         \bigl(p_i^X(c)+p_i^Y(c)\bigr).
    \]
    Choose
    \[
        Z_v\in\operatorname*{arg\,min}_{
                     c\in\mathcal A_v\setminus\mathcal F_v}L_v(c),
    \]
    breaking ties by taking the smallest color. This choice
    favors colors that the child roots are unlikely to use.
    Since each $p_i^X,p_i^Y$ is a probability distribution,
    $\sum_{c\in[q]}L_v(c)=2R_v$. Averaging over the candidate
    set therefore gives
\begin{equation}
        \sum_{\substack{i\in[m]\\u_{\mathrm{par}_{\mathcal C}(i)}=v}}(p_i^X(Z_v)+p_i^Y(Z_v))
            \leq\frac{2R_v}{|\mathcal A_v|-6}\leq2.
        \label{eq:lowtargetload}
    \end{equation}
    The resulting target $Z$ is proper and, at each vertex, avoids the
    initial colors of its later neighbors.

    Now, we construct the coupling. To this end, sample the first $P_0$ update independently in the two copies.
    Require every child attached to a hub to have its root avoid
    $\mathcal P$ in both copies. For a child attached to a
    low-degree vertex $v$, require both roots to avoid $Z_v$.
    More formally, let the forbidden colors for child $i$ be
    \[
        \mathcal D_i\coloneqq
        \begin{cases}
            \mathcal P,&u_{\mathrm{par}_{\mathcal C}(i)}\in\mathcal H,\\
            \{Z_{u_{\mathrm{par}_{\mathcal C}(i)}}\},
                &u_{\mathrm{par}_{\mathcal C}(i)}\notin\mathcal H,
        \end{cases}
        \qquad b_i^W\coloneqq p_i^W(\mathcal D_i),\quad W\in\{X,Y\}.
    \]
    Write $\mathcal E_{\mathrm{clr}}$ for the event that every
    child root avoids its forbidden set in both copies.
    With $X,Y,\mathcal P,Z$ fixed, the children are independent
    within each copy, and the two copies are sampled independently.
    Consequently,
    \[
        \Prlr{\mathcal E_{\mathrm{clr}}}
           =\prod_{i=1}^m(1-b_i^X)(1-b_i^Y).
    \]
    Equation~\eqref{eq:paletteload} bounds the sum of the
    forbidden masses for children attached to hubs by
    $14r_{\mathcal H}/q$. For each low-degree separator vertex
    $v$, equation~\eqref{eq:lowtargetload} bounds the corresponding
    sum by $2$. Since there are $\ell$ such vertices,
    \[
        \sum_{i=1}^m(b_i^X+b_i^Y)
            \leq\frac{14r_{\mathcal H}}q+2\ell.
    \]
    Moreover, every $b_i^W$ is at most $7/a\leq1/2$:
    a palette contains seven colors, each of mass at most $1/a$,
    while a set containing one color has mass at most $1/a$. Applying
    $\log(1-x)\geq-2x$ to each factor of the product gives
    \begin{align*}
        \log\Prlr{\mathcal E_{\mathrm{clr}}}
          &=\sum_{i=1}^m\bigl(\log(1-b_i^X)+\log(1-b_i^Y)\bigr)\\
          &\geq-2\sum_{i=1}^m(b_i^X+b_i^Y)\\
          &\geq-\frac{28r_{\mathcal H}}q-4\ell.
    \end{align*}
    Exponentiating yields
    $\Prlr{\mathcal E_{\mathrm{clr}}}
       \geq\exp(-28r_{\mathcal H}/q-4\ell)$.
    On that event, recoloring $u_1,\ldots,u_k$ successively
    to their target colors is legal in both copies. Earlier
    neighbors already have their proper target colors, later
    neighbors retain their excluded initial colors, and the
    children and fixed exterior neighbors permit the targets.
    At each step a maximal coupling of the two uniform
    heat-bath laws assigns their common target in both copies
    with probability at least $1/q$.
    Thus the separator agrees after $P_1,\ldots,P_k$ with
    probability at least
    \[
        \epsilon_{\mathcal C}
          =q^{-k}\exp\left(-4\ell-\frac{28r_{\mathcal H}}q\right).
    \]
    The following $P_0$ update can then couple the entire
    configurations identically.

    \paragraph{Step 2: a reversible chain and its spectral gap}
    To turn this sweep into a reversible comparison, consider
    \[
        Q=P_0P_1\cdots P_kP_0P_k\cdots P_1P_0.
    \]
    Each factor is a self-adjoint conditional-expectation
    projection on $L^2(\nu)$. Setting
    $M=P_0P_1\cdots P_kP_0$, we have $Q=MM^*$,
    since $P_0^2=P_0$. Thus $Q$ is reversible and positive
    semidefinite. The coupling just constructed coalesces by
    the middle $P_0$ with probability at least
    $\epsilon_{\mathcal C}$; subsequent updates preserve
    agreement. Consequently,
    \[
        \sup_{X,Y}\|Q(X,\cdot)-Q(Y,\cdot)\|_{\mathrm{TV}}
                       \leq1-\epsilon_{\mathcal C}.
    \]
    In particular, the oscillation of $Qg$ is at most
    $(1-\epsilon_{\mathcal C})$ times that of $g$.
    For any nonconstant eigenfunction $Qg=\lambda g$, we have
    $\operatorname{osc}(g)>0$, so the contraction gives
    $|\lambda|\leq1-\epsilon_{\mathcal C}$ after dividing by
    $\operatorname{osc}(g)$.
    Thus the eigenvalue $1$ of $Q$ has only constant eigenfunctions,
    and every nonzero eigenvalue of $I-Q$ is at least
    $\epsilon_{\mathcal C}$, proving
    $\operatorname{gap}(I-Q)\geq\epsilon_{\mathcal C}$.

    \paragraph{Step 3: comparison of Dirichlet forms}
    Start with $X_0\sim\nu$ and run the $2k+3$ successive
    updates defining $Q$, obtaining $X_1,\ldots,X_{2k+3}$.
    Each update preserves $\nu$, so every $X_j$ has law $\nu$.
    In particular, expanding the square gives
    \begin{align*}
        \tfrac12\Expected{(f(X_{2k+3})-f(X_0))^2}
          &=\Expectedsub{\nu}{f^2}
               -\Expected{f(X_0)f(X_{2k+3})}\\
          &=\langle f,(I-Q)f\rangle_\nu.
    \end{align*}
    Write $\Delta_j=f(X_j)-f(X_{j-1})$. The intermediate
    terms cancel in the telescoping sum, giving
    \[
        f(X_{2k+3})-f(X_0)=\sum_{j=1}^{2k+3}\Delta_j.
    \]
    Applying Cauchy--Schwarz to the vectors
    $(1,\ldots,1)$ and $(\Delta_1,\ldots,\Delta_{2k+3})$
    gives, for every realization of the updates,
    \[
        \left(\sum_{j=1}^{2k+3}\Delta_j\right)^2
           \leq\left(\sum_{j=1}^{2k+3}1^2\right)
                \left(\sum_{j=1}^{2k+3}\Delta_j^2\right)
           =(2k+3)\sum_{j=1}^{2k+3}\Delta_j^2.
    \]

    To evaluate the expected squared increments, suppose
    step $j$ resamples a set $U\subseteq T$. Conditional on
    the unchanged colors
    $\boldsymbol{\eta}=(X_{j-1})_{T\setminus U}=(X_j)_{T\setminus U}$,
    the old and new configurations have independent,
    identical laws under $\nu(\cdot\mid\sigmabf_{T\setminus U}=\boldsymbol{\eta})$.
    Thus,
    \begin{align*}
        &\Expected{(f(X_j)-f(X_{j-1}))^2\mid\boldsymbol{\eta}}
          \\ &\hspace{2cm}=\Expected{f(X_j)^2\mid\boldsymbol{\eta}}
               +\Expected{f(X_{j-1})^2\mid\boldsymbol{\eta}} -2\Expected{f(X_j)\mid\boldsymbol{\eta}}
               \Expected{f(X_{j-1})\mid\boldsymbol{\eta}}\\
          &\hspace{2cm}=2\Varsub{\nu}{f\mid\sigmabf_{T\setminus U}=\boldsymbol{\eta}}.
    \end{align*}
    Averaging over $\boldsymbol{\eta}$, a $P_0$ update therefore contributes
    $2\mathcal E_0(f)$, while an update at $v$ contributes
    $2\mathcal E_v(f)$. The sequence defining $Q$ contains
    three $P_0$ updates and two updates at each separator
    vertex. Taking expectations in the Cauchy--Schwarz bound
    consequently yields
    \begin{align*}
        \langle f,(I-Q)f\rangle_\nu
          &\leq\frac{2k+3}{2}\sum_{j=1}^{2k+3}
                  \Expected{(f(X_j)-f(X_{j-1}))^2}\\
          &=\frac{2k+3}{2}\left(6\mathcal E_0(f)
                       +4\sum_{v\in\mathcal C}\mathcal E_v(f)\right)\\
          &=(2k+3)\left(3\mathcal E_0(f)
                       +2\sum_{v\in\mathcal C}\mathcal E_v(f)\right).
    \end{align*}
    Combining this with the gap from Step 2 gives
    \[
        \epsilon_{\mathcal C}\Varsub{\nu}{f}
          \leq\langle f,(I-Q)f\rangle_\nu
          \leq3(2k+3)\left(\mathcal E_0(f)
                           +\sum_{v\in\mathcal C}\mathcal E_v(f)\right).
    \]
    Dividing by $\epsilon_{\mathcal C}$ proves
    \eqref{eq:separatorpoincare}, with precisely the stated
    value of $K_{\mathcal C}$.

    \paragraph{Conclusion: recursion over the child components}
    Conditional independence of the children and variance
    tensorization imply
    \[
        \mathcal E_0(f)
           \leq\sum_i\Expectedsub{\nu}{\Varsub{T_i}{f}}
           \leq\max_i F_i
                \sum_{v\in T\setminus\mathcal C}\mathcal E_v(f).
    \]
    Here $\Varsub{T_i}{f}$ denotes the conditional variance given
    all other variable colors. This proves \eqref{eq:separatorrecurrence}.
    If there are no children, take $P_0=I$ and $\mathcal E_0=0$.
\end{proof}

\subsubsection{Iteration along a path of separators}

The recurrence in \Cref{lem:separatorcomparison} allows to take the
maximum relaxation time over all children, rather than their sum. Its total cost is
therefore measured along one root-to-leaf path of separators,
not over the whole block. This is why a diameter bound, together
with control of connected hub sets, suffices.

Let $\mathcal T_B$ be the rooted tree of separators from
\Cref{lem:comparisonseparators}. For a root-to-leaf path
$\mathcal Q$ in $\mathcal T_B$, write $U(\mathcal Q)$ for the
union of its separators, and put
\[
    J(\mathcal Q)=
        \sum_{\mathcal C\in\mathcal Q}
           \left(|\mathcal C|\log q+
                    \frac{r_{\mathcal H}(\mathcal C)}q\right).
\]

\begin{lemma}[Cost and geometry of a separator path]
\label{lem:comparisonseparatorpaths}
    Let $B$ be a tree or unicyclic block with an admissible
    boundary condition, and put $b_B=\operatorname{diam}(G[B])$.
    Uniformly over feasible exterior colors,
\begin{equation}
        \log\tau_{\mathrm{rel}}(B,\taubf)
            \leq28\max_{\mathcal Q}J(\mathcal Q).
        \label{eq:separatorpathcomparison}
\end{equation}
    For each such path, $G[U(\mathcal Q)]$ is connected and
    is a tree or unicyclic, and
\begin{equation}
    |U(\mathcal Q)\setminus\mathcal H|\leq3b_B+2,\qquad
    \sum_{\mathcal C\in\mathcal Q}r_{\mathcal H}(\mathcal C)
        \leq\sum_{u\in U(\mathcal Q)\cap\mathcal H}\deg_G(u).
    \label{eq:separatorpathshape}
\end{equation}
\end{lemma}
\begin{proof}
    Writing $k=|\mathcal C|\geq1$ and using $\ell\leq k$,
    $q\geq7$, and $3(2k+3)\leq15^k$, we obtain
    \[
        \log K_{\mathcal C}
           \leq k\log15+k\log q+4k+
                         \frac{28r_{\mathcal H}(\mathcal C)}q
           \leq28\left(k\log q+
                         \frac{r_{\mathcal H}(\mathcal C)}q\right).
    \]
    Iterate \eqref{eq:separatorrecurrence}, starting at the
    leaves. The maximum-child recurrence bounds the total
    cost by the product of $K_{\mathcal C}$ along one
    root-to-leaf path. Taking logarithms proves
    \eqref{eq:separatorpathcomparison}.

    Consecutive separators are joined by an edge, so their
    union is connected. It induces a tree or unicyclic
    subgraph because it lies in $B$. For a tree block, the
    low-degree roots of successive separators lie along
    a single simple path in $G[B]$. There are at most
    $b_B+1$ such vertices. For a unicyclic block, the initial
    separator can also contain low-degree cycle vertices.
    The cycle has length at most $2b_B+1$, since distances
    between its vertices are their cycle distances.
    All later low-degree roots lie on one path in a tree
    attached to that cycle, contributing at most $b_B+1$
    more vertices. This proves the first inequality in
    \eqref{eq:separatorpathshape}. For the second, charge
    every counted child edge to its hub endpoint in its
    separator. Separators are disjoint, so no edge is
    charged twice.
\end{proof}

\subsubsection{Uniform random-graph estimates}

The remaining task is to bound the cost of a separator path.
A path with few hubs has few vertices by
\Cref{lem:comparisonseparatorpaths}. A path with many hubs has
few vertices because connected hub-rich sets are rare.
The following estimates make this dichotomy uniform over
every possible partition.

\begin{lemma}[Connected hub sets and degree sums]
\label{lem:comparisonrandomgeometry}
    There are absolute constants $c_3,c_4>0$ such that, for fixed
    $0<\eta<1/2$ and sufficiently large fixed $d$, with probability
    $1-o(1)$ over $G\sim G(n,d/n)$, the following hold for every
    nonempty connected
    vertex set $U$ for which $G[U]$ is a tree or unicyclic:
    \begin{enumerate}[(i)]
        \item If $|U\cap\mathcal H|\geq|U|/2$, then
        \[
            |U|\leq c_3\frac{\log n}{d^{2\eta}}.
        \]
        \item In all cases,
        \[
            \sum_{u\in U}\deg_G(u)
                      \leq c_4(d|U|+\log n).
        \]
    \end{enumerate}
\end{lemma}

The idea is to balance the number of possible sets against the
probability of unusually large exterior degrees. We first count the
possible tree and unicyclic internal patterns; conditional on any
such pattern, the exterior degrees are independent. For (i), the
sparse internal pattern forces a positive fraction of the hubs to
have exterior degree at least $d+\tfrac12d^{1/2+\eta}$. Each such
deviation has probability at most $\exp(-c d^{2\eta})$, so a union
bound rules out large sets with many hubs. For (ii), we instead
bound the total number of exterior edges by an exponential tail
estimate and use the same pattern count. The $\log n$ term allows
this bound to hold simultaneously for every set $U$.

\begin{proof}
    Put $m=|U|$, $p=d/n$, $t_0=d^{1/2+\eta}$, and
    $\xi=t_0^2/d=d^{2\eta}$. For $m\geq2$, fix a vertex set
    $U$ of size $m$. By Cayley's formula, there are
    $m^{m-2}$ labeled trees on $U$. For each such tree, the
    probability that it is exactly $G[U]$ is at most $p^{m-1}$:
    its $m-1$ edges must be present, and the requirement that
    all other internal edges be absent only decreases the
    probability.

    Every unicyclic graph on $U$ can be obtained from a tree
    by adding one edge. There are at most $\binom m2\leq m^2$
    choices for that edge, and each resulting graph has
    probability at most $p^m$. A unicyclic graph may be
    counted more than once, which is harmless for an upper
    bound. Thus the two types of internal patterns contribute
    at most $m^{m-2}p^{m-1}$ and $m^{m-2}m^2p^m$,
    respectively. Let $N_m$ count the sets $U$ of size $m$
    for which $G[U]$ is a tree or unicyclic. Summing over
    the $\binom nm$ choices of $U$ gives
\begin{equation}
        \Expected{N_m}
          \leq\binom nm m^{m-2}p^{m-1}(1+m^2p)
                         \leq n^3(c_5d)^m
        \label{eq:comparisonpatterncount}
    \end{equation}
    for an absolute constant $c_5$. To see the last inequality,
    substitute $p=d/n$ and use
    $\binom nm\leq(\mathrm e n/m)^m$ to obtain
    \begin{align*}
        \binom nm m^{m-2}p^{m-1}(1+m^2p)
          &\leq(\mathrm e d)^m
                       \left(\frac{n}{dm^2}+1\right)\\
          &\leq2n(\mathrm e d)^m
           \leq n^3(\mathrm e d)^m,
    \end{align*}
    where $d\geq1$ and $n\geq2$ were used. Thus we may take
    $c_5=\mathrm e$. For $m=1$, $N_1=n$, so the same final
    bound holds directly.
    Conditional on an exact internal pattern, the edges from
    $U$ to its complement remain independent.

    Since $e(G[U])\leq m$, fewer than $4m/t_0$ vertices can have
    internal degree greater than $t_0/2$. If $t_0\geq32$ and
    at least $m/2$ vertices are hubs, at least $m/4$ vertices
    consequently have exterior degree at least $d+t_0/2$.
    The exterior degrees of distinct vertices of $U$ depend
    on disjoint edge indicators and are independent
    $\operatorname{Bin}(n-m,p)$ variables of mean at most $d$.
    To bound their tails, take $\theta=t_0/(4d)\leq1$.
    For such a variable $Y$, exponential Markov and
    $\mathrm e^\theta-1-\theta\leq\theta^2$ give
    \[
        \Prlr{Y\geq d+t_0/2}
          \leq\exp\left(d(\mathrm e^\theta-1-\theta)
                              -\theta t_0/2\right)
          \leq\exp(-\xi/16).
    \]
    If at least $m/2$ vertices are hubs, there is therefore a
    witnessing subset of $\lceil m/4\rceil$ vertices whose
    exterior degrees all exceed $d+t_0/2$. For any specified
    subset, independence bounds its probability by
    $\exp(-\xi\lceil m/4\rceil/16)\leq\exp(-\xi m/64)$.
    There are at most $2^m$ choices of this subset. Hence,
    conditional on any exact tree or unicyclic internal
    pattern, the probability of hub density at least one half
    is at most $2^m\exp(-\xi m/64)$.

    Let $\mathcal F$ be the event that (i) fails, and put
    $m_0=\lfloor c_3\log n/\xi\rfloor+1$. Failure requires
    a set of some size $m\geq m_0$. The preceding conditional
    estimate is uniform over the exact internal patterns, so
    summing it against their probabilities and using
    \eqref{eq:comparisonpatterncount} gives
    \begin{align*}
        \Prlr{\mathcal F}
          &\leq n^3\sum_{m=m_0}^{n}
                         (c_5d)^m2^m\exp(-\xi m/64)\\
          &=n^3\sum_{m=m_0}^{n}
                 \exp\bigl(m\log(2c_5d)-\xi m/64\bigr).
    \end{align*}
    Take $d$ large enough that
    $\xi\geq128\log(2c_5d)$. Each summand is then at most
    $\exp(-\xi m/128)$. Extending the sum to infinity and
    evaluating the geometric series yields
    \begin{align*}
        \Prlr{\mathcal F}
          &\leq n^3\frac{\exp(-\xi m_0/128)}
                             {1-\exp(-\xi/128)}\\
          &\leq\frac{n^{3-c_3/128}}{1-\exp(-\xi/128)},
    \end{align*}
    since $m_0>c_3\log n/\xi$. In particular, choosing
    $c_3=1024$ makes the exponent $3-c_3/128=-5$.
    The denominator is at least $1/2$ for sufficiently large
    $d$, so the failure probability in (i) is at most
    $2n^{-5}=o(1)$.

    For (ii), conditional on the internal pattern, let
    $X=e(U,V\setminus U)$. Then
    $X\sim\operatorname{Bin}(m(n-m),p)$ and
    $\Expected{\mathrm e^X}\leq\exp((\mathrm e-1)dm)$.
    For any fixed sufficiently large absolute $c_6$,
    exponential Markov gives 
    \[
        \Prlr{X>c_6(dm+\log n)}
            \leq n^{-c_6}
                    \exp\bigl(-(c_6-\mathrm e+1)dm\bigr).
    \]
    A union bound using \eqref{eq:comparisonpatterncount}
    bounds the probability that this fails for any such $U$ by
    \[
        n^{3-c_6}\sum_{m=1}^n
                (c_5d)^m
                    \exp\bigl(-(c_6-\mathrm e+1)dm\bigr)=o(1).
    \]
    Finally,
    $\sum_{u\in U}\deg_G(u)=2e(G[U])+X\leq2m+X$.
    Increasing an absolute constant gives (ii).
\end{proof}

\subsubsection{Proof of \Cref{thm:withinblockcomparison}}

Combining all the previous estimates, we can prove \Cref{thm:withinblockcomparison}.

\begin{proof}[Proof of \Cref{thm:withinblockcomparison}]
    Work on the event in \Cref{lem:comparisonrandomgeometry}.
    For a separator path and its union $U$, put $m=|U|$.
    If at least half of $U$ consists of hubs, part (i) gives
    $m\leq c_3\log n/d^{2\eta}$. Otherwise,
    \eqref{eq:separatorpathshape} gives
    $m<2|U\setminus\mathcal H|\leq6b_B+4$.
    In either case,
    \[
        m\leq\max\left\{6b_B+4,\,
                             c_3\frac{\log n}{d^{2\eta}}\right\}.
    \]
    The separators are disjoint, and by part (ii) and
    \eqref{eq:separatorpathshape},
    \[
        J(\mathcal Q)
            \leq m\log q+\frac{c_4}{q}(dm+\log n)
            \leq c_7\left(m\log q+\frac{\log n}{q}\right),
    \]
    where $q>d$ and $q\geq7$ were used. Together with
    \eqref{eq:separatorpathcomparison}, this proves
    \eqref{eq:localrelaxationgeneral}, uniformly over all
    feasible exterior colors.
    If $b_B\leq D\log n/d^{2\eta}$, then
    $d^{2\eta}<d<q$, so $\log n/q$ is absorbed by
    $(\log n/d^{2\eta})\log q$.
    Taking $n$ sufficiently large depending on $d$ and $D$,
    but not on $q$, proves \eqref{eq:withinblockloss}.

    Finally, \Cref{thm:validblockpartition} supplies the stated
    diameter bound with $D=16(K+1)$.
\end{proof}

\subsection{Comparison with the full Glauber dynamics}

Having established a within-block relaxation estimate, we can bound
the spectral gap of Glauber dynamics by comparing Dirichlet forms.
Similar block-to-single-site comparisons for colorings of sparse
random graphs are used by
\citet[Proposition~57 in the full version]{efthymiou2018sampling}
and \citet[Lemma~5.2 and Section~5.5]{li2026rapid}.

Before using the block gap as a Poincare inequality, we check
that the coloring space is nonempty and that no nonconstant
function lies in the kernel. This avoids an implicit
irreducibility assumption.

\begin{lemma}[Nonemptiness and irreducibility]
\label{lem:comparisonirreducibility}
    Suppose $G$ has a valid block partition and
    $q\geq\lfloor\lowdegree\rfloor+2$, with $\lowdegree\geq2$.
    Then $G$ has a proper $q$-coloring. Its single-site
    heat-bath dynamics and its exact block dynamics are
    irreducible, and both generators have only constants
    in their kernel.
\end{lemma}
\begin{proof}
    Put $k=\lfloor\lowdegree\rfloor$. Delete the singleton
    vertices first, in any order. At deletion, each has
    degree at most $k$, because it is a global breakpoint.
    The remaining components are trees or unicyclic, so
    their vertices can then be deleted successively with
    degree at most two. This gives a $k$-degeneracy ordering.
    Greedy coloring in the reverse order proves nonemptiness.

    We prove connectivity of the recoloring graph of any
    $k$-degenerate graph by induction on its number of vertices.
    Remove a vertex $v$ of degree at most $k$, and connect the
    restrictions of two proper colorings by an inductively
    obtained sequence on $G-v$. Lift each move to $G$.
    If a neighbor is about to take the current color of $v$,
    first recolor $v$ to a color avoiding all its neighbors'
    current colors and the neighbor's proposed color.
    At most $k+1$ colors are forbidden, so a spare color
    exists. After lifting the sequence, recolor $v$ to its
    final color. This is legal because its neighbors already
    have their final colors.

    Every legal single-site recoloring has positive heat-bath
    probability. It also has positive probability under an
    exact update of the block containing that vertex, by
    leaving the other block colors unchanged. Thus both
    chains are irreducible. The kernel assertions follow
    from their Dirichlet forms.
\end{proof}

The local Poincare inequalities can now be averaged over the
exterior colors. No factor involving the number or sizes of the
blocks is lost in this step.

Let $\mu$ be the uniform measure on the nonempty space of proper
$q$-colorings of $G$. Extend the single-site notation to every
$v\in V$ by
\[
    L_v f=f-\Expectednop{f\mid\sigmabf_{V\setminus\{v\}}},
    \qquad L_{\mathrm{Gl}}=\sum_{v\in V}L_v.
\]
As in the preceding section, $L$ denotes the exact block generator
with rate one for each singleton and each block.
Here $\Varsub{B}{f}$ denotes the variance conditional on all colors outside $B$.

\begin{theorem}[Block-to-single-site comparison]
\label{prop:blockglaubercomparison}
    Suppose that the block dynamics satisfies
    \[
        \Varsub{\mu}{f}
                       \leq\gamma^{-1}\langle f,Lf\rangle_\mu
        \qquad\text{for every }f,
    \]
    and the uniform within-block relaxation bound is $T_*$.
    Then
    \[
        \operatorname{gap}(L_{\mathrm{Gl}})\geq\frac{\gamma}{T_*}.
    \]
    For the discrete-time chain
    $P_{\mathrm{Gl}}=n^{-1}\sum_{v\in V}(I-L_v)$,
\begin{equation}
        \operatorname{gap}(I-P_{\mathrm{Gl}})
                  \geq\frac{\gamma}{nT_*}.
        \label{eq:globalglaubergap}
    \end{equation}
\end{theorem}
\begin{proof}
    For a heat-bath update, the Dirichlet form is the expected
    variance conditional on the colors that are not updated. Hence
    \begin{align*}
        \langle f,Lf\rangle_\mu
           &=\sum_{s\in S}\Expectedsub{\mu}{\Varsub{s}{f}}
               +\sum_{B\in\mathcal B}
                    \Expectedsub{\mu}{\Varsub{B}{f}},\\
        \langle f,L_{\mathrm{Gl}}f\rangle_\mu
           &=\sum_{v\in V}\Expectedsub{\mu}{\Varsub{v}{f}}.
    \end{align*}

    \paragraph{Within-block Poincare inequalities}
    Fix $B\in\mathcal B$ and a feasible exterior coloring $\taubf$.
    The conditional law of $\sigmabf_B$ under $\mu$, given
    $\sigmabf_{V\setminus B}=\taubf$, is $\nu_{B,\taubf}$.
    Write $f^{\taubf}(\sigmabf_B)=f(\sigmabf_B,\taubf)$ for the
    restriction of $f$ to this conditional state space.
    By \Cref{def:withinblockrelaxation} and the uniform bound
    $\tau_{\mathrm{rel}}(B,\taubf)\leq T_*$,
    \[
        \Varsub{\nu_{B,\taubf}}{f^{\taubf}}
          \leq T_*\sum_{v\in B}
              \Expectedsub{\nu_{B,\taubf}}{\Varsub{v}{f^{\taubf}}}.
    \]
    Here the variance at $v$ conditions on the other colors in $B$
    as well as the fixed exterior $\taubf$, and therefore conditions
    on exactly the same colors as $\Varsub{v}{f}$ in the full chain.

    \paragraph{Averaging and summing over the partition}
    Multiply this inequality by
    $\mu(\sigmabf_{V\setminus B}=\taubf)$ and sum over feasible
    $\taubf$. The tower property gives
    \begin{align*}
        \Expectedsub{\mu}{\Varsub{B}{f}}
          &=\sum_{\taubf}\mu(\sigmabf_{V\setminus B}=\taubf)
                         \Varsub{\nu_{B,\taubf}}{f^{\taubf}}\\
          &\leq T_*\sum_{v\in B}
                         \Expectedsub{\mu}{\Varsub{v}{f}}.
    \end{align*}
    Apply this to each block term in the Dirichlet form of $L$.
    Since the blocks together with the singletons partition $V$,
    every vertex is counted once. Using $T_*\geq1$ also for the
    singleton terms, we obtain
    \begin{align*}
        \langle f,Lf\rangle_\mu
          &\leq\sum_{s\in S}\Expectedsub{\mu}{\Varsub{s}{f}}
              +T_*\sum_{B\in\mathcal B}\sum_{v\in B}
                         \Expectedsub{\mu}{\Varsub{v}{f}}\\
          &\leq T_*\sum_{v\in V}\Expectedsub{\mu}{\Varsub{v}{f}}
           =T_*\langle f,L_{\mathrm{Gl}}f\rangle_\mu.
    \end{align*}

    \paragraph{Applying the block spectral gap}
    The assumed block Poincare inequality is the spectral gap
    bound $\operatorname{gap}(L)\geq\gamma$ applied to the centered
    function $f-\Expectedsub{\mu}{f}$, whose squared norm is
    $\Varsub{\mu}{f}$. Combining it with the comparison above yields
    \[
        \Varsub{\mu}{f}
          \leq\gamma^{-1}\langle f,Lf\rangle_\mu
          \leq\frac{T_*}{\gamma}
                     \langle f,L_{\mathrm{Gl}}f\rangle_\mu
          =\frac{T_*}{\gamma}\sum_{v\in V}
                     \Expectedsub{\mu}{\Varsub{v}{f}}.
    \]
    This is the Poincare inequality for continuous-time Glauber
    dynamics with each vertex updated at rate one, so its spectral
    gap is at least $\gamma/T_*$. In particular, zero Glauber
    Dirichlet form forces zero variance, and hence its kernel
    consists only of constants. Finally,
    $I-P_{\mathrm{Gl}}=L_{\mathrm{Gl}}/n$ gives the discrete-time gap.
\end{proof}

\subsection{Proof of the main theorem}

We now combine the established estimates. Choosing
$\eta=\varepsilon/2$ makes the exponent in the within-block
comparison cost of order $\log q/d^\varepsilon$, as required
by \Cref{thm:main}.

\begin{proof}[Proof of \Cref{thm:main}]
    Fix $0<\varepsilon\leq1/2$ and choose
    \[
        \eta=\frac{\varepsilon}{2}.
    \]
    Choose $A\geq1$ and $C>0$ sufficiently large as required
    by \Cref{thm:spectralgapblockcoloring}, and then take
    $d$ sufficiently large. Intersect its high-probability
    event with those of
    \Cref{thm:validblockpartition,thm:withinblockcomparison}.
    The partition exists with a diameter constant
    $D=D_\varepsilon$ independent of $d,n,q$, and its
    within-block relaxation bound is
    \[
        T_*\leq n^{\zeta_d},
        \qquad
        \zeta_d=c_1(D_\varepsilon+1)
                       \frac{\log q}{d^\varepsilon}.
    \]
    By \Cref{lem:comparisonirreducibility}, the block gap
    estimate is a Poincare inequality with
    $\gamma\geq c_\varepsilon d^{\varepsilon-1/2}$.
    By \Cref{prop:blockglaubercomparison}, the discrete-time
    chain $P_{\mathrm{GD}}=P_{\mathrm{Gl}}$ therefore satisfies
    \[
        g\coloneqq\operatorname{gap}(I-P_{\mathrm{GD}})
            \geq\frac{\gamma}{nT_*}
            \geq c_\varepsilon d^{\varepsilon-1/2}
                            n^{-1-\zeta_d}.
    \]
    The heat-bath kernel $P_{\mathrm{GD}}$ is an average of
    orthogonal projections, so its spectrum is nonnegative.
    Since $\mu_{q,G}$ is uniform on at most $q^n$ proper
    colorings, its smallest atom satisfies
    $\mu_{\min}\geq q^{-n}$. The standard spectral-gap
    bound on the mixing time now gives, for every $0<\delta<1$,
    \[
        \tmix(\delta)\leq
           \left\lceil\frac1g
                  \left(\frac12\log\frac1{\mu_{\min}}
                              +\log(1/\delta)\right)\right\rceil
           \leq C_{d,\varepsilon}\,
                    n^{2+\zeta_d}\log q\,
                         \bigl(1+\log(1/\delta)\bigr).
    \]
    For $\log n\geq d^\varepsilon$, absorb $\log q$ using
    $\log q\leq n^{\log q/d^\varepsilon}$. Thus, with
    $C_\varepsilon=c_1(D_\varepsilon+1)+1$, we obtain
    \[
        \tmix(\delta)\leq C_{d,\varepsilon}\,
                    n^{2+C_\varepsilon\log q/d^\varepsilon}
                         \bigl(1+\log(1/\delta)\bigr).
    \]
    The event has probability $1-o(1)$, and both the event
    and the lower threshold on $n$ are independent of $q$ and $\delta$.
    The constant $C_\varepsilon$ depends only on
    $\varepsilon$, which proves \Cref{thm:main}.
\end{proof}

%% file: references.bib
@article{le2017concentration,
  author={Le, Can M. and Levina, Elizaveta and Vershynin, Roman},
  title={Concentration and regularization of random graphs},
  journal={Random Structures \& Algorithms},
  volume={51},
  number={3},
  pages={538--561},
  year={2017},
  doi={10.1002/rsa.20713}
}

@article{feige2005spectral,
  author={Feige, Uriel and Ofek, Eran},
  title={Spectral techniques applied to sparse random graphs},
  journal={Random Structures \& Algorithms},
  volume={27},
  number={2},
  pages={251--275},
  year={2005},
  doi={10.1002/rsa.20089}
}

@incollection{bakry2006diffusions,
  title={Diffusions hypercontractives},
  author={Bakry, Dominique and {\'E}mery, Michel},
  booktitle={S{\'e}minaire de Probabilit{\'e}s XIX 1983/84: Proceedings},
  pages={177--206},
  year={2006},
  publisher={Springer}
}

@book{bakry2013analysis,
  title={Analysis and Geometry of {M}arkov Diffusion Operators},
  author={Bakry, D. and Gentil, I. and Ledoux, M.},
  isbn={9783319002279},
  lccn={2013952461},
  series={Grundlehren der mathematischen Wissenschaften},
  url={https://books.google.de/books?id=gU3ABAAAQBAJ},
  year={2013},
  publisher={Springer International Publishing}
}

@incollection {chen2025rapid,
    AUTHOR = {Chen, Xiaoyu and Chen, Zejia and Chen, Zongchen and Yin,
              Yitong and Zhang, Xinyuan},
     TITLE = {Rapid mixing on random regular graphs beyond uniqueness},
 BOOKTITLE = {2025 {IEEE} 66th {A}nnual {S}ymposium on {F}oundations of
              {C}omputer {S}cience---{FOCS} 2025},
     PAGES = {2170--2193},
 PUBLISHER = {IEEE Comput. Soc. Press, Los Alamitos, CA},
      YEAR = {2025},
      ISBN = {979-8-3315-7132-0},
   MRCLASS = {68Q87},
  MRNUMBER = {5048082},
       DOI = {10.1109/FOCS63196.2025.00115},
       URL = {https://doi-org.turing.library.northwestern.edu/10.1109/FOCS63196.2025.00115},
}

@article{goebel2026simple,
  title={A simple proof of rapid mixing on random regular graphs beyond uniqueness},
  author={G{\"o}bel, Andreas and Jenssen, Matthew and Michelen, Marcus and Pappik, Marcus and Perkins, Will and Schiller, Leon},
  journal={arXiv preprint arXiv:2606.27545},
  year={2026}
}

@article{haqi2026sampling,
  title={Sampling Line-Graph Colorings with Constant Extra Colors},
  author={Haqi, Alireza},
  journal={arXiv preprint arXiv:2609.27440},
  year={2026}
}

@article{chen2026spectral,
  title={A Spectral Local-to-Global Principle for Spin Systems on Graphs with Girth At Least Five},
  author={Chen, Xiaoyu and Liu, Kuikui},
  journal={arXiv preprint arXiv:2608.25491},
  year={2026}
}

@article{chen2026spectral2,
  title={Spectral Gap of Down-Up Walks via Trickle-Down: A Simplified and Sharpened Analysis},
  author={Chen, Xiaoyu and Liu, Kuikui},
  journal={arXiv preprint arXiv:2609.19514},
  year={2026}
}

@article{boban2026rank,
  title={Rank-1-perturbed trickledown theorems: Mixing time of Glauber dynamics for the Sherrington-Kirkpatrick model up to beta less than 1/2 + epsilon},
  author={Boban, Mathews and Li, Anqi and Gharan, Shayan Oveis},
  journal={arXiv preprint arXiv:2609.13138},
  year={2026}
}

@article{wang2026optimal,
  title={Optimal Mixing of Glauber Dynamics for the Sherrington-Kirkpatrick Model at beta less than 1/2},
  author={Wang, Sihan},
  journal={arXiv preprint arXiv:2608.22159},
  year={2026}
}

@inproceedings{efthymiou2018sampling,
  title={Sampling random colorings of sparse random graphs},
  author={Efthymiou, Charilaos and Hayes, Thomas P and {\v{S}}tefankovi{\v{c}}, Daniel and Vigoda, Eric},
  booktitle={Proceedings of the Twenty-Ninth Annual ACM-SIAM Symposium on Discrete Algorithms},
  pages={1759--1771},
  year={2018},
  organization={SIAM}
}

@article{chen2026high,
  title={High-Dimensional Ultra-Log-Concave Distributions},
  author={Chen, Zongchen and Wang, Sihan},
  journal={arXiv preprint arXiv:2609.23994},
  year={2026}
}

@article{guo2026approximating,
  title={Approximating spin systems on planar graphs},
  author={Guo, Heng and Zhang, Xinyuan},
  journal={arXiv preprint arXiv:2608.06172},
  year={2026}
}

@article{miao2026sampling,
  title={Sampling Matchings in Near-linear Time},
  author={Miao, Tianshun and Yin, Yitong},
  journal={arXiv preprint arXiv:2609.21936},
  year={2026}
}

@article{li2026rapid,
  title={Rapid Mixing for Sampling Colorings of Random Graphs up to Tree Uniqueness Threshold},
  author={Li, Zhidan and Yang, Kuan},
  journal={arXiv preprint arXiv:2610.25391},
  year={2026}
}

@article{benaych2019largest,
  title={Largest eigenvalues of sparse inhomogeneous {Erd{\H{o}}s--R{\'e}nyi} graphs},
  author={Benaych-Georges, Florent and Bordenave, Charles and Knowles, Antti},
  journal={The Annals of Probability},
  volume={47},
  number={3},
  pages={1653--1676},
  year={2019},
  doi={10.1214/18-AOP1293},
  url={https://doi.org/10.1214/18-AOP1293}
}

@article{mossel2010gibbs,
  title={Gibbs rapidly samples colorings of {$G(n,d/n)$}},
  author={Mossel, Elchanan and Sly, Allan},
  journal={Probability Theory and Related Fields},
  volume={148},
  number={1--2},
  pages={37--69},
  year={2010},
  doi={10.1007/s00440-009-0222-x}
}

@inproceedings{efthymiou2014mcmc,
  title={{MCMC} sampling colourings and independent sets of {$G(n,d/n)$} near uniqueness threshold},
  author={Efthymiou, Charilaos},
  booktitle={Proceedings of the Twenty-Fifth Annual ACM-SIAM Symposium on Discrete Algorithms},
  pages={305--316},
  year={2014},
  publisher={SIAM},
  doi={10.1137/1.9781611973402.22}
}

@article{efthymiou2016simple,
  title={A simple algorithm for sampling colorings of {$G(n,d/n)$} up to the {G}ibbs uniqueness threshold},
  author={Efthymiou, Charilaos},
  journal={SIAM Journal on Computing},
  volume={45},
  number={6},
  pages={2087--2116},
  year={2016},
  doi={10.1137/140977643}
}

@article{chen2023strong,
  title={Strong spatial mixing for colorings on trees and its algorithmic applications},
  author={Chen, Zongchen and Liu, Kuikui and Mani, Nitya and Moitra, Ankur},
  journal={arXiv preprint arXiv:2304.01954},
  year={2023},
  url={https://arxiv.org/abs/2304.01954}
}

@article{jerrum1995simple,
  title={A very simple algorithm for estimating the number of {$k$}-colorings of a low-degree graph},
  author={Jerrum, Mark},
  journal={Random Structures \& Algorithms},
  volume={7},
  number={2},
  pages={157--165},
  year={1995},
  doi={10.1002/rsa.3240070205}
}

@article{vigoda2000improved,
  title={Improved bounds for sampling colorings},
  author={Vigoda, Eric},
  journal={Journal of Mathematical Physics},
  volume={41},
  number={3},
  pages={1555--1569},
  year={2000},
  doi={10.1063/1.533196}
}

@inproceedings{carlson2025flip,
  title={Flip dynamics for sampling colorings: Improving {$(11/6-\varepsilon)$} using a simple metric},
  author={Carlson, Charlie and Vigoda, Eric},
  booktitle={Proceedings of the 2025 Annual ACM-SIAM Symposium on Discrete Algorithms},
  pages={2194--2212},
  year={2025},
  publisher={SIAM},
  doi={10.1137/1.9781611978322.71},
  url={https://arxiv.org/abs/2407.04870}
}

@article{metropolis1953equation,
  title={Equation of state calculations by fast computing machines},
  author={Metropolis, Nicholas and Rosenbluth, Arianna W. and Rosenbluth, Marshall N. and Teller, Augusta H. and Teller, Edward},
  journal={The Journal of Chemical Physics},
  volume={21},
  number={6},
  pages={1087--1092},
  year={1953},
  doi={10.1063/1.1699114}
}

@article{zdeborova2007phase,
  title={Phase transitions in the coloring of random graphs},
  author={Zdeborov{\'a}, Lenka and Krzakala, Florent},
  journal={Physical Review E},
  volume={76},
  number={3},
  pages={031131},
  year={2007},
  doi={10.1103/PhysRevE.76.031131}
}

@article{krzakala2007gibbs,
  title={{G}ibbs states and the set of solutions of random constraint satisfaction problems},
  author={Krzakala, Florent and Montanari, Andrea and Ricci-Tersenghi, Federico and Semerjian, Guilhem and Zdeborov{\'a}, Lenka},
  journal={Proceedings of the National Academy of Sciences},
  volume={104},
  number={25},
  pages={10318--10323},
  year={2007},
  doi={10.1073/pnas.0703685104}
}
